\documentclass[12pt]{amsart}
\usepackage{amsbsy}
\usepackage{graphicx,epsfig,subfigure}
\usepackage[scanall]{psfrag}
\begin{document}

\newtheorem{theorem}{Theorem}
\newtheorem{proposition}{Proposition}
\newtheorem{lemma}{Lemma}
\newtheorem{corollary}{Corollary}
\newtheorem{definition}{Definition}
\newtheorem{remark}{Remark}
\newcommand{\tex}{\textstyle}
\numberwithin{equation}{section} \numberwithin{theorem}{section}
\numberwithin{proposition}{section} \numberwithin{lemma}{section}
\numberwithin{corollary}{section}
\numberwithin{definition}{section} \numberwithin{remark}{section}
\newcommand{\ren}{\mathbb{R}^N}
\newcommand{\re}{\mathbb{R}}
\newcommand{\n}{\nabla}
\newcommand{\iy}{\infty}
\newcommand{\pa}{\partial}
\newcommand{\fp}{\noindent}
\newcommand{\ms}{\medskip\vskip-.1cm}
\newcommand{\mpb}{\medskip}
\newcommand{\AAA}{{\bf A}}
\newcommand{\BB}{{\bf B}}
\newcommand{\CC}{{\bf C}}
\newcommand{\DD}{{\bf D}}
\newcommand{\EE}{{\bf E}}
\newcommand{\FF}{{\bf F}}
\newcommand{\GG}{{\bf G}}
\newcommand{\oo}{{\mathbf \omega}}
\newcommand{\Am}{{\bf A}_{2m}}
\newcommand{\CCC}{{\mathbf  C}}
\newcommand{\II}{{\mathrm{Im}}\,}
\newcommand{\RR}{{\mathrm{Re}}\,}
\newcommand{\eee}{{\mathrm  e}}
\newcommand{\LL}{L^2_\rho(\ren)}
\newcommand{\LLL}{L^2_{\rho^*}(\ren)}
\renewcommand{\a}{\alpha}
\renewcommand{\b}{\beta}
\newcommand{\g}{\gamma}
\newcommand{\G}{\Gamma}
\renewcommand{\d}{\delta}
\newcommand{\D}{\Delta}
\newcommand{\e}{\varepsilon}
\newcommand{\var}{\varphi}
\newcommand{\lll}{\l}
\renewcommand{\l}{\lambda}
\renewcommand{\o}{\omega}
\renewcommand{\O}{\Omega}
\newcommand{\s}{\sigma}
\renewcommand{\t}{\tau}
\renewcommand{\th}{\theta}
\newcommand{\z}{\zeta}
\newcommand{\wx}{\widetilde x}
\newcommand{\wt}{\widetilde t}
\newcommand{\noi}{\noindent}
\newcommand{\uu}{{\bf u}}
\newcommand{\xx}{{\bf x}}
\newcommand{\yy}{{\bf y}}
\newcommand{\zz}{{\bf z}}
\newcommand{\aaa}{{\bf a}}
\newcommand{\cc}{{\bf c}}
\newcommand{\jj}{{\bf j}}
\newcommand{\ggg}{{\bf g}}
\newcommand{\UU}{{\bf U}}
\newcommand{\YY}{{\bf Y}}
\newcommand{\HH}{{\bf H}}
\newcommand{\GGG}{{\bf G}}
\newcommand{\VV}{{\bf V}}
\newcommand{\ww}{{\bf w}}
\newcommand{\vv}{{\bf v}}
\newcommand{\hh}{{\bf h}}
\newcommand{\di}{{\rm div}\,}
\newcommand{\ii}{{\rm i}\,}
\newcommand{\inA}{\quad \mbox{in} \quad \ren \times \re_+}
\newcommand{\inB}{\quad \mbox{in} \quad}
\newcommand{\inC}{\quad \mbox{in} \quad \re \times \re_+}
\newcommand{\inD}{\quad \mbox{in} \quad \re}
\newcommand{\forA}{\quad \mbox{for} \quad}
\newcommand{\whereA}{,\quad \mbox{where} \quad}
\newcommand{\asA}{\quad \mbox{as} \quad}
\newcommand{\andA}{\quad \mbox{and} \quad}
\newcommand{\withA}{,\quad \mbox{with} \quad}
\newcommand{\orA}{,\quad \mbox{or} \quad}
\newcommand{\atA}{\quad \mbox{at} \quad}
\newcommand{\onA}{\quad \mbox{on} \quad}
\newcommand{\ef}{\eqref}
\newcommand{\ssk}{\smallskip}
\newcommand{\LongA}{\quad \Longrightarrow \quad}
\def\com#1{\fbox{\parbox{6in}{\texttt{#1}}}}
\def\N{{\mathbb N}}
\def\A{{\cal A}}
\newcommand{\de}{\,d}
\newcommand{\eps}{\varepsilon}
\newcommand{\be}{\begin{equation}}
\newcommand{\ee}{\end{equation}}
\newcommand{\spt}{{\mbox spt}}
\newcommand{\ind}{{\mbox ind}}
\newcommand{\supp}{{\mbox supp}}
\newcommand{\dip}{\displaystyle}
\newcommand{\prt}{\partial}
\renewcommand{\theequation}{\thesection.\arabic{equation}}
\renewcommand{\baselinestretch}{1.1}
\newcommand{\Dm}{(-\D)^m}

\title
{\bf
Eigenfunctions and very singular similarity solutions of odd-order
nonlinear dispersion PDEs}



\author {R.S.~Fernandes and V.A.~Galaktionov}

\address{Department of Mathematical Sciences, University of Bath,
 Bath BA2 7AY, UK}
\email{vag@maths.bath.ac.uk}

\address{Department of Mathematical Sciences, University of Bath,
 Bath BA2 7AY, UK}
\email{rsf21@maths.bath.ac.uk}



\keywords{Odd-order nonlinear dispersion equations,
 fundamental solutions, nonlinear eigenfunctions,
 very singular  self-similar solutions, branching, bifurcation}

 \subjclass{35K55, 35K40, 35K65}
\date{\today}

\pagenumbering{arabic}

\begin{abstract}


Asymptotic properties of solutions of the nonlinear dispersion
equations
 \be
  \label{jjj1}
 u_t=(|u|^n u)_{xxx}
  \quad \mbox{and} \quad
 u_t=(|u|^n u)_{xxx} -|u|^{p-1} u \quad \mbox{in}
\quad \re \times \re_+,
  \ee
 where $n>0$ and $p>n+1$ are fixed exponents,  and their higher-odd-order analogies
  are studied.
  The global in time similarity solutions, which
lead to
  ``nonlinear
 eigenfunctions" of the rescaled ODEs, are constructed.
 The basic mathematical  tools include
 a ``homotopy-deformation" approach, where the limit $n \to 0^+$ in
 the first equation in (\ref{jjj1})
  turns out to be fruitful by reducing the problems to
  the linear one at $n=0$, i.e.,
 $$
 v_t=v_{xxx},
 $$
  for which Hermitian spectral theory was developed
 earlier in \cite{RayGI}, as well as other nonlinear operator and
  numerical methods.
  The nonlinear limit $n \to + \iy$, with the limit PDE
   $$
 ({\rm sign}\, v)_t=v_{xxx},
\quad \mbox{
 in terms of the variable $v=|u|^n u$},
  $$
 admitting three almost explicit nonlinear eigenfunctions,
    has also been described.

  For the
 second equation in (\ref{jjj1}), {\em very singular similarity
 solutions} (VSSs)
 are constructed. A ``nonlinear bifurcation" phenomenon at
 critical values $\{p=p_l(n)\}_{l \ge 0}$ of the absorption exponents
 is  discussed.
 In fact, regardless their wide spread in applications and clear significance,
  such third-order nonlinear dispersion
 PDEs are quite poorly represented in general PDE theory, unlike
 their direct ``neighbouring" second and fourth-order quasilinear
 parabolic equations  (the so-called {\em porous medium equations}, PMEs)
  $$
  u_t=(|u|^nu)_{xx} \,\,\,\mbox{(the PME--2)} \andA u_t=-(|u|^n u)_{xxxx}
  \,\,\,\mbox{(the PME--4)}.
  $$
  The main goal of the
 present paper is to essentially fill such a gap.


\end{abstract}

\maketitle







\section{Introduction: nonlinear dispersion PDEs}
 \label{S1}

\subsection{NDEs: global smooth similarity patterns}

In the present paper, we study  asymptotic properties of {\em
nonlinear dispersion equations} (NDEs) of the following form:
\begin{equation}
 \label{nde33}
u_t = (-1)^{k+1}D^{2k+1}_x(|u|^nu)- |u|^{p-1}u \quad \mbox{in}
\quad \re \times \re_+, \quad k=1,2,...\,,
\end{equation}
where $n >0$ and $p>n+1$ are  fixed exponents.
 This
 continues the study began in
\cite{RayGI} for $n=0$, i.e., for the {\em semilinear} dispersion
equation. A continuous ``homotopy" path $n \to 0^+$ connecting
the present study with that in \cite{RayGI} turns out to be rather
effective.

Let us clarify from the beginning a general role and the main
purpose of the present study in the theory of higher-order NDEs
such as \ef{nde33}:

{\bf (i)}  We will study here a countable family of {\em
sufficiently smooth}
  continuous similarity solutions of VSS (very singular) type, which are defined for all
  $t>0$.

 {\bf (ii)} In addition, (\ref{nde33}), as an equation
  {\em with nonlinear dispersion mechanism}, contains and
  describes
  other key
 singularity phenomena such as a complicated formation of various
{\em shock} and smoother
 {\em rarefaction} waves, which appear from discontinuous data, as
 well as a general {\em nonuniqueness} of ``entropy solutions" after
 shocks.
  We do not touch these difficult, even mathematically obscure, phenomena  and refer to \cite{GPndeII, GPnde, Gal3NDENew} for further details.

 Overall, it may be said that the smooth similarity behaviour and asymptotic patterns studied here
  occur in the NDE \ef{nde33} for large times, when the shock wave influence
 has already been settled down via evolution, and becomes negligible.

\ssk

Concerning local existence and uniqueness theory, including
``smoothing results" for NDEs (i.e., for solutions without shocks
and stronger singularities), see references and results in
\cite{GPndeII}. Since we are dealing with some special exact
similarity solutions of (\ref{nde33}), we do not use any
advanced results of local and/or global regularity and any
shock-entropy theory here, and always tackle continuous solutions.

Various applications of nonlinear dispersion equations are
characterized in detail in \cite{GPnde}, so we will minimize extra
references to previous physical and more formal investigations of
such PDE models, and concentrate on their quite unusual
mathematical aspects of our interest.
 It is worth  mentioning here
 a couple of remarkable  examples of NDEs from integrable
 PDE theory. The first one is the third-order  {\em Harry Dym}
 {\em equation}
  \be
 \label{HD0}
  u_t = u^3 u_{xxx} \, ,
  \ee
  a  most exotic integrable soliton equation; see
\cite[\S~4.7]{GSVR} for survey and references
 therein.  The second one is the fifth-order  {\em
Kawamoto equation} \cite{Kaw85}
 $$
 u_t = u^5 u_{xxxxx}+ 5 \,u^4 u_x u_{xxxx}+ 10 \,u^5 u_{xx}
 u_{xxx},
$$
 which
  has higher-degree algebraic terms.
Shock waves for the pure fifth-order NDEs
 $$
 u_t=u^5 u_{xxxxx} \andA u_t=(|u|^n u)_{xxxxx},
 $$
 are studied in \cite{GalNDE5}.
   Let us discuss further
  necessary aspects of NDEs and their standing in general PDE
  theory.



\subsection{Compactons in NDEs: compactly supported travelling waves}

Consider the higher odd-order nonlinear dispersion equation  with
another divergent lower-order term:
\begin{equation}
\label{GenNonlin} u_t = (-1)^{k+1}D^{2k+1}_x(|u|^nu)+ (|u|^n u)_x.
\end{equation}
Equation \eqref{GenNonlin} is a generalization of the  third-order
{\em Rosenau--Hyman} (RH) {\em equation}
\begin{equation}
 \label{RH11}
u_t = (u^2)_{xxx} + (u^2)_x,
\end{equation}
which models the effect of nonlinear dispersion in the pattern
formation of liquid drops (see \cite{Compactons}).  For $n=k=1$,
\eqref{GenNonlin} is the RH equation in classes of nonnegative
solutions.

It is well known that the RH equation \ef{RH11} possesses explicit
moving compactly supported soliton-type solutions known as {\em
compactons}. Compactons have the same structure as travelling wave
solutions, given by
\begin{equation*}
u_{\rm c}(x,t)= f(z) \whereA  z=x-\lambda t.
\end{equation*}
So looking at compacton solutions for \eqref{GenNonlin}, on
substitution we have that
\begin{equation*}
-\lambda
f^{\prime}=(-1)^{k+1}D^{2k+1}_z(|f|^nf)+(|f|^nf)^{\prime}.
\end{equation*}
After integrating once with zero constant (i.e., zero ``flux"), we
find that $f(z)$ satisfies
\begin{equation}
\label{intnontwlode} -\lambda f=(-1)^{k+1}D^{2k}_z(|f|^nf)+|f|^nf.
\end{equation}
 For $k=1$, the ODE \ef{intnontwlode} can be solved explicitly
 (see the expression \ef{f11} below), while for $k \ge 2$, this is a difficult
 variational problem, which admits various countable families  of
 compactly supported oscillatory solutions $f(z)$ of changing
 sign, \cite{GMPSob}.

\ssk

Whilst we note that compacton solutions may be found for nonlinear
dispersion equations, for non-conservative and non-divergent  NDEs
such as \ef{nde33}, TW solutions may be irrelevant for classes of
bounded solutions of the Cauchy problem.
 Therefore, instead we look to find other
more complicated similarity solutions.

 \subsection{A relation to blow-up in reaction-diffusion theory}

Surprisingly, the NDE \eqref{GenNonlin} is somehow related to the
parabolic even-order equations of reaction-diffusion type:
\begin{equation}
\label{ParaEven} u_t=(-1)^{k+1}D^{2k}_x(|u|^nu)+|u|^nu
 \quad \big(u_t=(u^{n+1})_{xx} + u^{n+1} \,\, \mbox{for} \,\,
 k=1, \,\, u \ge 0\big).
\end{equation}
Equation \eqref{ParaEven} admits blow-up self-similar solutions of
the separate form
\begin{equation}
 \label{SS1}
\mbox{$u(x,t)=(T-t)^{-\frac{1}{n}}f(x)$},
\end{equation}
where $T$ is the finite blow-up time.  This self-similarity
reduces the PDE to the ODE for the similarity profile $f$ easily
obtained by substituting \ef{SS1} to \ef{GenNonlin}:
\begin{equation}
\label{EvenSS}
\mbox{$(-1)^{k+1}D_x^{2k}(|f|^nf)+|f|^nf=\frac{1}{n}\, f$}.
\end{equation}
For $k=1$, equation \eqref{EvenSS} possesses the explicit
compactly supported solution
\begin{equation}
 \label{f11}
f(x)=
 \left\{
 \begin{matrix}
\big[\frac{2(n+1)}{n(n+2)}\cos^2\big(\frac{nx}{2(n+1)}\big)\big]^{\frac{1}{n}}
 \,\,\, \mbox{for} \,\,\, |x| < \frac {\pi(n+1)}n,\smallskip \\
 0  \,\,\, \mbox{for} \,\,\, |x| \ge \frac {\pi(n+1)}n.
 \qquad\qquad\qquad\qquad\quad
 \end{matrix}
 \right.
\end{equation}
Thus, \eqref{SS1}, \ef{f11} yields  the so-called {\em standing
wave blow-up solution} (S-regime of blow-up) of \eqref{GenNonlin},
which  have compact support for all times $t \in (0,T)$. This
exact {\em Zmitrenko--Kurdyumov blow-up solution} has been known
since the middle of 1970s; see details in \cite[Ch.~4]{quasilin}.
As we have mentioned,  for $k \ge 2$, \eqref{EvenSS} cannot be
solved explicitly, but the ODE is shown to admit a countable set
of compactly supported solutions obtained by variational
Lusternik--Schnirel'man and fibering methods, \cite{GMPSob}.

Thus, comparing the ODEs \ef{EvenSS} and \ef{intnontwlode}, we see
that these coincide provided the compacton velocity $\l$ is given
by
\begin{equation}
 \label{l11}
\mbox{$\lambda =-\frac{1}{n}$}.
\end{equation}
In other words, this means that some principles of compacton
propagation for such NDEs \ef{GenNonlin} are directly related to
blow-up formation mechanisms for reaction-diffusion equations
\ef{ParaEven}. Moreover, in both equations, such asymptotic
structures are expected to be structurally stable (for $k=1$ in
\ef{ParaEven}, this has been proved, \cite[p.~260]{quasilin}).
 This is indeed surprising, since both classes of
nonlinear PDEs seem to be responsible for entirely different
physical phenomena (to say nothing of their different mathematical
essence).

The lower order case of equation \eqref{ParaEven} for $k=1$, which
is just a standard reaction-diffusion PDE, is fairly well
understood. However, the third-order nonlinear dispersion equation
in \eqref{GenNonlin}, also for the minimal value $k=1$, has not
been studied extensively and some of its basic governing
mathematical principles are still relatively unknown.

 \subsection{Two main model NDEs and layout of the paper}

 Thus, we will construct  global similarity solutions of the
 following two NDEs. The first one is the pure NDE
\begin{equation}
\label{linnonl}
  u_t = (-1)^{k+1}D^{2k+1}_x(|u|^nu) \quad \text{in}
\quad \mathbb{R}\times\mathbb{R}_+, \quad n>0,
\end{equation}
 which is studied in Sections \ref{S2}--\ref{S4}.
 Our goal therein is to show that (\ref{linnonl}) admits an infinite {\em countable}
 family of self-similar solutions governed by ``nonlinear eigenfunctions" of a rescaled
 operator. Moreover, in Section \ref{S3.5}, we show that this
 countable family
   as $n \to 0^+$ can be
 described by eigenfunctions from Hermitian spectral theory
 of the non-self adjoint operator \cite[\S~4]{RayGI}
  \be
  \label{HermB}
   \tex{
    \BB= (-1)^{k+1}D^{2k+1}_y + \frac 1{2k+1}\, y D_y + \frac
    1{2k+1} \, I.
    }
    \ee
The above operator occurs after similar scaling of the
{\em linear dispersion equation} (LDE)
\begin{equation}
\label{airyfunc}
 u_t = (-1)^{k+1}D_x^{2k+1}u \quad \text{in} \quad
\mathbb{R}\times\mathbb{R}_+; \quad k \ge 1,
\end{equation}
so we suggest to describe nonlinear eigenfunctions of the NDE
\ef{linnonl} by using linear theory.

 In Section \ref{S4}, we perform the opposite
 ``nonlinear" limit $ n \to + \infty$ in \ef{linnonl}.
Namely, we observe that the natural change of the independent
variable leads to the following PDE:
 \be
 \label{v11}
\tex{
 |u|^n u=v \LongA \big(|v|^{-\frac n{n+1}} v\big)_t=
 (-1)^{k+1}D^{2k+1}_x v.
 }
 \ee
The formal limit $n \to +\infty$ leads to the following {\em limit
NDE}:
 \be
 \label{v12}
 ({\rm sign}\, v)_t= (-1)^{k+1}D^{2k+1}_x v,
 \ee
 which admits analogous similarity patterns. It turns out that, at least
 for $k=1$, where \ef{v12} takes a particularly simple form
  \be
 \label{v12k1}
 ({\rm sign} \, v)_t= v_{xxx},
  \ee
  first three occurring ODEs for \ef{v12k1} admit
 reducing to an algebraic system, and we develop  a
 geometric-algebraic approach to constructing first nonlinear
 eigenfunctions.


In Section \ref{S5}, we return to the study of VSSs of the NDE
with absorption \ef{nde33}, where our main goal is to justify
existence of the so-called $p$-bifurcation branches of VSS
solutions, which appear at some critical exponents $p=p_l(n)>n+1$,
$l=0,1,2,...\,$.


\section{Similarity solutions of the NDE (``nonlinear eigenfunctions")}
 \label{S2}

 Here, we consider the pure NDE \ef{linnonl}, which is connected to
\eqref{GenNonlin}, but now we do not have the convection-like term
$(|u|^nu)_x$, which is negligible in our asymptotics.  This
nonlinear dispersion equation may be compared with the even-order
model, which represents the general higher-order porous medium
equation (the PME--2$m$)
\begin{equation}
 \label{PME2m}
u_t = (-1)^{m+1} \Delta^{m}(|u|^{n}u) \quad \text{in} \quad
\mathbb{R}^N\times\mathbb{R}_+, \quad m \ge 2.
\end{equation}
The classic PME--2, for $m=1$, appears in a number of physical
applications, such as fluid and gas flows, heat transfer, or
(nonlinear) diffusion. Other applications have been proposed in
mathematical biology, lubrication and  boundary layer theory, and
various other fields. Papers exploring delicate asymptotics for
the PME include \cite{EvoCompNonlEigPME}, where extra references
are available. For nonlinear eigenfunctions of (\ref{PME2m}) with
$m=2$, $N=1$,
 see \cite{GalPMEn} and key references therein. Overall, the higher-order case $m \ge 2$ in
 (\ref{PME2m}) has been studied much less in the mathematical
 literature and represents a number of difficult
  open problems.

 \subsection{Self-similar solutions: towards a ``nonlinear eigenvalue problem"}

Our NDE \eqref{linnonl}  admits standard similarity solutions
given by
\begin{equation}
\label{linnonlss} u_{\rm gl}(x,t)=t^{-\alpha}f(y) \whereA
 y=x/t^{\beta},
\end{equation}
for some unknown real parameters $\alpha$ and $\beta$.  After
substitution \ef{linnonlss} into the NDE, we obtain the ODE for
the rescaled profile $f$,
\begin{equation}
\label{SSNonlOde} -\alpha t^{-\alpha-1}f - \beta
t^{-\alpha-1}f^{\prime}y =
(-1)^{k+1}t^{-\alpha(n+1)-\beta(2k+1)}D_y^{2k+1}(|f|^nf).
\end{equation}
By equating powers of $t$, the parameter $\beta$ can be found in
terms of $\alpha$ and is given by
\begin{equation}
  \label{b112}
\mbox{$\beta = \frac{1-\alpha n}{2k+1} > 0$}
 \quad \mbox{for} \quad \mbox{$\alpha < \frac{1}{n}$},
\end{equation}
so that now $\a \in \re$ remains the only unknown {\em nonlinear
eigenvalue}.

 Thus, the ODE \eqref{SSNonlOde} takes the form
\begin{equation}
\label{rednonlOde}
 \AAA(f,\a) \equiv
    \mbox{$(-1)^{k+1}D_y^{2k+1}(|f|^nf)
+\frac{1-\alpha n}{2k+1}\, f^{\prime}y + \alpha f = 0$} \inB \re.
\end{equation}
In order to formulate a suitable {\em nonlinear eigenvalue
problem} for the pairs $\{\a_l,f_l\}_{l \ge 0}$ for the ODE
\ef{rednonlOde}, one should specify the ``boundary" conditions at
$y = \pm \iy$, which is one of the main goals of this paper. It
turns out that, loosely speaking, such conditions can be
formulated as follows: for some special values of the real
eigenvalues $\{\a_l=\a_l(n)>0\}_{l \ge 0}$, the corresponding
nonlinear eigenfunctions
 \be
 \label{BC1}
  \left\{
  \begin{matrix}
  \mbox{ $f_l(y)$ have: \, (i) finite left-hand
  interface ($f_l(y) \equiv 0$ for $y \ll -1$),}\ssk \ssk \\
   \mbox{and (ii)  ``minimal oscillatory" behaviour as $y \to
   +\iy$. \qquad\quad\,\,\,\,}
    \end{matrix}
  \right.
   \ee
Both such conditions will get proper explanations later on. Of
course, at least for sufficiently large $l \ge 2k+1$, the problem
\ef{rednonlOde}, \ef{BC1} exhibits all typical features of
self-similarity of the {\em second kind} (the {\em first kind}
corresponds to standard dimensional analysis of PDEs), a notion
introduced by Ya.B.~Zel'dovich in 1956 \cite{Zel56}. Unlike
previous and known examples, in general, the eigenvalues
$\alpha_l(n)$ (and hence $\beta_l(n)$) cannot be found explicitly
from any dimensional analysis. Thus, admitted (real) values of
$\alpha= \a_l(n)>0$ at this stage are still unknown and play a
role of ``nonlinear eigenvalues".

In particular, our goal is to show by a combination of analytic,
formal, and numerical tools that, for any $n>0$, the eigenvalue
problem (\ref{rednonlOde}), (\ref{BC1})
 \be
 \label{BC2}
 \mbox{
 has a countable set of pairs $\{\a_l,f_l(y)\}_{l \ge 0}$: \,\,
  \,\, $\AAA(f_l,\a_l)=0$}.
  \ee

For the linear case with $n=0$, the corresponding linear
non-self-adjoint Hermitian spectral problem was developed in
\cite[\S~4]{RayGI}; we present and use these results below.

\subsection{Towards blow-up patterns}

The present analysis admits a natural duality: since the NDE
\ef{linnonl} is invariant under {\em reflections},
 \be
 \label{refl1}
 \left\{
 \begin{matrix}
  x \mapsto -x,\quad \\
  t \mapsto T-t,
  \end{matrix}
  \right.
  \ee
 the global patterns \ef{linnonlss} lead to the corresponding reflected
 {\em blow-up} ones:
 \be
 \label{bl11}
  \tex{
  u_{\rm bl}(x,t)=(T-t)^{-\a}f(y), \quad y=- \frac
  x{(T-t)^\b} \whereA \b=\frac{1-\a n}{2k+1}>0.
  }
  \ee
  Therefore, the nonlinear eigenvalue problem \ef{rednonlOde},
  \ef{BC1} is assumed to describe both countable families of
  {global} and {blow-up} patterns for the NDE \ef{linnonl}.
  In what follows, for simplicity, we will mainly use a global treatment of such
  asymptotic patterns.

\subsection{A homotopy path $n \to 0^+$ to Hermitian spectral
theory: a route to countability of nonlinear eigenfunctions}

 For $n=0$, $\alpha$ gives  the eigenvalues
$\{\lambda_l\}$ in linear Hermitian spectral theory,
\cite[\S~4]{RayGI}. Indeed, for $n=0$, \ef{rednonlOde} reads
\begin{equation}
\label{N1}
 \mbox{$(-1)^{k+1}D_y^{2k+1}f +\frac{1}{2k+1}\,
f^{\prime}y +\a f \equiv \BB f+ \big(\a- \frac 1{2k+1}\big) f
=0$}.
\end{equation}
 Therefore, this gives the eigenvalue equation for the linear  operator
 $\BB$ in \ef{HermB}:
  \be
   \label{N2}
   \tex{
    \BB \psi= \l \psi \whereA \l= -\a +\frac 1{2k+1}.
    }
    \ee
By Hermitian spectral theory \cite{RayGI}, this defines a
countable set of eigenvalues for $n=0$:
 \be
 \label{N3}
  \tex{
  \a_l(0) = \frac {1+l}{2k+1}, \quad l=0,1,2,...\,.
   }
   \ee
The corresponding eigenfunctions are then the normalized
derivatives of the rescaled kernel $F(y)$ of the fundamental
solution of the dispersion operator $D_t+(-1)^k D_x^{2k+1}$:
 \be
 \label{N4}
  \tex{
  \psi_l(y) = \frac{(-1)^l}{\sqrt{l !}}\, D^l_y F(y), \quad l \ge
  0.
  }
  \ee

Moreover,  the ``adjoint" (not in the standard $L^2$-sense)
operator
 \be
 \label{N5}
 \tex{
 \BB^*= (-1)^{k+1} D_y^{2k+1}- \frac 1{2k+1} \, y D_y
 }
 \ee
 has the same spectrum and eigenfunctions $\{\psi_l^*\}$, which are {\em generalized
 Hermite polynomials} given by (see \cite[\S~5.2]{RayGI})
\begin{equation}
\label{poleigfunc}
 \mbox{$
  \mbox{$\psi ^{\ast}_{l}(y) =
\frac{1}{\sqrt{l !}}$}\Big[\mbox{$y^{l} +$} (-1)^{k+1}\sum\limits
_{j=1}^{\lfloor \frac{|l
|}{2k+1}\rfloor}\mbox{$\frac{1}{j!}D_y^{(2k+1)j}y^{l}$}\Big].
 $}
\end{equation}
 Finally, the operator pair
 $\{\BB,\BB^*\}$ has a bi-orthonormal sets of eigenfunctions.

\ssk

As a key conclusion,
 we expect that there exists a uniform and pointwise
convergence as $n \to 0^+$ of the nonlinear eigenfunctions to the
linear ones for $\BB$ defined by
 \ef{N4}.
 We postpone explanations of main aspects of such a branching
 until
 Section \ref{S3.5}, when we will have gained enough understanding
of basic nonlinear eigenfunction theory involved.
 However, a rigorous proof of such a
branching of eigenfunctions at $n=0$ is a hard open problem.
Nevertheless, this $n$-branching approach still remains a unique
analytically convincing fact towards existence of an {\em infinite
countable} discrete
 set of nonlinear eigenfunctions of the problem
\ef{rednonlOde}. In this case, branching theory as $n \to 0^+$ can
be developed along the same lines as for the PME--4 \ef{PME2m}
with $m=2$ \cite[\S~6]{GalPMEn}. Such a construction then uses a
specific Hermitian spectral theory created in \cite{RayGI}; see
Section \ref{S3.5}.

\subsection{Conservation laws: explicit values of  first
nonlinear eigenvalues}

As a first pleasant surprise,  it turns out that some first
eigenvalues for any $n>0$ can be calculated explicitly by
conservation laws for the NDE.
Assuming that the solution $u(x,t)$ is integrable in $x$ over
$\re$ (this is a formal assumption that can be got rid of), we
have that \eqref{linnonl} is conservative in mass, and so
\begin{equation}
\label{cons}
\mbox{$\frac{\mathrm{d}}{\mathrm{d}t}$} \,
 \mbox{$
\int\limits_{\mathbb{R}} \mbox{$u(x,t)\, \mathrm{d}x$} = 0.
 $}
\end{equation}
For similarity solutions \eqref{linnonlss}, we have that
\begin{equation*}
\mbox{$ \int\limits_{\mathbb{R}} $} u(x,t)\, \mathrm{d}x =
t^{\beta - \alpha} \mbox{$ \int\limits_{\mathbb{R}} $} f(y)\,
\mathrm{d}y.
\end{equation*}
This satisfies \eqref{cons}, provided that, in addition to
\ef{b112},
\begin{equation}
 \label{al0}
\mbox{$ \beta-\a = 0 \quad \Longrightarrow \quad \alpha_0(n) =
\frac{1}{(2k+1)+n}$},
\end{equation}
for non-zero rescaled mass $\int f\neq 0$. So, on substitution
into \eqref{rednonlOde},
\begin{equation*}
\mbox{$(-1)^{k+1}D_y^{2k+1}(|f|^nf) +\frac{1}{(2k+1)+n}\,
f^{\prime}y + \frac{1}{(2k+1)+n}\, f = 0$}.
\end{equation*}
   Integrating once with the zero constant (a zero-flux condition), we end up with the ODE
\begin{equation}
 \label{ff1}
\mbox{$(-1)^{k+1}D_y^{2k}(|f|^nf) + \frac{1}{(2k+1)+n}\, fy = 0$}.
\end{equation}
Note that, for $n=0$, we have exactly the linear ODE (see
\cite{RayGI})
\begin{equation}
\label{lineqn} \mbox{$(-1)^{k+1}F^{(2k)}+ \frac{1}{2k+1}\,Fy = 0$}
\quad \text{for} \quad y\in\mathbb{R}.
\end{equation}

For convenience, we use in (\ref{ff1}) the natural substitution
\begin{equation}
\label{nonlsub} Y=|f|^nf \LongA \mbox{$f=|Y|^{-\frac{n}{n+1}}Y$},
\end{equation}
in order to remove nonlinearities in the highest differential.
Substitution yields
\begin{equation}
\label{massconsY} \mbox{$(-1)^{k+1}D_y^{2k}Y +
\frac{1}{(2k+1)+n}\, y|Y|^{-\frac{n}{n+1}}Y = 0$}.
\end{equation}

Similarly, we have conservation of the first moment, with
\begin{equation*}
\mbox{$ \int\limits_{\mathbb{R}} $} x u(x,t)\, \mathrm{d}x =
t^{2\beta - \alpha} \mbox{$ \int\limits_{\mathbb{R}} $} y f(y)\,
\mathrm{d}y.
\end{equation*}
Hence, we have that
\begin{equation}
 \label{al1}
  2\b-\a=0 \LongA
\mbox{$
\alpha_1(n) =\frac{2}{(2k+1)+2n}$}.
\end{equation}
This then gives the ODE
\begin{equation}
\label{1stmoment} \mbox{$(-1)^{k+1}D_y^{2k+1}(|f|^nf)
+\frac{1}{(2k+1)+2n}\, f^{\prime}y + \frac{2}{(2k+1)+2n}\, f =
0$}.
\end{equation}
However, we cannot simply integrate this equation, as we could
before, to reduce the order of the ODE.  Instead, we multiply
\eqref{1stmoment} by $y$, so that
\begin{equation*}
\mbox{$(-1)^{k+1}D_y^{2k+1}(|f|^nf)y +\frac{1}{(2k+1)+2n}\,
f^{\prime}y^2 + \frac{2}{(2k+1)+2n}\, fy = 0$},
\end{equation*}
and now it is possible to integrate by parts, to obtain
\begin{equation}
\label{1stmomintcons} \mbox{$(-1)^{k+1}D_y^{2k}(|f|^nf)y +
(-1)^{k}D_y^{2k-1}(|f|^nf)+ \frac{1}{(2k+1)+2n}\, fy^2 = 0$}.
\end{equation}


We also look at conservation of the second moment,
\begin{equation*}
\mbox{$ \int\limits_{\mathbb{R}} $} x^2 u(x,t)\, \mathrm{d}x =
t^{3\beta - \alpha} \mbox{$ \int\limits_{\mathbb{R}} $} y^2 f(y)\,
\mathrm{d}y.
\end{equation*}
This gives $3\b-\a=0$, so that, by \ef{b112},
\begin{equation}
\label{2ndmoment} \mbox{$\alpha_2(n) =\frac{3}{(2k+1)+3n}$} \quad
\mbox{and} \quad
 \mbox{$(-1)^{k+1}D_y^{2k+1}(|f|^nf) + \frac{1}{(2k+1)+3n}\,
 f^{\prime}y+
\frac{3}{(2k+1)+3n}\, f  = 0$}.
\end{equation}
Similarly,  we multiply by $y^2$ and integrate to reduce the
order, to obtain the ODE
\begin{equation}
\begin{split}
\label{2ndmomintcons} \mbox{$(-1)^{k+1}$}&\mbox{$D_y^{2k}(|f|^nf)y^2
+
2(-1)^{k}D_y^{2k-1}(|f|^nf)y$}\\[1mm]
&\mbox{$+2(-1)^{k+1}D_y^{2k-2}(|f|^nf)+ \frac{1}{(2k+1)+3n}\, fy^3 =
0$}.
\end{split}
\end{equation}

These three conservation laws in particular are important, as we
can explicitly find the first three (second-order) equations for
the case $k=1$ (corresponding to the first three nonlinear
eigenvalues), but not for others. The case $k=1$ is essential,
since it is much easier to develop theory for the lower-order
case, as well as it is
 easier to solve numerically (see Section
\ref{NumConstrNonl}).

\ssk

 In general, for arbitrary $k \ge 1$, for all $l<2k+1$, where $l$ is the eigenvalue index
as before, we have our moments conservation given by
\begin{equation*}
\mbox{$ \int\limits_{\mathbb{R}} $} x^l u(x,t)\, \mathrm{d}x =
t^{(l+1)\b-\alpha} \mbox{$ \int\limits_{\mathbb{R}} $} y^l f(y)\,
\mathrm{d}y.
\end{equation*}
Therefore, our nonlinear eigenvalues are represented by
\begin{equation}
 \label{alk}
  (l+1)\b-\a=0 \LongA
\mbox{$\alpha_l(n) = \frac{l+1}{(2k+1)+(l+1)n}$} \quad \mbox{for
any} \quad  0\leq l < 2k+1.
\end{equation}
 The corresponding ODEs, which can be integrated once for all such eigenvalues,
  are
\begin{equation}
 \label{D110}
\mbox{$(-1)^{k+1}D_y^{2k+1}(|f|^nf)  +\frac{1}{(2k+1)+(l+1) n}\,
f^{\prime}y + \frac{l+1}{(2k+1)+(l+1)n}\, f = 0$}.
\end{equation}
All of them  admit reducing the order by multiplying by $y^l$ and
integrating by parts $l$ times.
 The function  $Y=|f|^n f$
then solved the semilinear PDE
 \be
 \label{D11}
  \tex{
  (-1)^{k+1} Y^{(2k+1)} +\frac{1}{(2k+1)+(l+1) n}\,
y (|Y|^{-\frac n{n+1}}Y)' + \frac{l+1}{(2k+1)+(l+1)n}\,
|Y|^{-\frac n{n+1}}Y = 0. }
 \ee

Thus, for such NDEs, one can always obtain explicitly $2k+1$
$n$-branches of nonlinear eigenvalues \ef{alk}, though of course
the solvability of the corresponding eigenvalue equations remains
a difficult open problem, especially for large $l$. However,
\ef{alk} shows one important feature of this nonlinear eigenvalue
theory: the $n$-branches \ef{alk} are {\em global} in $n>0$, i.e.,
exist for all $n>0$. Existence of such a global continuation, with
no turning points, is one of the most difficult question in
general nonlinear operator theory; see e.g., \cite{GeoMethNonAn,
VainbergTr}.
 Recall that the eigenvalue problem \ef{rednonlOde} is by no means
 variational, neither contains any monotone or coercive operators.

 Next, we must admit that, for any larger  $l\geq2k+1$, we cannot find $\alpha_l(n)$
explicitly using conservation laws, and then, searching for
$n$-branches of these nonlinear eigenfunctions, we either should
rely on the above local homotopy approach as $n \to 0^+$, or use
advanced numerical methods.


\subsection{Numerical construction of nonlinear eigenfunctions}
\label{NumConstrNonl}

We look to use a shooting method in order to find reliable
profiles for the NDE \eqref{linnonl}.  First, considering
solutions satisfying  conservation of mass, we have our first
``nonlinear eigenvalue" (for $l=0$), where $n=0$ corresponds to
the linear kernel, $\psi_0=F(y)$, i.e., $\l=0$ in \ef{N2}. Here
the rescaled equation is given by \eqref{massconsY}.

\ssk

Consider first  the simpler case $k=1$, when we have the
second-order equation:
\begin{equation}
\label{2ndOrdNonl}
Y^{\prime\prime}=\mbox{$-\frac{1}{n+3}\,|Y|^{-\frac{n}{n+1}}Yy$}.
\end{equation}
Then, it can be shown that the left-hand interface is not
oscillatory (in fact, the same is true for $n=0$, i.e., for the
Airy function). Hence,
  for small solutions of $Y(y)$ as $y \to y_0^+$, with some
interface point $y_0<0$, we can approximate it by
\begin{equation}
\label{shtsol} Y(y)= C_0 (y-y_0)_+^{\tilde{\alpha}}\big(1+o(1)\big),
\end{equation}
for some constant $C_0>0$ and exponent $\tilde{\alpha}>0$.  Here
we have used the  standard notation
\begin{equation*}
(\cdot)_+ = \max\{0,\cdot\}.
\end{equation*}
Substituting \eqref{shtsol} into the ODE \eqref{2ndOrdNonl}, we
obtain the leading equality
\begin{equation}
\label{alphainserteqn}
\tilde{\alpha}(\tilde{\alpha}-1)C_0(y-y_0)^{\tilde{\alpha}-2}=
\mbox{$-\frac{1}{n+3}\,C_0^{\frac{1}{n+1}}(y-y_0)^{\frac{\tilde{\alpha}}{n+1}}y_0$},
\end{equation}
where we use $(y-y_0)$ to mean $(y-y_0)_+\big(1+o(1)\big)$, as
defined before.  Hence, we must have from \eqref{alphainserteqn}
that
\begin{equation*}
\tilde{\alpha} = \mbox{$\frac{2(n+1)}{n}$},
\end{equation*}
with $2<\tilde{\alpha}\leq 4$    for $n\geq 1$.  From this, the
constant $C_0$ is given by
\begin{equation*}
C_0=\mbox{$\big(\frac{n^2|y_0|}{2(n+1)(n+2)(n+3)}
\big)^{\frac{n+1}{n}}$} \quad \text{for any interface} \quad
y_0<0.
\end{equation*}
Our solution then has the expansion
\begin{equation}
 \label{Yfbp}
Y(y)=\mbox{$\big(\frac{n^2|y_0|}{2(n+1)(n+2)(n+3)}
\big)^{\frac{n+1}{n}}\,(y-y_0)^{\frac{2(n+1)}{n}}$}(1+o(1)) \asA y
\to y_0^+,
\end{equation}
with the derivative expansion obtained similarly,
\begin{equation}
\label{Yprime}
Y^{\prime}(y)=\mbox{$\frac{2(n+1)}{n}\big(\frac{n^2|y_0|}{2(n+1)(n+2)(n+3)}
\big)^{\frac{n+1}{n}}\,(y-y_0)^{\frac{n+2}{n}}$}(1+o(1)).
\end{equation}

We use the {\tt MatLab} IVP solver {\tt ode15s} to plot our
profiles.
 Taking an arbitrary initial (interface) point $y=y_0<0$, our solution and
first derivative will be zero here.  Since both initial values are
zero at $y=y_0$, we will often find the solution $Y=0$, whilst
trying to solve numerically.  In order to overcome this problem,
we must look at some point $y_0+\delta$ (for small $\delta\sim
10^{-3}$), close to this point, and after finding the derivative
there, this is used as the initial condition.

Obviously due to the nature of the expansion  \eqref{Yprime}, we
must take a relatively large initial point (in our case we take
$y_0=-10$), in order for us to have an initial condition that is
not negligible. Hence, the solution is a large rescaling of the
``nonlinear fundamental solution" with the unit mass.
 Below are a few profiles that have been found, in which we have
taken $\delta=10^{-3}$.  For $n\sim 0.5$,
 the derivative
$Y^\prime$ is very small and larger negative interface
 points must
be used to find reliable profiles.  This makes comparison, between
different values of $n$, more difficult.

\begin{figure}[htbp]
\begin{center}
\includegraphics[scale=0.70]{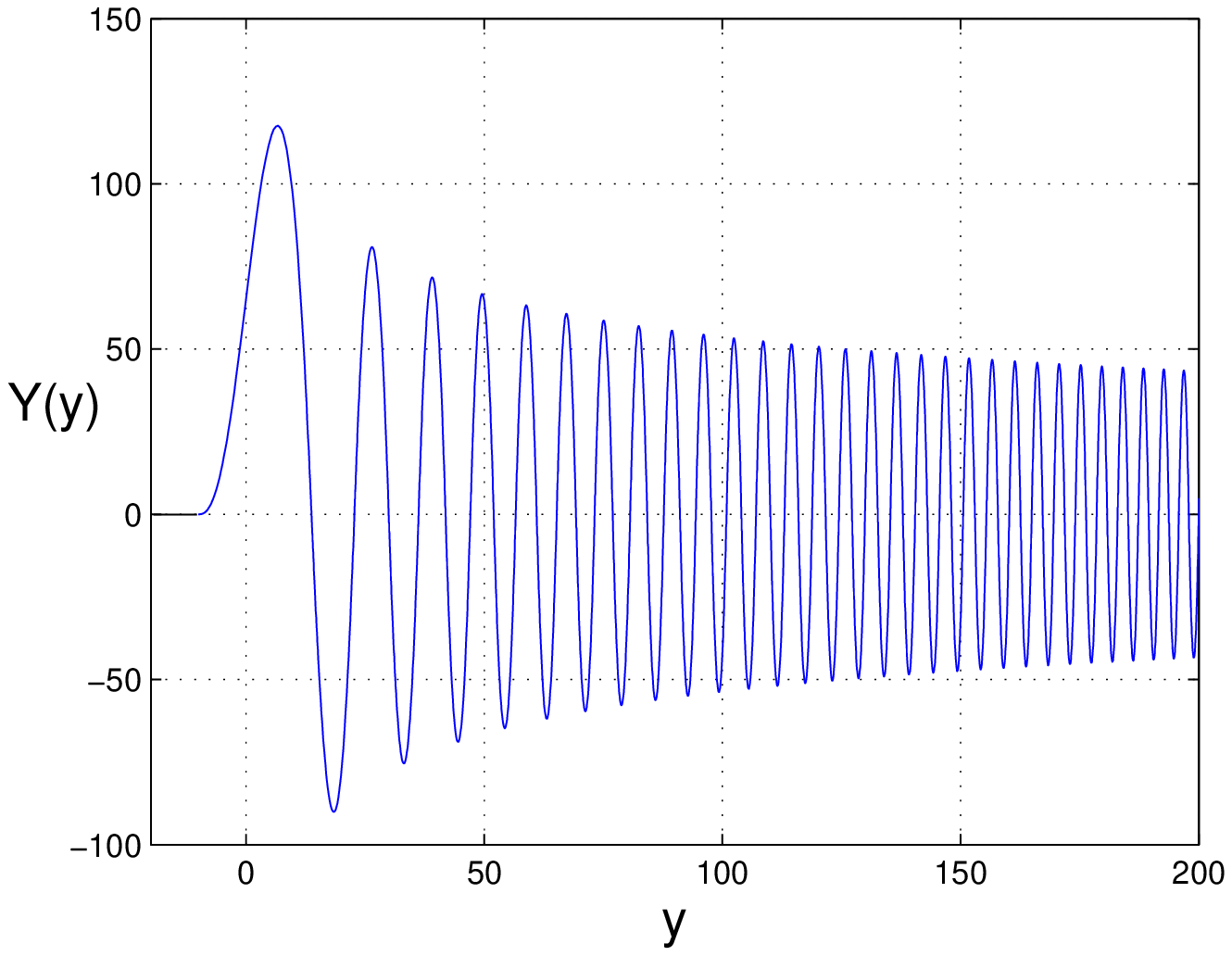}
\caption{\small Rescaled solution $Y(y)$ of the ODE
(\ref{2ndOrdNonl}) for $k=1$, $l=0$,  with $n=3$.}
 \label{F1}
\end{center}
\end{figure}

\begin{figure}[htbp]
\begin{center}
\includegraphics[scale=0.70]{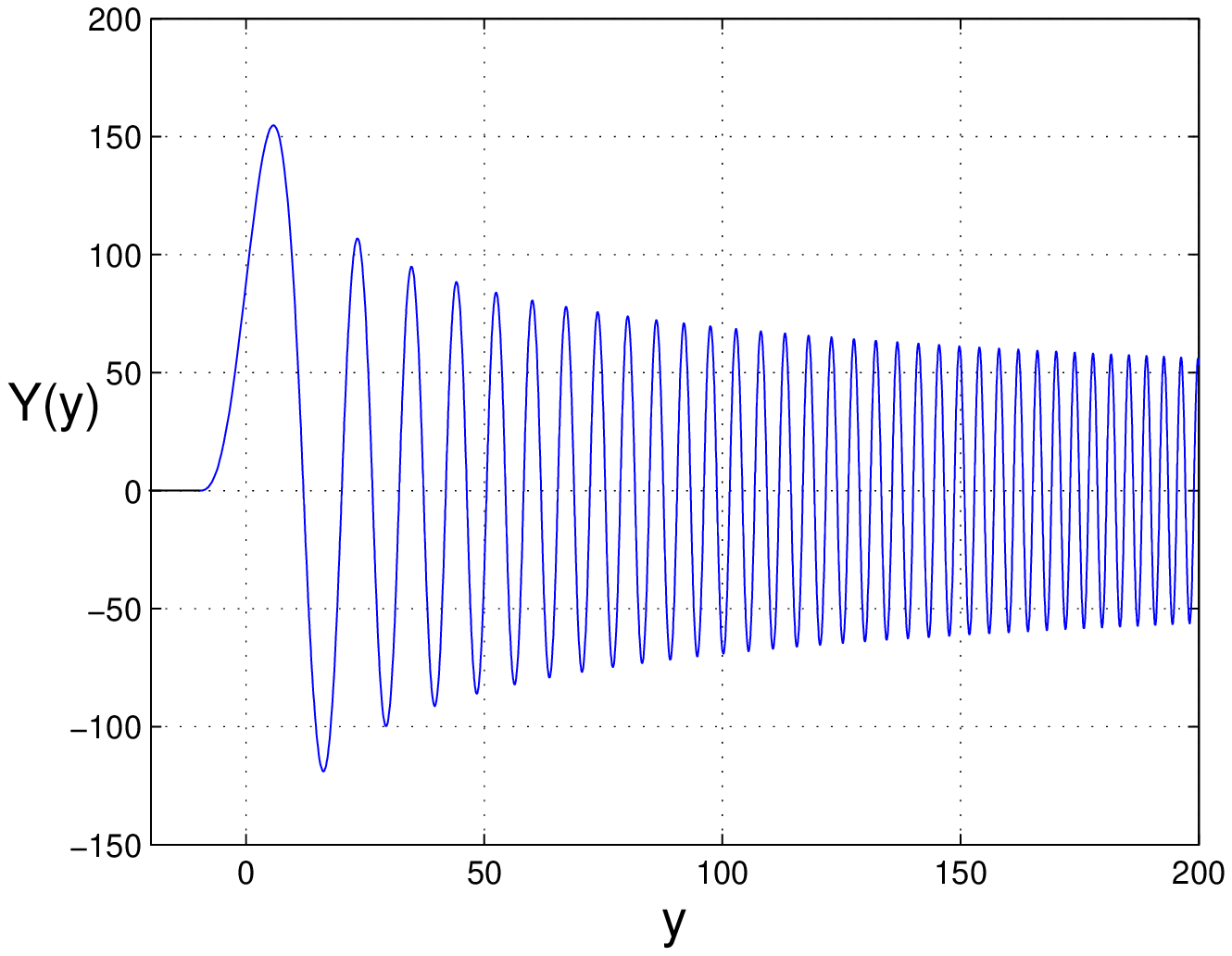}
\caption{\small Rescaled solution $Y(y)$ of the ODE
(\ref{2ndOrdNonl}) for $k=1$, $l=0$,  with $n=2$.}
 \label{F2}
\end{center}
\end{figure}

\begin{figure}[htbp]
\begin{center}
\includegraphics[scale=0.70]{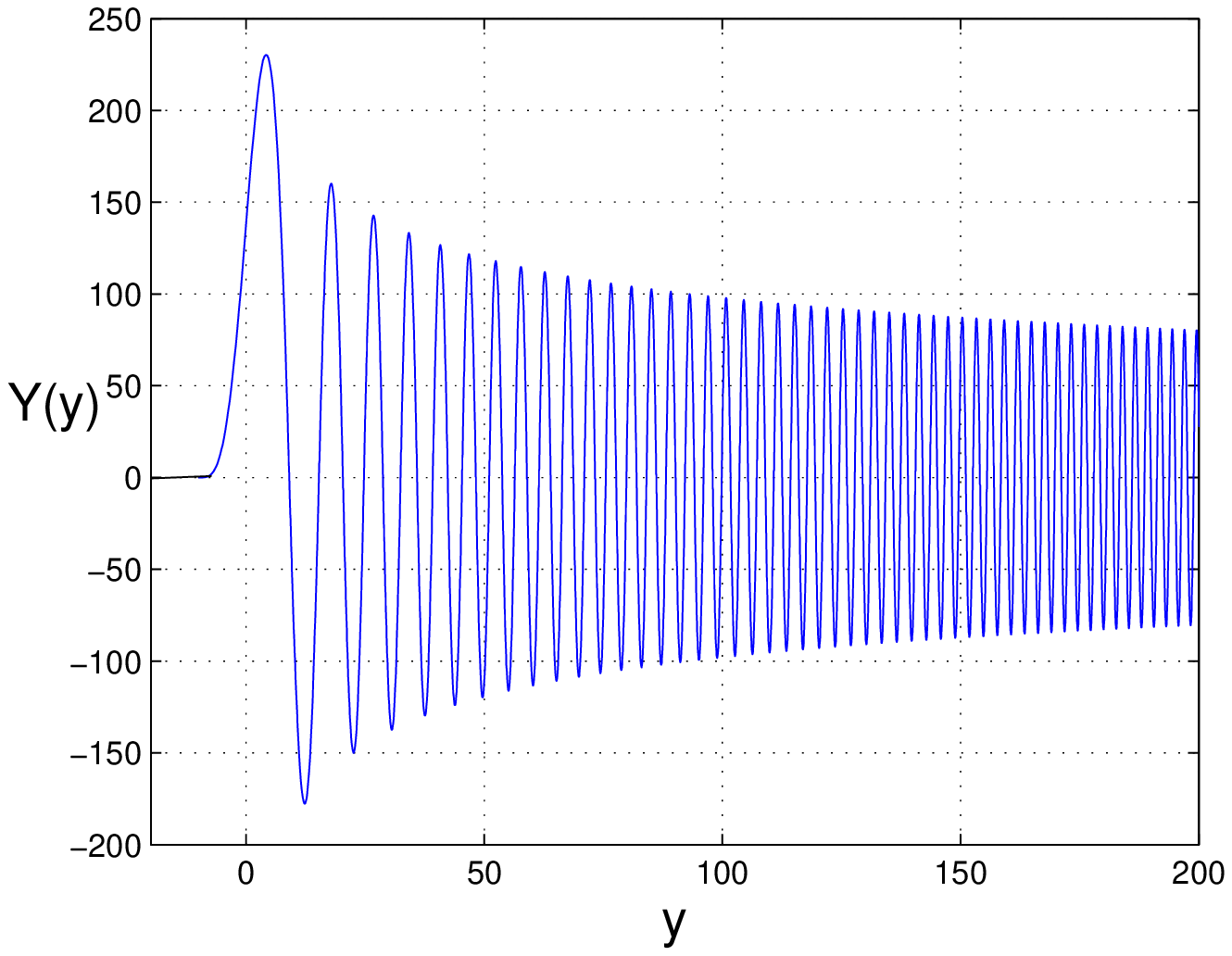}
\caption{\small Rescaled solution $Y(y)$ of the ODE
(\ref{2ndOrdNonl}) for $k=1$, $l=0$,  with $n=1$.}
 \label{F3}
\end{center}
\end{figure}

\begin{figure}[htbp]
\begin{center}
\includegraphics[scale=0.60]{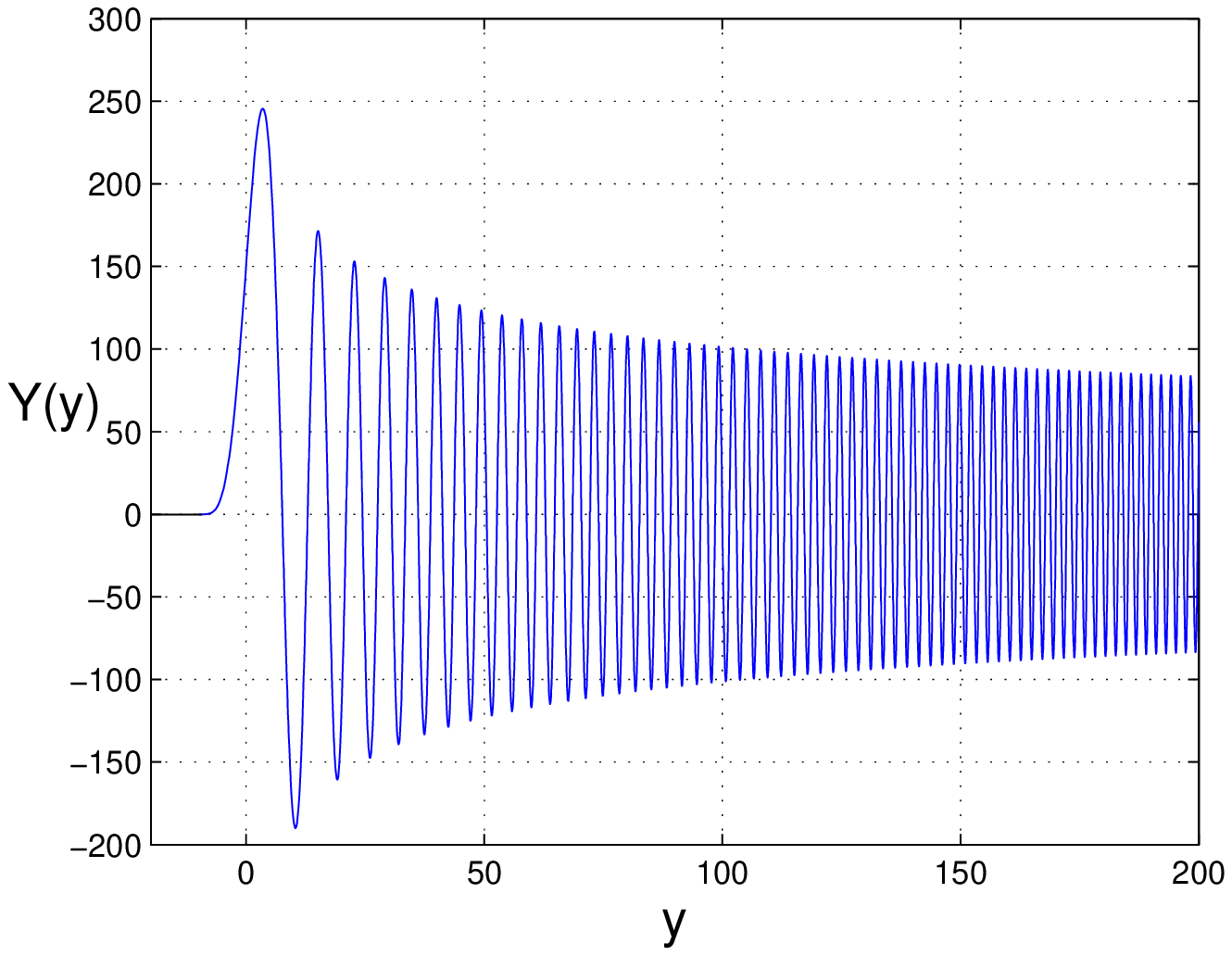}
\caption{\small Rescaled solution $Y(y)$ of the ODE
(\ref{2ndOrdNonl}) for $k=1$, $l=0$, with $n=0.7$.}
 \label{F4}
\end{center}
\end{figure}

In Figures \ref{F1}--\ref{F4}, we show the first nonlinear
eigenfunction $Y_0(y)$, with $k=1$, for various values of $n=3, 2,
1,$ and $0.7$. Note that this is precisely the profile that, as $n
\to 0^+$, must converge to the rescaled kernel $F(y) \equiv {\rm
Ai}\,(y)$, representing
 Airy's classic function satisfying (cf.
 (\ref{2ndOrdNonl}) for $n=0$)
  \be
  \label{Ai1}
   \tex{
  F''+ \frac 13\, F y=0, \quad \int F =1,
  }
  \,\,\, \mbox{where} \,\,\, F \in L^2_\rho(\re), \quad \rho(y)=
  \left\{
   \begin{matrix}
    {\mathrm e}^{a|y|^{3/2}} \,\, \mbox{for} \,\,y<0,\\
{\mathrm e}^{-a y^{3/2}} \,\, \mbox{for} \,\,y>0,
 \end{matrix}
  \right.
 \ee
 where $a>0$ is a small enough constant.

In general, for all values of $l$, we have, for the lower-order
case of $k=1$, that the same expansion \ef{Yfbp} and \ef{Yprime}
 hold near finite interfaces.
Here $l<2k+1$ and hence, for $k=1$, we can only take $l=0,1,2$ for
explicitly given nonlinear eigenvalues.

 For $l=1$, we have from \eqref{1stmomintcons}, with $Y(y)=|f|^n f$
and $k=1$, that the eigenfunction $Y_1(y)$ solves the ODE
\begin{equation}
 \label{Y1}
\mbox{$Y^{\prime\prime} = \frac{1}{y}\,Y -
\frac{1}{2n+3}\,|Y|^{-\frac{n}{n+1}}Yy$}.
\end{equation}
For $l=2$ and $k=1$, we have from \eqref{2ndmomintcons} that
$Y_2(y)$ solves
\begin{equation}
 \label{Y2}
\mbox{$Y^{\prime\prime} = \frac{2}{y}\,Y^\prime - \frac{2}{y^2}\,Y
-\frac{1}{3n+3}\,|Y|^{-\frac{n}{n+1}}Yy$}.
\end{equation}
We use these to plot the profiles of our NDE with $l=1$ and $l=2$,
which, for $n=0$,  coincide with the derivatives $F^\prime(y)$ and
$F^{\prime\prime}(y)$ respectively, of the linear kernel $F(y)$.

\begin{figure}[htbp]
\begin{center}
\includegraphics[scale=0.55]{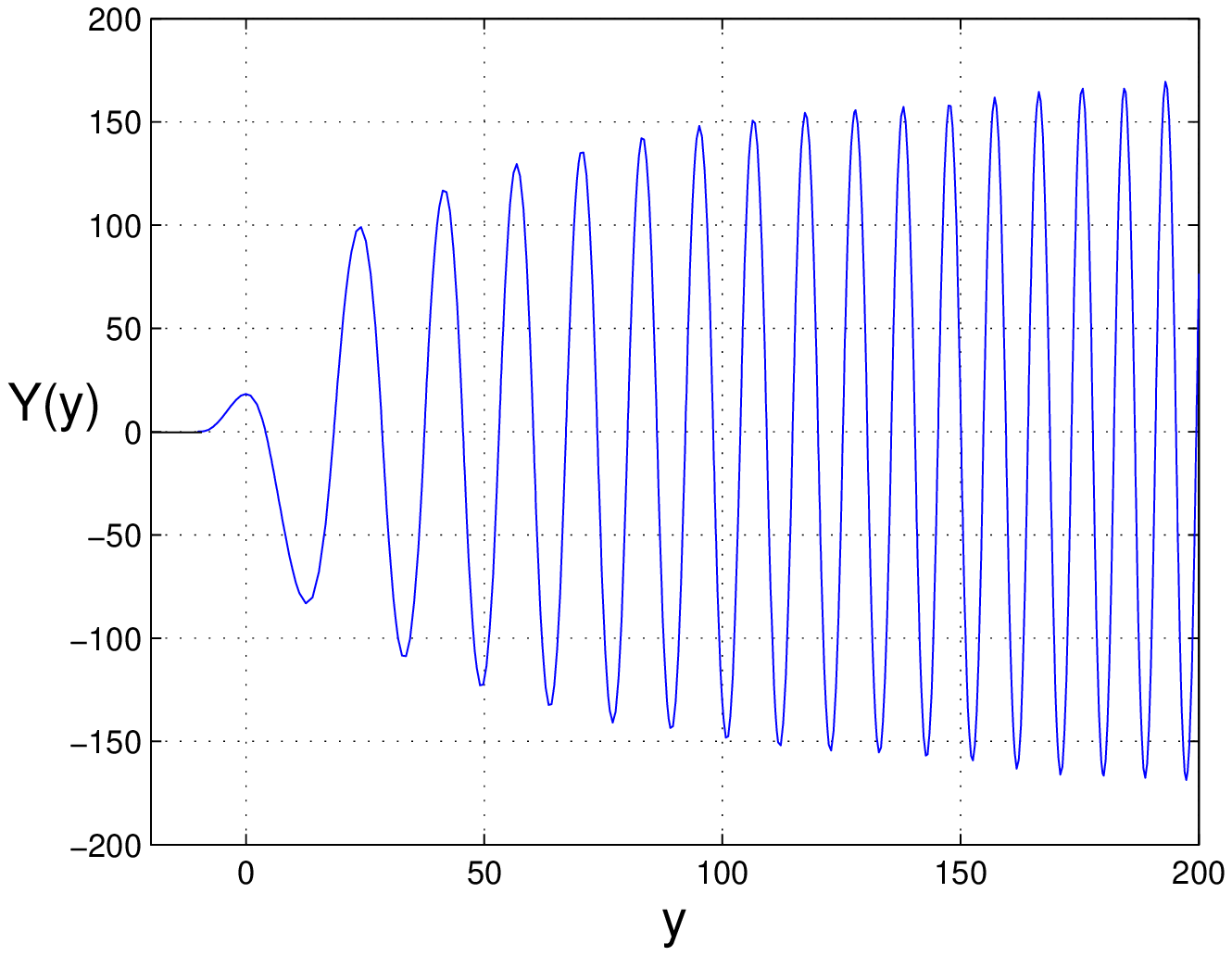}
\caption{\small Rescaled solution $Y(y)$ of the ODE (\ref{Y1}) for
$k=1$, $l=1$, with $n=3$.}
 \label{F5}
\end{center}
\end{figure}

\begin{figure}[htbp]
\begin{center}
\includegraphics[scale=0.55]{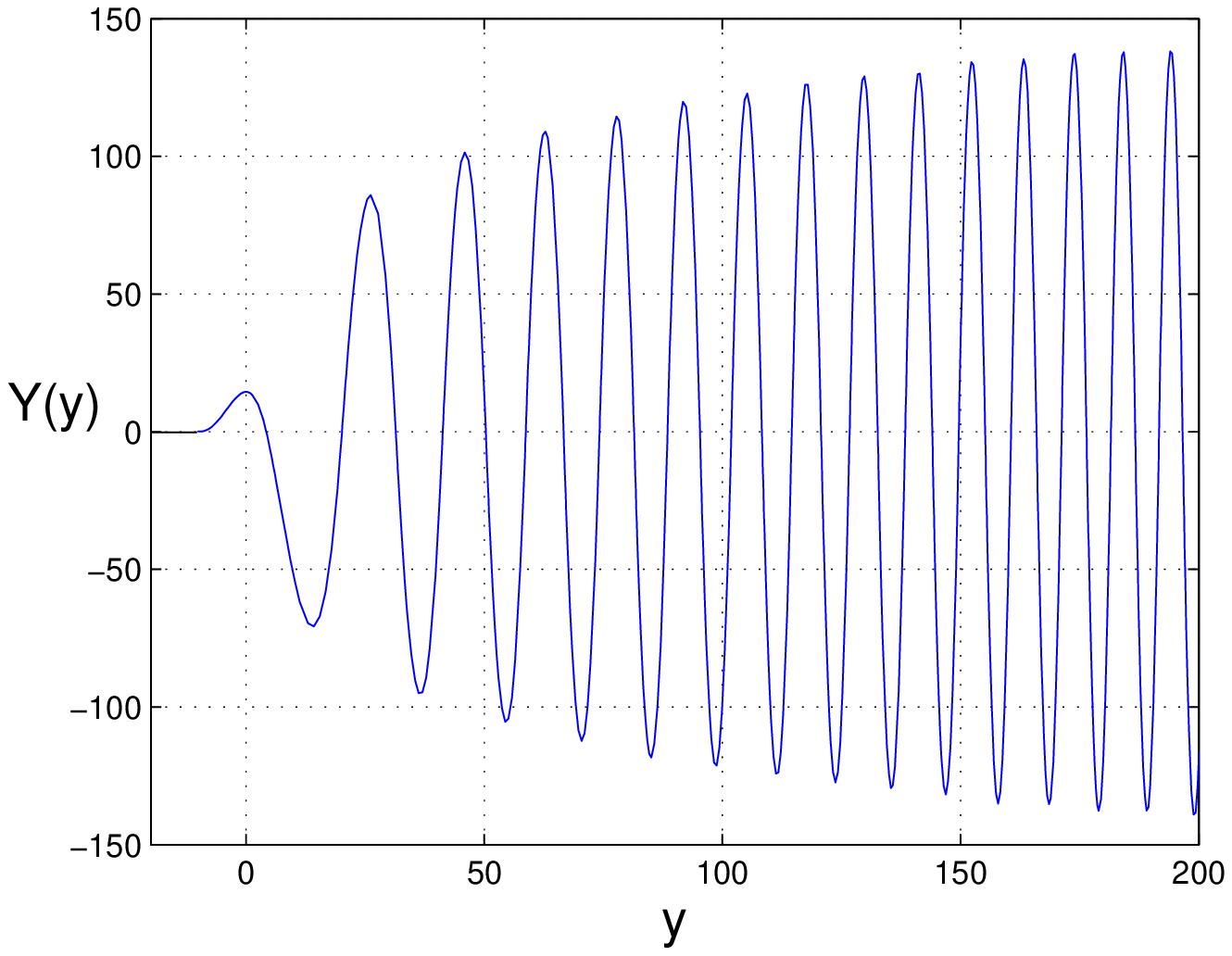}
\caption{\small Rescaled solution $Y(y)$ of the ODE (\ref{Y1}) for
$k=1$, $l=1$, with $n=4$.}
 \label{F6}
\end{center}
\end{figure}

\begin{figure}[htbp]
\begin{center}
\includegraphics[scale=0.55]{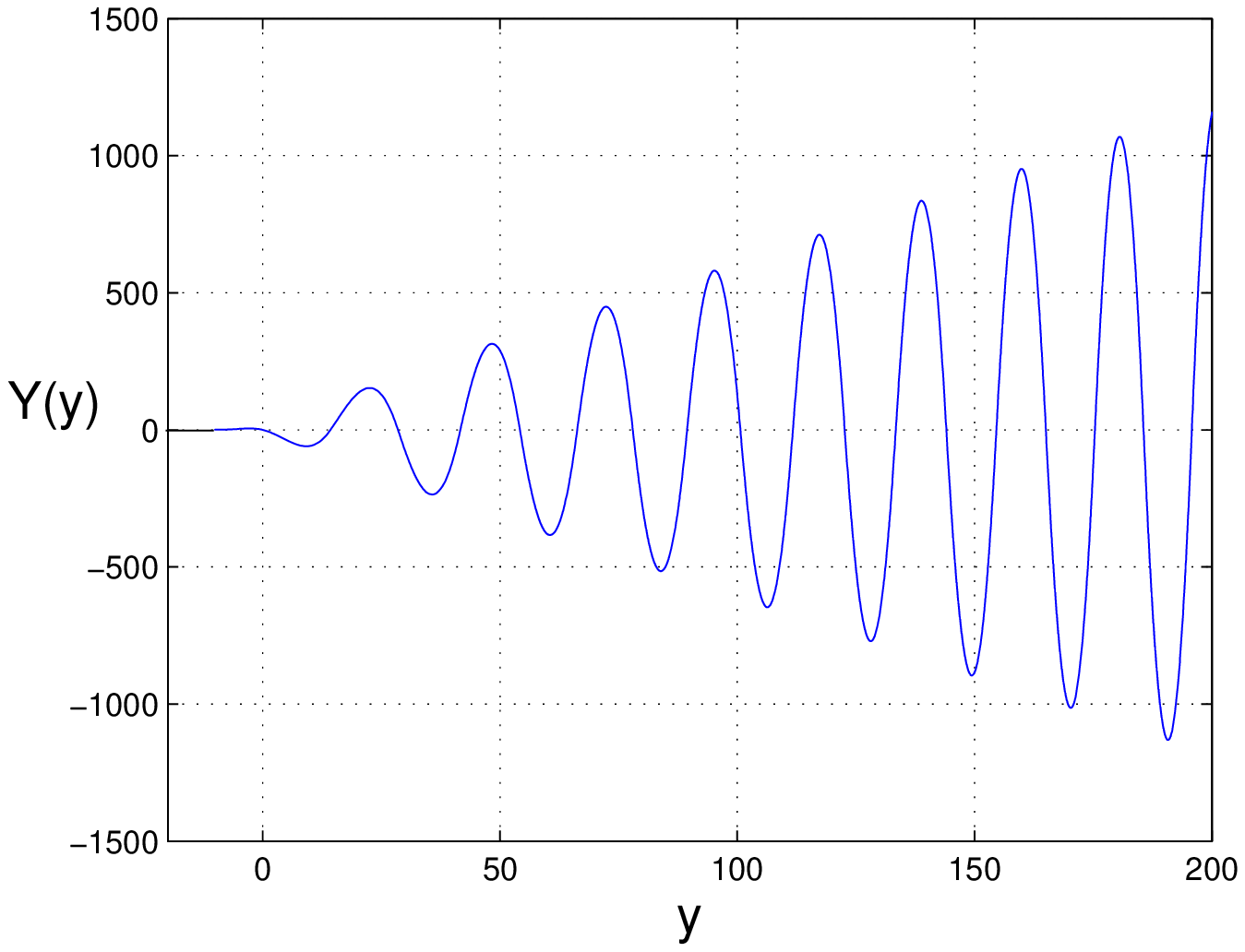}
\caption{\small Rescaled solution $Y(y)$ for $k=1$, $l=2$ of the
ODE (\ref{Y2})
 with $n=3$.}
  \label{F7}
\end{center}
\end{figure}

Figures \ref{F5} and \ref{F6} show  the ``dipole-like" profiles
for $l=1$ as solutions of \ef{Y1}, while Figure \ref{F7}
represents even more oscillatory third eigenfunction $f_2(y)$ for
$l=2$. Recall that, by a homotopy path as $n \to 0^+$, this
function is expected to converge to the highly oscillatory linear
eigenfunction of the operator $\BB$ in \ef{HermB}
 $$
 \tex{
  \psi_2(y)= \frac 1{\sqrt 2} \, F''(y).
   }
  $$


\ssk

For $k>1$, the solutions $Y(y)$ are oscillatory (changing sign) as
$y \to y_0^+$, so that expansions such as \ef{shtsol} must include
{\em oscillatory components} as an extra multiplier:
 \be
 \label{comp1}
Y(y)=(y-y_0)^{\hat \a}( \phi(s)+o(1)) \whereA s = \ln(y-y_0),
 \ee
 and $\phi(s)$ is a periodic solution of a $(2k+1)$th-order ODE.
 Examples of such oscillatory structures \ef{comp1} for $k=2$ are
 presented in \cite[pp.~186-192]{GSVR} for $k=2$ and $k=3$.

 We will use a similar asymptotic approach in Section \ref{S3.4} in the
 opposite limit $y \to +\iy$.



\section{On some mathematical aspects of similarity profiles}
 \label{S3}

\subsection{Local existence and uniqueness for $k=1$}

We need to show that the above numerical construction can be
justified, by carefully verifying some rigorous aspects of the
asymptotic analysis. We first apply a
 fixed point approach for the equivalent nonlinear integral equation
 to prove the expansion \ef{Yfbp} near the finite left-hand
 interface.

To this end, we look at our integrated second-order equation
\ef{2ndOrdNonl} for $k=1$, assuming that $Y(y)>0$ is {\em strictly
monotone increasing} (and hence non-oscillatory)
 sufficiently
close to the interface at $y=y_0^+<0$:
\begin{equation}
\label{nonl2nd}
\mbox{$Y^{\prime\prime}=-\frac{1}{n+3}\,|Y|^{-\frac{n}{n+1}}Y\,
y$} \quad \mbox{or} \quad
\mbox{$Y^{\prime\prime}=-\frac{1}{n+3}\,Y^{\frac{1}{n+1}}y$}
\,\,\,\, \mbox{for} \,\,\,\,
 Y>0.
\end{equation}
Therefore, we can rewrite our derivatives of $Y(y)$ in terms of
the inverse function $y(Y)$:
\begin{equation*}
\mbox{$Y^{\prime}=\frac{\mathrm{d}Y}{\mathrm{d}y}=\frac{1}{y^\prime(Y)}$}
\quad \mbox{and}\quad
\mbox{$Y^{\prime\prime}=\frac{\mathrm{d}}{\mathrm{d}y}\big(\frac{1}{y^\prime}\big)=-\frac{y^{\prime\prime}}{(y^\prime)^3}$}.
\end{equation*}
So, now we can reduce our differential equation \eqref{nonl2nd} to
its equivalent integral form,
\begin{equation}
\begin{split}
\label{baninit}  & \mbox{$-\frac{y^{\prime\prime}}{(y^\prime)^3}$}
= \mbox{$-\frac{1}{n+3}\,Y^{\frac{1}{n+1}}y$} \quad
\mbox{$\Longleftrightarrow \quad-\frac{1}{2(y^\prime)^2}$}
=\mbox{$-\frac{1}{n+3}$}
 \mbox{$
\int\limits_0^Y \mbox{$y(s)s^{\frac{1}{n+1}}$}\, \mathrm{d}s
 $}
\\
& \Longleftrightarrow \quad\quad y(Y)
 \mbox{$
=y_0 + \int\limits_0^Y \sqrt{\frac{n+3}{2\int\limits_0^r
y(s)s^{\frac{1}{n+1}}\,\mathrm{d}s}} \, \mathrm{d}r \equiv M(y).
 $}
\end{split}
\end{equation}

We next prove the following property of the integral operator $M$
in (\ref{baninit}):

\begin{proposition}
 \label{Pr.1}
For small $\d>0$, $M(y)$ is a contraction in $C[0,\delta]$, with
the sup-norm, and therefore admits a unique fixed point $y(Y)>0$
on $(0,\d)$ giving the unique positive solution of the ODE
 $(\ref{nonl2nd})$ on $(y_0,y_0+ \e)$ with some sufficiently small $\e=\e(\d)>0$.
\end{proposition}

\noindent {\em Proof.} We need to show that, for $M(y)$ to be a
contraction, in the $C$-metric,
\begin{equation}
\label{BanIneqdelta} \|M(\zeta_2)-M(\zeta_1)\|<\mu
\|\zeta_2-\zeta_1\|,
\end{equation}
for some constant $\mu=\mu(\d)\in(0,1)$. It is easy to see that
$M: Z_\delta \to Z_\delta$,
for the space $Z_\delta$ of continuous functions,
$Z_\delta = \{\zeta(Y)\in C[0,\delta], \, \zeta(0)=y_0\}$,
with the sup-norm:
\begin{equation*}
 \tex{
\|\zeta\|:= \sup_{Y\in(0,\delta)} |\zeta(Y)|.
  }
\end{equation*}

Now, take arbitrary $\zeta_1(Y), \zeta_2(Y) \in Z_\d$.  Then from
\eqref{baninit} we have that
\begin{equation*}
 \mbox{$
\mbox{$\|M(\zeta_2)-M(\zeta_1)\|=\sqrt{\frac{n+3}{2}}$}\int\limits_0^Y
\mbox{$\big\|\big(\int\zeta_2(s)s^{\frac{1}{n+1}}\,\mathrm{d}s\big)^{-\frac{1}{2}}-\big(\int
\zeta_1(s)s^{\frac{1}{n+1}}\,\mathrm{d}s\big)^{-\frac{1}{2}}\big\|$}
\,\mathrm{d}r,
 $}
\end{equation*}
where we use the simplified notation for the integral $\int$,
without any limits of integration, to mean $\int_0^r$.  This
equality can now be written as
\begin{equation*}
 \mbox{$
\mbox{$\|M(\zeta_2)-M(\zeta_1)\|=\sqrt{\frac{n+3}{2}}$}\int\limits_0^Y
\Big\|\frac{\big(\int\zeta_2(s)s^{\frac{1}{n+1}}\,\mathrm{d}s\big)^{-\frac{1}{2}}-\big(\int
\zeta_1(s)s^{\frac{1}{n+1}}\,\mathrm{d}s\big)^{-\frac{1}{2}}}{\big(\int\zeta_2(s)s^{\frac{1}{n+1}}\,\mathrm{d}s\big)^{-\frac{1}{2}}+\big(\int
\zeta_1(s)s^{\frac{1}{n+1}}\,\mathrm{d}s\big)^{-\frac{1}{2}}}\Big\|\,\mathrm{d}r.
 $}
\end{equation*}
Denoting the exponent $\nu=\frac{1}{n+1}$, we see that
\begin{equation*}
 \mbox{$
\mbox{$\|M(\zeta_2)-M(\zeta_1)\|\leq
\sqrt{\frac{n+3}{2}}$}\int\limits_0^Y \Big\|\frac{\int
[\zeta_2(s)- \zeta_1(s)]s^{\nu}\,\mathrm{d}s}{\int
\zeta_2(s)s^{\nu}\,\mathrm{d}s\int
\zeta_1(s)s^{\nu}\,\mathrm{d}s\big[\big(\int
\zeta_2(s)s^{\nu}\,\mathrm{d}s\big)^{-\frac{1}{2}}+\big(\int
\zeta_1(s)s^{\nu}\,\mathrm{d}s\big)^{-\frac{1}{2}}\big]}\Big\|\,\mathrm{d}r.
 $}
\end{equation*}
Since we deal with sufficiently small values of $Y$, so that
always $\zeta_{1,2}(s) \approx y_0$,
 it is easy to estimate in the denominator to get that
\begin{equation*}
\begin{split}
&\mbox{$\|M(\zeta_2)-M(\zeta_1)\|$}
 \mbox{$
\leq \mu_0\int\limits_0^Y \Big\|\frac{\|\zeta_2-\zeta_1\| \int
s^{\nu}\,\mathrm{d}s}{\int s^{\nu}\,\mathrm{d}s\int
s^{\nu}\,\mathrm{d}s\big[\big(\int
s^{\nu}\,\mathrm{d}s\big)^{-\frac{1}{2}}+\big(\int
s^{\nu}\,\mathrm{d}s\big)^{-\frac{1}{2}}\big]}\Big\|\,\mathrm{d}r
 $}
\\
&
 \mbox{$
\leq \mu_0 \|\zeta_2-\zeta_1\|\int\limits_0^Y
\frac{r^{\frac{n+2}{n+1}}}{r^{\frac{n+2}{n+1}}r^{\frac{n+2}{n+1}}r^{-\frac{n+2}{2(n+1)}}}\,
\mathrm{d}r
 $}
 \mbox{$
 \leq \mu_0 \|\zeta_2-\zeta_1\|\int\limits_0^Y r^{-\frac{n+2}{2(n+1)}}\,\mathrm{d}r
 $}
 \mbox{$
\leq \mu_0 \|\zeta_2-\zeta_1\|\,Y^{\frac n{2(n+1)}}.
 $}
\end{split}
\end{equation*}
Here, $\mu_0$ is a constant dependent on $n$ and $y_0$.  Since we
take $\frac{1}{2}\,|y_0|\leq y \leq |y_0|$, then we have that
$|Y|<1$, and so, fixing $Y\in[0,Y_0]$, with $\mu=\mu_0\,|Y_0|<1$,
we have that \eqref{BanIneqdelta} holds true.  Hence by {Banach's
Fixed Point Theorem} \cite[p.~39]{NonFuncAnaly}, $M(y)$ has a
unique fixed point in $Z_{\delta}$.
 $\qed$

\ssk

Similarly, we prove local existence  and uniqueness of a positive
solution for $l=1, \, 2$.

 For  $k \ge 2$, the
solutions are oscillatory (changing sign) close to interfaces, so
that a contraction approach does not apply, and, as we have
mentioned, we need to use other techniques of asymptotic analysis
in both limits $y \to y_0^+$ and $y \to +\iy$ (the latter one
includes $k=1$). Such oscillatory structures near interfaces have
been thoroughly studied for higher-order thin film equations; see
\cite{Gl4, Gl6}. Examples of such oscillatory patterns for various
NDEs can be found in \cite[Ch.~4]{GSVR}, so we do not address
these questions here anymore, and concentrate on another and more
difficult limit.


\subsection{Global existence and uniqueness for $k=1$}

We restrict our attention to first nonlinear eigenfunction $Y_0$
and prove the following:

\begin{theorem}
 \label{Th.1}
 For any $n>0$ and a
 fixed interface point $y_0<0$, the problem
 $\ef{2ndOrdNonl}$, $\ef{Yfbp}$ admits a unique solutions
 $Y_0(y)$, which is infinitely oscillatory as $y \to +\iy$.
  \end{theorem}

  \noi {\em Proof.} Once the local existence and uniqueness have
  been established in Proposition \ref{Pr.1}, its unique existence
  on the whole interval $(y_0,+\iy)$ follows from elementary
  checked local extension properties for the ODE \ef{2ndOrdNonl},
  which is shown not to admit strong singularities (``blow-ups")
  at any finite point. Oscillatory character of such a solution
  will be shown below. $\qed$


\subsection{Local and global behaviour of nonlinear eigenfunctions as $y \to +\iy$}

At this moment, we do not know the behaviour of nonlinear
eigenfunctions $f_l(y)$ (or $Y_l(y)$) for $y\gg 1$. Since the
first eigenfunction $Y_0(y)$ is expected to converge to our linear
kernel $F(y)$ as $n\to 0^+$ (the Airy function (\ref{Ai1}) for
$k=1$),  we also expect to have reasonably similar behaviour as $y
\to +\iy$. Recall that
 the rescaled kernel $F(y) \equiv \psi_0(y)$  has
 the following oscillatory slow algebraic decay as $y \to + \iy$
\cite[\S~2.2]{RayGI}:  for some constant $\hat c \in \re$,
\begin{equation}
\label{posasympt}
 \mbox{$F(y) \sim y^{-\frac{2k-1}{4k}}\cos
{\big(d_ky^{\frac{2k+1}{2k}}+\hat{c}\big)}$} \quad \text{as} \quad
y \to +\infty, \quad \mbox{where} \quad
\mbox{$d_k=2k\,\big(\frac{1}{2k+1}\big)^{\frac{2k+1}{2k}}$}.
\end{equation}
 Further eigenfunctions $\psi_l(y)$, given by the derivatives \ef{N4}, have
 the corresponding asymptotics via differentiating \ef{posasympt},
 so these are much more oscillatory and unbounded for $y \gg 1$.
  On the other hand, the eigenvalue problem \ef{N2} admits
  polynomial solutions (cf. \ef{poleigfunc} for $\BB^*$)
 \be
 \label{N7}
  \tex{
   \tilde \psi_l(y) \sim y^l +... \forA \tilde \l_l=
   \frac{l+1}{2k+1}, \quad l=0,1,2,...\,.
   }
   \ee
 Since $\tilde \psi_l \not \in L^2_\rho$, these are not proper
 eigenfunctions. However, this shows that the equation \ef{N2}
 admits polynomially growing non-oscillatory solutions as $y \to +\iy$.

 Due to the complicated nature of the nonlinear equations for $n>0$,
 it is not that easy to predict possible behaviour of solutions as
 $y \to + \iy$, where we expect them to be oscillatory as in the
 linear case $n=0$. We first study existence and nonexistence of
 non-oscillatory solutions, which mimic the polynomials \ef{N7}.
 Without loss of generality, we consider the ODE \ef{massconsY}
 for the first nonlinear eigenfunction $Y_0(y)$.

\begin{proposition}
 \label{Pr.Pol}
  For any $n>0$,
 the ODE $\eqref{massconsY}$ for even $k=2,4,...$
 admits  a bundle of algebraically growing solutions as $y \to + \iy$
  \be
  \label{N8}
   \tex{
   Y(y) \sim  \pm A y^m \whereA m= 2k + \frac{n+2}{n+1} \andA
   A=A(n,k) \ne 0,
   }
   \ee
   and does not admit such non-oscillatory solutions for odd
   $k=1,3,...\,$.
\end{proposition}

\noindent {\em Proof.} As a formal calculus, we substitute into
\ef{massconsY} the asymptotic expression from \ef{N8} to get by
balancing the leading terms:
 \be
 \label{N9}
 \tex{
 |A|^{-\frac n{n+1}} =(-1)^k m(m-1)...(m-2k+1)[(2k+1)+n],
 }
  \ee
  whence the explicit expression for $A$ and the result. $\qed$

  \ssk

Thus, for odd $k$, all the admitted behaviour as $y \to +\iy$ are
not of  algebraic growth, and, more plausibly, are oscillatory (as
well as for all even $k$). To show this, we separately consider
the particularly interesting case $k=1$.

\begin{proposition}
 \label{Pr.k1}
 All the solutions of the ODE $(\ref{massconsY})$ for $k=1$ are
 oscillatory as $y \to +\iy$, i.e., have infinitely many sign
 changes in any neighbourhood of $+\iy$.
  \end{proposition}

\noi {\em Proof.} Assume first that, for $k=1$,
\begin{equation*}
Y(y)\to+\infty \quad \text{as} \quad y\to+\infty.
\end{equation*}
Then \ef{massconsY} yields
 \be
 \label{N10}
 \tex{
  Y''= - \frac 1{n+3}\, y Y^{\frac 1{n+1}}\, \ll - \frac y{n+3} \forA
  y \gg 1 \LongA Y(y) \ll - \frac {y^3}{6(n+1)} \to -\iy,
  }
  \ee
  whence the contradiction. Further cases are studied similarly. The oscillatory character
 of all the solutions follows form the ODE in (\ref{N10}), since
  ${\rm sign} \,Y'' =- {\rm sign}\, Y$ (cf. $Y''=-Y$ for harmonic oscillations).
   $\qed$

  \ssk






\subsection{A priori bounds for  $Y(y)$ for $k=1$: nonlinear oscillatory tail}

We continue to study oscillatory properties of nonlinear
eigenfunctions $Y_l(y)$, as $y \to +\iy$, for the basic case $k=1$.
Let us fix two successive local extremum points  $1 \ll y_1<y_2$
of $Y(y)$, where $Y^\prime(y_1)=Y^\prime(y_2)=0$.
 Let us characterize the size of oscillations of $Y(y)$ at those
 points.

We initially look at the equation for the first nonlinear
eigenfunction, with the eigenvalue $\alpha_0(n)$, for $k=1$.
Multiplying \eqref{2ndOrdNonl} by $Y^\prime$ and integrating
over $(y_1, y_2)$ yields
\begin{equation*}
\mbox{$\int$}_{y_1}^{y_2} Y^{\prime\prime}Y^\prime =
\mbox{$-\frac{1}{n+3}$}\mbox{$ \int $}_{y_1}^{y_2} \,
\mbox{$|Y|^{-\frac{n}{n+1}}YY^\prime y$}.
\end{equation*}
Simplifying this, we see that
\begin{equation*}
\mbox{$\frac{1}{2}\,\big[(Y^\prime)^2\big]^{y_2}_{y_1}$} =
\mbox{$-\frac{(n+1)}{(n+2)(n+3)}$}\,\mbox{$ \int $}_{y_1}^{y_2}
\mbox{$(|Y|^{\frac{n+2}{n+1}})^\prime y$}.
\end{equation*}
After integrating by parts, we obtain
\begin{equation*}
\mbox{$\frac{1}{2}\,\big[(Y^\prime)^2\big]^{y_2}_{y_1}$} =
\mbox{$-\frac{(n+1)}{(n+2)(n+3)}\,\big[|Y|^{\frac{n+2}{n+1}}y\big]^{y_2}_{y_1}$}+
\mbox{$\frac{(n+1)}{(n+2)(n+3)}$}\,\mbox{$ \int $}_{y_1}^{y_2}
\mbox{$|Y|^{\frac{n+2}{n+1}}$}.
\end{equation*}
However, since we are looking at extremum points, where
$Y^\prime(y_1)=Y^\prime(y_2)=0$,
\begin{equation*}
\mbox{$\big[(Y^\prime)^2\big]^{y_2}_{y_1}$}=0
\quad \mbox{and} \quad
\mbox{$|Y(y_2)|^{\frac{n+2}{n+1}}y_2-|Y(y_1)|^{\frac{n+2}{n+1}}y_1
=$} \mbox{$ \int $}_{y_1}^{y_2} \mbox{$|Y|^{\frac{n+2}{n+1}} \,
\mathrm{d}y$}.
\end{equation*}
Since
 $\mbox{$ \int $}_{y_1}^{y_2} \mbox{$|Y|^{\frac{n+2}{n+1}} \,
\mathrm{d}y>0$}$,
we have
\begin{equation}
 \label{N11}
\mbox{$|Y(y_2)|^{\frac{n+2}{n+1}}y_2>|Y(y_1)|^{\frac{n+2}{n+1}}y_1$},
\quad \mbox{or, on rearranging,} \quad 
\mbox{$\big(\frac{|Y(y_2)|}{|Y(y_1)|}\big)^{\frac{n+2}{n+1}}>
\big( \frac{y_1}{y_2}\big)$}.
\end{equation}
This gives a lower estimate on the character of oscillations.


Let us now look at our second nonlinear eigenfunction with
$\alpha_1(n)$, where our equation is given by
\eqref{1stmomintcons}, $k=1$. Multiplying by $Y^\prime$ and
integrating yield
\begin{equation*}
\begin{split}
\mbox{$Y^\prime Y^{\prime\prime}y - (Y^\prime)^2 +\frac{1}{3+2n}
|Y|^{-\frac{n}{n+1}}YY^\prime y^2$}& = 0\\[1mm] \LongA
  \quad
\mbox{$\frac{1}{2}[(Y^\prime)^2]^\prime y - (Y^\prime)^2 +
\frac{n+1}{(n+2)(3+2n)}\,(|Y|^{\frac{n+2}{n+1}})^\prime y^2$}& =
0.
\end{split}
\end{equation*}
Integrating between $y_1$ and $y_2$ again, we have that
\begin{equation*}
\mbox{$\frac{1}{2}$}[(Y^\prime)^2y]_{y_1}^{y_2} -
\mbox{$\frac{1}{2}$}\mbox{$ \int $}_{y_1}^{y_2} (Y^\prime)^2
+\mbox{$\frac{n+1}{(n+2)(3+2n)}$}\,[|Y|^{\frac{n+2}{n+1}}y^2]_{y_1}^{y_2}
-\mbox{$\frac{2(n+1)}{(n+2)(3+2n)}$}\, \mbox{$ \int $}_{y_1}^{y_2}
|Y|^{\frac{n+2}{n+1}}y=0
 \end{equation*}
 \begin{equation*}
 \begin{split}
 \LongA \quad
-\mbox{$\frac{3}{2}$}\mbox{$ \int $}_{y_1}^{y_2} (Y^\prime)^2 -
\mbox{$\frac{2(n+1)}{(n+2)(3+2n)}$}\, \mbox{$ \int $}_{y_1}^{y_2}
|Y|^{\frac{n+2}{n+1}}y&\\[1mm]
+\mbox{$\frac{n+1}{(n+2)(3+2n)}$}&\,\mbox{$(|Y(y_2)|^\frac{n+2}{n+1}y_2^2-|Y(y_1)|^\frac{n+2}{n+1}y_1^2)$}=0.
\end{split}
\end{equation*}
One can see that
\begin{equation*}
\mbox{$\frac{3}{2}$}\mbox{$ \int $}_{y_1}^{y_2}
(Y^\prime)^2+\mbox{$\frac{2(n+1)}{(n+2)(3+2n)}$}\, \mbox{$ \int
$}_{y_1}^{y_2} |Y|^{\frac{n+2}{n+1}}y>0,
\end{equation*}
since $y>0$. Hence, we have  a similar estimate:
\begin{equation}
 \label{N12}
|Y(y_2)|^\frac{n+2}{n+1}y_2^2>|Y(y_1)|^\frac{n+2}{n+1}y_1^2,
\quad \mbox{or} \quad
\mbox{$\big(\frac{|Y(y_2)|}{|Y(y_1)|}\big)^{\frac{n+2}{n+1}}>\big(\frac{y_1}{y_2}\big)^2$}.
\end{equation}

Finally, for the third eigenfunction with $\alpha_2(n)$, using
\eqref{2ndmomintcons}, $k=1$ and multiplying by $Y^\prime$ yields
\begin{equation*}
\begin{split}
&\mbox{$Y^\prime Y^{\prime\prime} y^2 - 2(Y^\prime)^2 y +
2YY^\prime +\frac{1}{3+3n}Y^{-\frac{n}{n+1}}YY^\prime y^3$} =
0\\[1mm] \LongA \quad &\mbox{$\frac{1}{2}[(Y^\prime)^2]^\prime y^2
- 2(Y^\prime)^2 + (Y^2)^\prime +
\frac{n+1}{(n+2)(3+3n)}\,(|Y|^\frac{n+2}{n+1})^\prime y^3 = 0$}.
\end{split}
\end{equation*}
Integrating over $(y_1,y_2)$ yields
\begin{equation*}
\begin{split}
\mbox{$\frac{1}{2}$}[(Y^\prime)^2y^2]_{y_1}^{y_2} &-
 \mbox{$ \int $}_{y_1}^{y_2} (Y^\prime)^2 y -
2\mbox{$ \int $}_{y_1}^{y_2} (Y^\prime)^2 + [Y^2]_{y_1}^{y_2} +
\mbox{$\frac{n+1}{(n+2)(3+3n)}$}[|Y|^{\frac{n+2}{n+1}}y^3]_{y_1}^{y_2}\\[1mm]
&-\mbox{$\frac{3(n+1)}{(n+2)(3+3n)}$}\mbox{$ \int $}_{y_1}^{y_2}
|Y|^{\frac{n+2}{n+1}}y^2=0.
\end{split}
\end{equation*}
It then follows that
\begin{equation}
 \label{N13}
\mbox{$\big(\frac{|Y(y_2)|}{|Y(y_1)|}\big)^{\frac{n+2}{n+1}}>\big(\frac{y_1}{y_2}\big)^3$},
\quad \mbox{and, for any $l$,} \quad
\mbox{$\big(\frac{|Y(y_2)|}{|Y(y_1)|}\big)^{\frac{n+2}{n+1}}>\big(\frac{y_1}{y_2}\big)^{l+1}$}.
\end{equation}

Thus, the three estimates \ef{N11}, \ef{N12}, and \ef{N13}
characterize the behaviour of the nonlinear oscillatory tail for
$k=1$,
 which can be compared with numerical evidence in
Section \ref{NumConstrNonl}.

\subsection{More on oscillatory structure and periodicity}
 \label{S3.4}

As mentioned before, at least for small $n>0$, we expect that
nonlinear eigenfunctions exhibit, as $y \to +\iy$, a behaviour
which is structurally similar to the linear kernel and its
derivatives. We therefore expect to have a special oscillatory
behaviour for $y\gg1$, as we have already seen before. Hence, we
look to describe this oscillatory structure and, as a first
natural attempt, we will try to find if these oscillations are
given by {\em periodic functions}. Let us introduce the {\em
oscillatory component} $\phi(s)$ such that, as $y \to +\iy$,
\begin{equation}
\label{OscillSol}
Y(y)= y^\gamma \phi(s), \quad \text{where} \quad
s=\ln y,
\end{equation}
for some power $\gamma\in\mathbb{R}$.  Here the term $y^\gamma$
gives the rate of any growth/decay of the oscillations and may be compared to the controlling factor $y^{-\frac{2k-1}{4k}}$ found
in the linear asymptotics \ef{posasympt} for $n=0$.

Let us begin with the simpler case $k=1$, where substituting into
\eqref{nonl2nd}, we obtain
\begin{equation}
 \label{M1}
  \tex{
 (y^\g \phi)'' + \frac 1{n+3}\, y^{1+\frac \g{n+1}}|\phi|^{-\frac
 n{n+1}} \phi=0.
  }
 \ee
Expanding this expression and equating powers of $y$, we find that
\begin{equation}
 \label{M2}
  \tex{
\g-2= 1+\frac \g{n+1} \LongA \mbox{$\gamma = \frac{3(n+1)}{n}$}.
 }
\end{equation}
This gives us a second-order ODE for the oscillatory component:
\begin{equation}
\label{wrongoscilleqn}
 \tex{
 P_2(\phi) \equiv \phi''+(2\g-1) \phi'+ \g(\g-1) \phi=
- \frac 1{n+3}\, |\phi|^{-\frac
 n{n+1}} \phi \inB \re.
  }
 \ee
A most typical orbit describing oscillations should be a periodic
orbit of \ef{wrongoscilleqn}. However, since $\g>0$ in \ef{M2},
the behaviour \ef{OscillSol} describes unbounded oscillations as
$y \to +\iy$, which are not acceptable for a source-type solution;
cf. the linear decaying one \ef{posasympt} for $n=0$.
Hence, we conclude that  oscillatory behaviour as $y \to +\iy$ is
not given by periodic oscillatory components $\phi(\ln y)$ as in
\ef{OscillSol}.

Similarly, we arrive at a contradiction, applying \ef{OscillSol}
to the general equation \ef{D11}:
\begin{equation*}
 \tex{
\mbox{$(-1)^{k}D^{2k+1}_y(y^\gamma\phi) =
\frac{1}{(2k+1)+(l+1)n}\,y
(y^{\frac{\gamma}{n+1}}|\phi|^{-\frac{n}{n+1}} \phi)^\prime$} +
\frac{l+1}{(2k+1)+(l+1)n}\,
y^{\frac{\gamma}{n+1}}|\phi|^{-\frac{n}{n+1}}\phi.
 }
\end{equation*}
Balancing the polynomial terms yields
 \be
 \label{M4}
 \tex{
 \g-(2k+1)= \frac \g{n+1} \LongA \g=\frac{(2k+1)(n+1)}n>0.
  }
  \ee
 This leaves us with an ODE of the order $2k+1$, for $\phi(s)$:
\begin{equation}
 \label{M5}
\mbox{$
(-1)^k
P_{2k+1}(\phi)=\frac{1}{(2k+1)+(l+1)n}\,(|\phi|^{-\frac{n}{n+1}}\phi)^\prime
+ \frac{1}{n}\,|\phi|^{-\frac{n}{n+1}} \phi.
 $}
\end{equation}
Here, $P_{2k+1}(\phi)$ is a polynomial operator on $\phi$, induced
by the term
 $ 
D^{2k+1}_y(y^\gamma\phi).
 $ 
For the case $k=1$, the polynomial operator is given by
\begin{equation*}
P_3(\phi) = \mbox{$\phi^{\prime\prime\prime} +
\frac{3(2n+3)}{n}\,\phi^{\prime\prime}+ \frac{9n^2 + 7n +27}{n}\,
\phi^{\prime} + \frac{3(n+1)(2n+3)(n+3)}{n^3}\, \phi$}.
\end{equation*}
In general, $P_{2k+1}(\phi)$ is defined by the recursion
\begin{equation*}
\begin{split}
P_{2k+1}(\phi) = \mbox{$
\frac{\mathrm{d}^2}{\mathrm{d}s^2}P_{2k-1}(\phi)$}&+
\mbox{$(\gamma-2k)\frac{\mathrm{d}}{\mathrm{d}s}P_{2k-1}(\phi)
+(\gamma -
2k+1)\frac{\mathrm{d}}{\mathrm{d}s}P_{2k-1}(\phi)$}\\[1mm]
&\mbox{$+ (\gamma-2k)(\gamma-2k+1)P_{2k-1}(\phi)$}.
\end{split}
\end{equation*}


As usual, periodic solutions of \ef{M5} are at most simple
oscillatory components. However, as before, since $\g>0$ in
\ef{M4}, oscillatory structures of the form \eqref{OscillSol}, for
any $k \ge 2$, are not applicable, at least for the first
nonlinear eigenfunction $Y_0(y)$, which is assumed to be
integrable as $y \to +\iy$.

For higher-order nonlinear eigenfunctions $Y_l(y)$, with $l \ge 1$,
proving existence of periodic solutions of the ODE \ef{M5} is the
first step in understanding the oscillatory behaviour. This is a
difficult mathematical problem, which nevertheless can be solved for
orders of $k$ that are not too large. We refer to \cite{Gl4, Gl6,
PetI} for key references and recent results on existence-uniqueness
of periodic solutions of even-order ODEs such as \ef{M5}, which
occur in parabolic thin film theory. We also refer to
\cite[\S~4]{GalPMEn} for existence results of periodic orbits for
oscillatory solutions for the PME--4 \ef{PME2m}, $m=2$.

Overall, we rule out the ``periodic" structures \ef{OscillSol} as
$y \to +\iy$ for the first nonlinear eigenfunction $Y_0(y)$, which
is expected to have an oscillatory decay and be integrable (not in
the absolute sense, i.e., it is not measurable there). Therefore,
in this case, the oscillatory behaviour may be more complicated
and corresponds to not that easy ``nonlinear focus", which we are
going to catch using numerical methods.
 For $l \ge 1$,
 such a behaviour \ef{OscillSol} with almost periodic oscillatory components $\phi(\ln y)$
 is still plausible (but remains rather suspicious).

\subsection{Branching of nonlinear eigenfunctions at $n=0$}
 \label{S3.5}

As we have promised, we now apply another classic idea to trace
out the behaviour of all the nonlinear eigenfunctions for small
$n>0$. Namely, we are going to show that
 there exists branching at $n=0^+$ of solutions with respect to the parameter
$n$.
 In other words, we show that,  as $n\to 0$, there
exists certain convergence to solutions (driven by the
eigenfunctions of the linear operator $\BB$ in \ef{HermB}) of the
the LDE \ef{airyfunc}.

\ssk

 To this end,
 let us look at the general ODE given by \eqref{rednonlOde}.
We first expand $|f|^n$ to formally get
\begin{equation}
 \label{fn0}
\mbox{$|f|^n=1+n\ln|f|+O(n^2)$}.
\end{equation}
This is pointwise and uniformly true in any bounded positivity
subset $\{|f|\geq\delta_0>0\}$.
 However, we are not
 at this moment going to discuss a rigorous functional meaning of
this expansion for changing sign functions $f(y)$ defined in the
whole $\mathbb{R}$.
 Note that \ef{fn0} can then  be understood in a
 weak sense, which may be  sufficed  for passing to the limit in the
 equivalent integral equations; see \cite[\S~7.6]{Gl4}
 for asymptotic details.

Thus, using the formal expansion \ef{fn0}, \eqref{rednonlOde}
reduces to
\begin{equation*}
\mbox{$(-1)^{k+1}D_y^{2k+1}\big[(1+n\ln|f|)f\big] + \frac{1-\alpha
n}{2k+1}\,f^{\prime}y +\alpha f+O(n^2)=0$}.
\end{equation*}
Expanding coefficients for small $n>0$ yields
\begin{equation}
\label{branorig} \mbox{$({\bf B}-\lambda_l I)f +
(-1)^{k+1}D_y^{2k+1}(n\ln|f|f)+(\alpha-\frac{l+1}{2k+1})\,f-\frac{\alpha
n}{2k+1}\,f^\prime y+O(n^2)=0$},
\end{equation}
where $\BB$ is the linear operator \ef{HermB} and
$\lambda_l=-\frac{l}{2k+1}$ is its $(l+1)$th eigenvalue.

For $l<2k+1$, we can find our eigenvalues $\alpha_l(n)$ explicitly
as in \ef{alk}, so that
\begin{equation*}
\alpha_l(n)=\mbox{$\frac{l+1}{(2k+1)+n(l+1)}$}
=\mbox{$\frac{l+1}{2k+1}\,\big[1+\frac{n(l+1)}{2k+1}\big]^{-1}$}
=\mbox{$\frac{l+1}{2k+1}\,\big[1-\frac{n(l+1)}{2k+1}\big]
+O(n^2)$}.
\end{equation*}
Then \eqref{branorig} reduces to
\begin{equation*}
\mbox{$({\bf B}-\lambda_l I)f +
(-1)^{k+1}D_y^{2k+1}(n\ln|f|f)-\frac{n(l+1)^2}{(2k+1)^2}\,f-\frac{n(l+1)}{(2k+1)^2}\,f^\prime
y+O(n^2) = 0$}.
\end{equation*}

Hence using the Lyapunov--Schmidt method \cite{VainbergTr} by
setting
\begin{equation}
\label{LSbranchred} f=\psi_l +n\phi_l +O(n^2),
\end{equation}
we obtain, within the order $O(n)$, the following inhomogeneous
equation:
\begin{equation}
 \label{ss1}
\mbox{$({\bf B}-\lambda_l I)\phi_l =
(-1)^{k}D_y^{2k+1}(\ln|\psi_l|\psi_l)+\frac{(l+1)^2}{(2k+1)^2}\,\psi_l
+\frac{(l+1)}{(2k+1)^2}\,\psi_l^\prime y \equiv h$}.
\end{equation}

Using  Hermitian spectral theory for $\BB$ and
completeness-closure of the eigenfunctions subset
$\Phi=\{\psi_l\}_{l\ge 0}$ \cite[\S~4]{RayGI}, for the unique
solvability of (\ref{ss1}) for $\phi_l$, it now remains to demand
that the right-hand side $h$ is orthogonal to $\psi_l^\ast$, i.e.,
 \be
 \label{ss2}
 \langle h ,\psi_l^\ast \rangle_\ast = 0.
   \ee
   Here we have to
use the corresponding indefinite metric $\langle \cdot, \cdot
\rangle_*$, in which the pair $\{\BB,\BB^*\}$ comprises the
operator $\BB$ and its adjoint (this metric can be reduced to the
standard dual $L^2$-one, \cite[\S~5]{RayGI}). Here, (\ref{ss2}) is
known as a scalar bifurcation equation in the classic
Lyapunov-Schmidt method \cite{VainbergTr}.

We then use the adjoint polynomial eigenfunctions $\psi_l^\ast$
given by \eqref{poleigfunc}.
  Then, for $l<2k+1$, we have that the generalized Hermite polynomials are simple
  \cite[\S~5]{RayGI},
$$
 \mbox{$
 \psi_l^\ast(y) = \mbox{$\frac{1}{\sqrt{l!}}\,y^l$}
  \quad (0 \le l
   < 2k+1).
  $}
  $$
Hence, for all $l<2k+1$, (\ref{ss2}) is indeed valid:
\begin{equation*}
\begin{split}
\langle h ,\psi_l^\ast \rangle_\ast &=
\mbox{$\frac{1}{\sqrt{l!}}$}\, \mbox{$ \int $}
\mbox{$\big[(-1)^{k}D_y^{2k+1}(\ln|\psi_l|\psi_l)y^l
+(-1)^l\,\frac{(l+1)^2}{(2k+1)^2}\,\psi_ly^l+(-1)^l\,\frac{(l+1)}{(2k+1)^2}\,\psi_l^\prime
y^{l+1}\, \big] \mathrm{d}y $}\\[1mm] &=
\mbox{$\frac{1}{\sqrt{l!}}$}\, \mbox{$ \int $} \mbox{$ \big[
(-1)^{k}D_y^{2k+1}(\ln|\psi_l|\psi_l)y^l +
(-1)^l\,\frac{(l+1)}{(2k+1)^2}\,(\psi_ly^{l+1})^\prime \big]$}\,
\mathrm{d}y=0.
\end{split}
\end{equation*}

 Recall that, for $l\geq 2k+1$, we do not know  nonlinear eigenvalues $\alpha_l(n)$
explicitly.  In this case, we expand $\a_l(n)$ as follows:
\begin{equation}
\label{alphan} \alpha_l(n)=\alpha_0+\alpha_1 n+O(n^2),
\end{equation}
 where $\a_0=\a_l(0)$ in \ef{N3} comes from linear Hermitian
 theory, and $\a_1$ is a new unknown. As before, we use \eqref{branorig} and now we substitute
\eqref{alphan}, as well as \eqref{LSbranchred}, to obtain
\begin{equation*}
\begin{split}
\mbox{$n({\bf B}-\lambda_l I)\phi_l=$}&
\mbox{$(-1)^{k}D_y^{2k+1}(n\ln|\psi_l|\psi_l)$}+
\mbox{$\big(\alpha_0+n\alpha_1-\frac{l+1}{2k+1}\big)\,\psi_l$}\\[1mm]&+
\mbox{$n\big(\alpha_0-\frac{l+1}{2k+1}\big)\,\phi_l$}
\mbox{$-\frac{n\alpha_0}{2k+1}\,\psi^\prime_ly+O(n^2)=0$}.
\end{split}
\end{equation*}
Equating as usual the terms of the order $O(n)$, we can find the
value of $\alpha_0$, with
\begin{equation}
\label{initialpha1} \begin{split}
 &\mbox{$\alpha_0-\frac{l+1}{2k+1}=0$}, \quad \mbox{and} \quad
 \mbox{$n({\bf B}-\lambda_l I)\phi_l $}
 \\
& =
 \mbox{$
 (-1)^{k}D_y^{2k+1}(n\ln|\psi_l|\psi_l)+ n\alpha_1\psi_l+
n\big(\alpha_0-\frac{l+1}{2k+1}\big)\,\phi_l-\frac{n\alpha_0}{2k+1}\,\psi^\prime_ly
$}.
 \end{split}
\end{equation}
Hence,  substituting into \eqref{initialpha1} and passing to  the
limit $n \to 0^+$  yield
\begin{equation*}
\mbox{$({\bf B}-\lambda_l I)\phi_l=
(-1)^{k}D_y^{2k+1}(\ln|\psi_l|\psi_l)+
\alpha_1\psi_l-\frac{l+1}{(2k+1)^2}\,\psi^\prime_ly \equiv h$}.
\end{equation*}

Then, the orthogonality condition (\ref{ss2}) becomes an algebraic
equation for $\a_1$ in (\ref{alphan}). Namely, taking the inner
product with $\psi_l^\ast$ and noting that $\langle h, \psi_l^\ast
\rangle_\ast=0$, $\langle \psi_l, \psi_l^\ast\rangle_\ast = 1$
yield
\begin{equation}
 \label{112}
\mbox{$\alpha_1=-\langle
(-1)^{k}D_y^{2k+1}(\ln|\psi_l|\psi_l)-\frac{l+1}{(2k+1)^2}\,\psi^\prime_ly
, \psi_l^\ast \rangle_\ast$}.
\end{equation}

 Thus,
this uniquely defines the second coefficient  $\a_1$ in the
expansion \ef{alphan} and then, as usual, \ef{112} gives a unique
function $\psi_l$ in the eigenfunction expansion \ef{LSbranchred},
etc.

In the analytic or even finite regularity cases, solvability
conditions and existence of expansions such as \ef{alphan} usually
rigorously justify the actual presence of branching. Our case is
more delicate in view of the ``weakness" of the expansion
\ef{fn0}.
 However, for the variable $Y=|f|^nf$, the expansion (\ref{fn0})
 is easier to justify, especially now the equation becomes
 semilinear and can be reduced to an integral equation with
 compact Hammerstein--Uryson-type operators. Therefore, in the present case,
   a rather full justification of the
$n$-branching method, though being rather technical, is
 doable and does not represent a principally non-solvable problem of nonlinear
 integral operator theory.





\section{The nonlinear limit $n\to\infty$: $k=1$ and $l=0$}
\label{S4}

While we dealt before with a ``homotopy path" construction of
nonlinear eigenfunctions in the limit of small $n \to 0^+$, we now
consider the opposite ``highly nonlinear" limit $n \to +\iy$,
which also helps to understand properties of the nonlinear
eigenvalue problem.

\subsection{Reducing to an algebraic problem}

As we have seen in Section \ref{S3.5},  as $n\to 0$,
 there appears a direct connection of all the nonlinear eigenfunctions
 with the linear ones associated with the LDE \ef{airyfunc}.
 Now, we are going to perform the opposite {\em highly nonlinear limit}
$n \to + \infty$ for nonlinear eigenfunctions of the NDE.
Rather surprisingly, it turns out that this ``limit nonlinear case"
 admits a more profound analysis, and for $l=0,1,2$ we are able to
 tackle  the nonlinear eigenvalue problems by a simpler geometric-algebraic,
allowing us to obtain a number of analytical and explicit
expressions.

  Consider the ODE \eqref{massconsY},
with $k=1$ and where $l=0$, i.e., \eqref{2ndOrdNonl}. Since we are
dealing with  $n\gg 1$, it is necessary to scale out any
coefficients containing large $n$'s.  In order to do this, we let
\begin{equation}
 \label{C11}
Y(y) = C \, \tilde{Y}(y) \whereA C=C(n)>0 \quad \mbox{is a
constant.}
\end{equation}
 Substituting \ef{C11} into
\eqref{2ndOrdNonl} yields:
\begin{equation}
\label{C11N}
 \mbox{$C\tilde{Y}^{\prime\prime} =
-\frac{1}{n+3}\,|C|^{-\frac{n}{n+1}}C|\tilde{Y}|^{-\frac{n}{n+1}}\tilde{Y}y$}
 \LongA  C = (n+3)^{-\frac{n+1}{n}},
\end{equation}
so that, on scaling out such a $C(n)$, we obtain the rescaled ODE
\begin{equation}
 \label{nInf}
\mbox{$\tilde{Y}^{\prime\prime}=-
|\tilde{Y}|^{-\frac{n}{n+1}}\tilde{Y}y$}.
\end{equation}
 Here we can pass to the limit $n \to + \iy$, since
 the only $n$-dependent exponent satisfies
 $-\frac{n}{n+1}\to-1$. At $n=+\iy$, the ODE becomes simpler
 and contains a bounded discontinuous nonlinearity:
\begin{equation}
  \fbox{$
\label{ninftyeqn}
  \BB_\iy(\tilde Y) \equiv
 \mbox{$\tilde{Y}^{\prime\prime}$} +
{\rm sign}\, \tilde{Y} \, y=0 \inB \re.
   $}
\end{equation}
 Of course, passing to the limit to arrive at \ef{ninftyeqn} might be a
 delicate mathematical problem. A potentially dangerous situation
 occurs in those subsets, where $Y$ vanishes. However, if both
 $Y(y)$ and $\tilde Y(y)$ have a.a. zeros transversal in the
 natural sense, then passage to the limit is straightforward.
 A sufficient ``transversality" of zeros of the limit function
 $\tilde Y(y)$ can be checked {\em a posteriori}, after completing
 our algebraic construction.

Solving ``almost linear" ODE \ef{ninftyeqn}, we find  expressions
dependent on the sign of $\tilde{Y}$:
\begin{equation}
\label{casesninfty}
\begin{cases}
\tilde{Y}>0: \quad \mbox{$\tilde{Y}_+(y) =
-\frac{1}{6}\,y^3+c_1y+c_2$},\\[1mm]
\tilde{Y}<0: \quad \mbox{$\tilde{Y}_-(y) =
\frac{1}{6}\,y^3+d_1y+d_2$}.
\end{cases}
\end{equation}
Here $c_1, \,c_2,\, d_1,\, d_2$ are all constants,  not
necessarily positive ones. Knowing the conditions of continuity
for the function $\tilde{Y}(y)$ and $\tilde Y'(y)$, we must have
that all one-sided limits coincide, i.e.,
 \be
 \label{B2}
 \tilde{Y}_+=\tilde{Y}_- \andA
\tilde{Y}^{\prime}_+=\tilde{Y}^{\prime}_-, \quad \mbox{at any
zero, where}
 \quad
\tilde{Y}=0.
 \ee

Let the points $\{y=y_i\}_{i \ge 0}$ be successive  zeros, i.e.,
$\tilde{Y}(y_i)=0$, for $i=0,1,2,...\,$.  Hence, for
$\tilde{Y}_+(y_i)$ given in \eqref{casesninfty}, for a fixed
isolated zero  with an  $i \ge 1$ (as usual, $y_0<0$ corresponds
to the left-hand interface),
\begin{equation*}
\mbox{$\tilde{Y}_+(y_i)=-\frac{1}{6}\,y^3_i + c_{1i}y_i + c_{2i} =
0$}.
\end{equation*}
We now rearrange this to find one of the unknown parameters, in
terms of $y_i$ and the parameter $c_{2i}$, such that
$\mbox{$c_{1i} = \frac{1}{6}\,y^2_i - \frac{c_{2i}}{y_i}$}$.
We also have that $\tilde{Y}^\prime_+(y_i)=\tilde{Y}^\prime_-(y_i)$,
hence
\begin{equation*}
\mbox{$-\frac{1}{2}\,y^2_i + c_{1i} = \frac{1}{2}\,y^2_i +
d_{1i}$}.
\end{equation*}
From this, we find a second parameter in terms of $y_i$ and
$c_{2i}$, where
$\mbox{$d_{1i} = -\frac{5}{6}\,y^2_i - \frac{c_{2i}}{y_i}$}$.
Now, we see that from $\tilde{Y}_+=\tilde{Y}_-$,
\begin{equation*}
\mbox{$c_{1i}y_i+d_{1i}y_i+c_{2i}+d_{2i} = 0$}.
\end{equation*}
So, substituting in known values, we have our third parameter
$d_{2i}$ given by
 $\mbox{$d_{2i} = \frac{2}{3}y_i^3 + c_{2i}$}$.

We now see that, after substituting values for $c_{1i}, d_{1i}$,
and $d_{2i}$, \eqref{casesninfty} can be written as
\begin{equation}
\label{ninftycasesconsts}
\begin{cases}
\tilde{Y}>0: \quad \mbox{$\tilde{Y}_+(y) =
-\frac{1}{6}\,y^3+(\frac{1}{6}\,y^2_i -
\frac{c_{2i}}{y_i})y+c_{2i}$},\\[1mm] \tilde{Y}<0: \quad
\mbox{$\tilde{Y}_-(y) = \frac{1}{6}\,y^3-(\frac{5}{6}\,y^2_i +
\frac{c_{2i}}{y_i})y+\frac{2}{3}\,y_i^3 + c_{2i}$}.
\end{cases}
\end{equation}
From the above, it is noted that $y_i\neq 0$, for any $i$, unless
$c_{2i}=0$.

\ssk

 \subsection{Existence, uniqueness, and zero properties of $\tilde Y_0(y)$}

We now resolve the algebraic system to get a rather complete
description of some important properties of the solution $\tilde
Y_0(y)$ obtained by such a ``geometric approach.
Recall first  the scaling invariance of the ODE \ef{ninftyeqn}:
 \be
 \label{tr2}
  \tex{
\tilde Y(y) \,\,\,\mbox{is a solution} \LongA  \pm \, a^3 \tilde
Y\big(\frac y a \big)
 \,\,\, \mbox{is a solution for any $a>0$}.
  }
  \ee
 Therefore, choosing the interface at $y_0=-1$, we put two
 conditions there
  \be
  \label{tr00}
  \tilde Y(-1)=\tilde Y'(-1)=0,
   \ee
   and prove  the following:

   \begin{theorem}
   \label{Th.2}
   The problem $\ef{ninftyeqn}$, $\ef{tr00}$ admits a unique
   nontrivial solution $\tilde Y_0(y)$, which has transversal
   zeros $\{y_i\}_{i \ge 1}$ such that
    \be
    \label{tr01}
     \begin{matrix}
    y_0=-1, \quad y_1=2|y_0|=2, \quad y_2=3 \sqrt 2-1=3.2426...\, ,
     \andA \ssk \\
y_{i+1}= \frac{\sqrt{17 y_i^2-4 y_i y_{i-1}-4 y_{i-1}^2}-y_i}2
     \quad \mbox{for any} \quad i \ge 2.
 \end{matrix} 
      \ee
 \end{theorem}

\noi{\em Proof.} As we have promised, we  construct such a
solution using pure algebraic manipulations. Let us begin with the
first interval of positivity $(y_0=-1,y_1)$, where according to
\ef{tr00}, the solution reads
 \be
 \label{tr10}
  \tex{
 \tilde Y_+(y)=- \frac 16 \,(y+1)^2 \big(y-y_1\big).
  }
  \ee
Since \ef{casesninfty} implies no quadratic term $\sim y^2$ in the
cubic polynomial, this uniquely gives $y_1=-2 y_0=2$, and hence
 \be
 \label{tr101}
  \tex{
 \tilde Y_+(y)=-\frac 16\,(y+1)^2 \big(y-2\big)>0 \quad \mbox{on} \quad (-1,2).
  }
  \ee

On the next interval $(y_1=2,y_2)$, the negative solution takes
the form:
 \be
 \label{tr1}
  \tex{
\tilde  Y_-(y)=\frac 16\, (y-2)(y-y_{2})\big(y+ c_{1}\big)<0.
  }
  \ee
Similarly to the above, we then conclude that
 \be
 \label{tr11}
 c_{1}=2+y_2.
 \ee
 Then matching of the first derivatives \ef{B2} at $y=y_1=2$ implies
  \be
  \label{tr12}
   \tex{
    -\frac 16\, 3^2= \frac 16 \,(2-y_2)(4+y_2)
    \LongA y_2^2+ 2 y_2-17=0 \LongA y_2=3 \sqrt 2-1,
     }
     \ee
 and this procedure can be continued.

\ssk

Consider an arbitrary interval $(y_{i-1},y_i)$ of positivity (or
negativity), $i \ge 2$, with
 \be
 \label{tr13}
  \tex{
   \tilde Y_+(y)=(y-y_{i-1})(y-y_i)\big(- \frac 16\, y+
   c_{i-1}\big) \whereA c_{i-1}=- \frac 16\,(y_i+y_{i-1}).
   }
   \ee
Similarly, on the next negativity (or resp. positivity) interval
$(y_i,y_{i+1})$,
 \be
 \label{tr14}
\tex{
   \tilde Y_-(y)=(y-y_{i})(y-y_{i+1})\big(\frac 16\, y+
   c_{i}\big) \whereA c_{i}= \frac 16\,(y_i+y_{i+1}).
   }
   \ee
Therefore, the matching at $y=y_i$ yields
 \be
 \label{tr15}
  \begin{matrix}
  (y_i-y_{i-1})\big( - \frac 16 \, y_i- \frac 16\,(y_i+y_{i-1})\big)=
(y_i-y_{i+1})\big(  \frac 16 \, y_i+ \frac 16\,(y_i+y_{i+1})\big)
 \qquad\qquad
 \ssk\ssk\\
 \LongA  y_{i+1}^2+ y_i y_{i+1}- 4 y_i^2 + y_i
 y_{i-1}+y_{i-1}^2=0, \qquad\qquad
\end{matrix} 
 \ee
 whence the final result in \ef{tr01}.
 Thus, the unique solution can be extended indefinitely for arbitrary $y \gg 1$ and
  is infinitely oscillatory
 as $y \to + \iy$. $\qed$

 \ssk

\noi{\bf Remark.} The quadratic equation in \ef{tr15} reduces to a
2D {\em linear discrete equation}:
 \be
 \label{alg1}
 \a_{i,j} \equiv y_i y_j \LongA
 \a_{i+1,i+1}+\a_{i,i+1}-4 \a_{i,i}+\a_{i,i-1}+ \a_{i-1,i-1}=0
 \forA i \ge 2.
  \ee
It is not difficult to find some particular solutions:
 \be
 \label{alg2}
 \tex{ \a_{i,j}= \mu^i \nu^j \LongA
 \mu(\mu+1)\nu^2-4 \mu \nu + \mu+1 \LongA
 \nu_\pm(\mu)= \frac{2\mu\pm (\mu-1)\sqrt{-\mu}}{\mu(\mu+1)}.
 }
  \ee
Therefore, denoting by ${\mathcal M}$ a proper subset of
parameters $\mu$ such that $\{\mu^i\}_{\mu \in{\mathcal M}}$ is
complete/closed,
 the
general solution of \ef{alg2} is represented by a converging
infinite series
 \be
 \label{alg3}
  \tex{
  \a_{i,j}= \sum_{\mu \in {\mathcal M}} C_\mu \mu^i \nu^j,
  }
  \ee
 where $\{C_\mu\}$ are constants and $\nu=\nu(\mu)$ take values
 according to \ef{alg2}. Here, \ef{alg3} is a discrete analogy
 of  eigenfunction expansions of solutions of a linear PDE with two
 independent variables $(x,t)$. A proper posing ``boundary
 conditions" for \ef{alg1} to specify the corresponding
 Sturm--Liouville problem and next initial conditions to get the
 corresponding eigenfunction expansion \ef{alg3} is a difficult
 and uncertain problem\footnote{It seems, a suitable behaviour at
 infinity, as $i,j \to +\iy$, of $\a_{i,j}$ might include something
 like the ``minimality" condition in \ef{BC1}, which is hard to
 take into account.}, which will unlikely provide us with any useful
 finite explicit formulae.

\ssk


However, the analytic relationships such as \ef{tr01} can promise
 to get extra asymptotic properties of the first rescaled eigenfunction
 $\tilde Y_0(y)$, especially the decay rate of the minimal
 behaviour indicated in \ef{BC1}. However, in the present case
 $l=0$, unlike the simpler one $l=2$ in Section \ref{S42}, some
 computations are not expected to be easy all the way, since such
 algebraic relations are quadratic and hence not always explicitly
 solvable.

\subsection{Numerics for $n  \gg 1$: $k=1$ and $l=0$}

As usual, very sharp proper numerics can help to detect further
properties of $\tilde Y_0(y)$, and hence avoid trying to get too
complicated and exhaustive  results concerning  the corresponding
algebraic system.
 Moreover, which is even more important, we can
also check the character of convergence of solutions as $n \to +
\iy$, which, after scaling \ef{C11}, turns out to be rather fast.
 To deal with the ODE
\eqref{ninftyeqn}, it is possible to just use a simple shooting
method using the ODE solver {\tt ode45}, to find suitable profiles
and nonlinear eigenfunctions.





We use a similar shooting method as that applied for the general
case of $n>0$, set out in Section \ref{NumConstrNonl}.  We recall
that, close to the interface at some $y=y_0<0$, we look for small
solutions of $\tilde Y(y)$ such that, as $y \to y_0^+$,
\begin{equation}
 \label{ninftyshtY}
  \begin{matrix}
\tilde Y(y)=C_0(y-y_0)_+^{2}\big(1+o(1)\big) \whereA C_0= \frac 12
\, |y_0|>0,\ssk\ssk\\
 \quad \mbox{and} \quad
\mbox{$\tilde{Y}^\prime(y) = |y_0|(y-y_0)(1+o(1))$} \quad(k=1).
\end{matrix}
\end{equation}
  The proof of this
expansion is similar and even simpler than that of Proposition
\ref{Pr.1}.

In Figure \ref{FF1}, we show the general view of the first
eigenfunction $\tilde Y_0(y)$ on the large interval
$[y_0=-1,100]$. The envelope of the decaying oscillations are
governed by the algebraic curve
 \be
 \label{L0}
  \tex{
  L_0(y) \approx \pm 0.7 \, y^{-\frac 13} \asA y \to + \iy.
   }
   \ee
   It seems that, definitely, this decay can be seen from the
  above  algebraic system. We have checked that, for  $n=100$,
the corresponding nonlinear eigenfucntion (after scaling) is
practically indistinguishable.


 \begin{figure}
\centering
\includegraphics[scale=0.85]{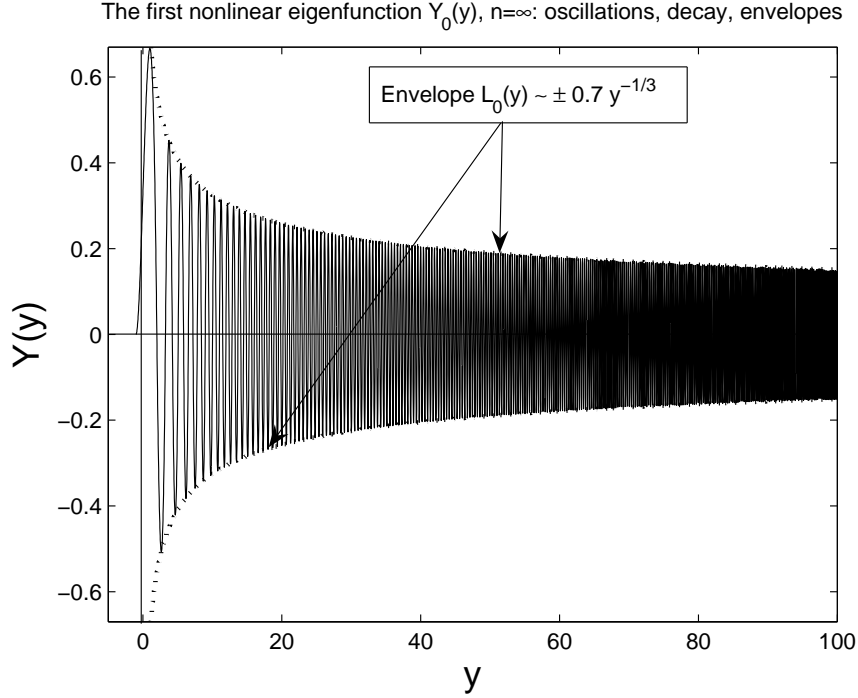}  
\vskip -.3cm
  \caption{The first nonlinear eigenfunction
 $\tilde Y_0(y)$ as a solution of the problem \ef{ninftyeqn}, \ef{ninftyshtY};
 $k=1$, $l=0$.}
 \label{FF1}
\end{figure}

In Figure \ref{FF2}, we show $\tilde Y_0(y)$ on a smaller interval
$y \in [-1,7]$, indicating all the first zeros, which will coincide
with those given by explicit algebraic expressions in \ef{tr01}.


 \begin{figure}
\centering
\includegraphics[scale=0.85]{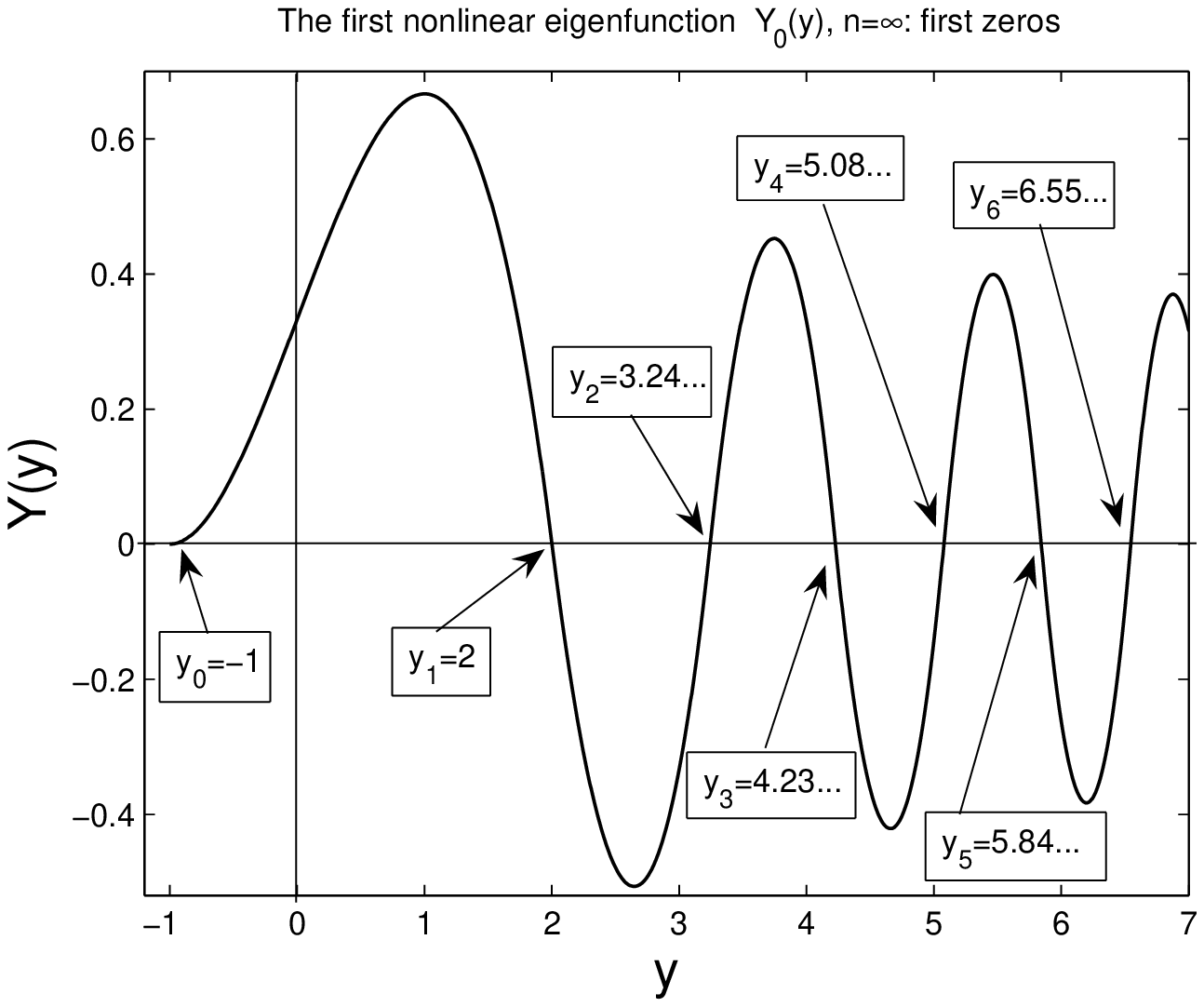}  
\vskip -.3cm
  \caption{The first nonlinear eigenfunction
 $\tilde Y_0(y)$ as a solution of the problem \ef{ninftyeqn}, \ef{ninftyshtY};
 $k=1$, $l=0$.}
 \label{FF2}
\end{figure}

All the computations have been performed with the enhanced
accuracy and tolerances $\sim 10^{-10}$, so these are quite
reliable. Let us present  first 15 zeros of $\tilde Y_0(y)$:
 $$
  \begin{matrix}
  y_0=-1, \,\, y_1=2|y_0|=2, \,\, y_3=3.2426...\, , \,\,
 y_4=5.0777...\, , \,\,  y_5=5.8426...\, ,
 \ssk\\
  y_6=6.5459...\, ,
y_7=7.2021...\, , \,\, y_8=7.8207...\, ,\,\,
   y_9=8.4083...\, ,
\,\, y_{10}=8.9697...\, ,
 \ssk\\
 y_{11}=9.5086...\, , \,\,
 y_{12}=10.0279...\,, \,\, y_{13}=10.5299...\, , \,\,
  y_{14}=10.0164...\, , \,\,...\, .
   \end{matrix}
  $$


\subsection{Branching at $n=+\iy$}
 \label{S4.4}

Introducing the small parameter in \ef{nInf} for $n \gg 1$,
 $$
 \tex{
 \e=1- \frac n{n+1} \to 0^+ \asA n \to + \iy,
 }
 $$
 and performing a standard (formal, as usual, at least in the
 pointwise sense) linearization yield the following problem:
  \be
  \label{x1}
   \tex{
   |\tilde Y|^{-\frac n{n+1}} \equiv |\tilde Y|^{\e-1} = \frac
   {1+\e \ln |\tilde Y|+O(\e^2)} {|\tilde Y|}
   \LongA
   \BB_\iy
   (\tilde Y)=-{\rm sign}\, \tilde Y\, y \, \ln|\tilde Y|+O(\e),
   }
   \ee
   where $\BB_\iy$ is the unperturbed operator in \ef{ninftyeqn}.
Next, studying the branching at $\e=0$ from the first nonlinear
eigenfunction $\tilde Y_0$, we obtain the corresponding linear
problem:
 \be
 \label{x2}
  \tex{
  \tilde Y= \tilde Y_0 + \e \phi + O(\e^2)
  \LongA \BB_\iy'(\tilde Y_0) \phi= h \equiv -({\rm sign}\, \tilde Y_0)y \,
  \ln|\tilde Y_0|,
  }
  \ee
  where the symmetric linearized operator $\BB_\iy'(\tilde Y_0)$
is given by
 \be
 \label{x3}
  \tex{
  \BB_\iy'(\tilde Y_0)=D^2_y + ({\rm sign}\, \tilde Y_0)' y \,I,
  \quad
({\rm sign}\, \tilde Y_0)'= \d(y-y_0)+ 2 \sum\limits_{(i \ge
1)}(-1)^i \d(y-y_i).
  }
   \ee
Here $\{y_i\}$ are zeros of $\tilde Y_0(y)$ as explained in
Theorem \ref{Th.2}.

At this stage, one then needs proper spectral theory for a
self-adjoint extension of the operator $\BB_\iy'(\tilde Y_0)$.
Then the branching condition reads as the orthogonality
 \be
 \label{x4}
 h \,\bot\, {\rm ker}\, \BB_\iy'(\tilde Y_0).
  \ee
  However, developing such a proper spectral theory faces some
  hard algebraic difficulties. Indeed, looking for eigenfunctions
  $\var_k$,
   \be
   \label{x5}
    \tex{
   \var_k''- \l_k \var_k=C_k={\rm const.} \equiv
   -\big(y_0 \var_k(y_0)+ 2 \sum (-1)^i y_i \var_k(y_i) \big),
    \quad \var_k(y_0)=0,
    }
    \ee
 for $\l_k<0$, we obtain the solution and an algebraic equation for such
 eigenvalues:
  \be
  \label{x6}
   \tex{
  \var_k(y)= \sin(\sqrt {|\l_k|}(y-y_0))
  \LongA
   \fbox{$
  \l_k: \quad
  \sum\limits_{(i \ge 1)}(-1)^i y_i \sin(\sqrt
  {|\l_k|}(y_i-y_0))=0.
  $}
 }
 \ee
 According to \ef{x6}, we require that each eigenfunction
 $\var_k(y)$ should be purely oscillatory as $y \to +\iy$ about
 zero, meaning a ``minimal" (``zero-average") behaviour or zero value in a weak
 sense. For $\l_0=0$, we  take $\var_0(y)=y-y_0$ (again up to a
 normalization multiplier), leading to another algebraic problem,
 which is expected to be non-proper since the eigenfunctions is
 not oscillatory and  not satisfying the ``zero condition"
 at infinity.

Proving that the algebraic equation in \ef{x6} has a discrete
family of solutions $\{\l_k < 0\}$ is a difficult open problem.
Nevertheless, at least, this analysis shows a principal
possibility to detect branching of nonlinear eigenfunctions at
$n=+\iy$, where a simpler algebraic treatment is available, and
more practical asymptotic and other properties of the nonlinear
eigenfunctions than in Theorem \ref{Th.1}, for any $n>0$, are known.

\subsection{The higher-order case $k \ge 2$ for $l=0$}

For the general case of the higher-order ODEs \ef{massconsY},
again  for $l=0$, the equations as $n\to\infty$ are much the same.
Here the scaling constant $C(n)$, for $Y=C(n)\tilde{Y}$ is given
by
\begin{equation*}
C(n) = \big[(2k+1)+n\big]^{-\frac{n+1}{n}}.
\end{equation*}
After scaling, this yields the ODE
\begin{equation}
 \label{MM1}
  \BB_\iy^{(2k)}(\tilde Y) \equiv
(-1)^{k+1}D^{2k}_y \tilde{Y}  +{\rm sign}\,{\tilde{Y}} \, y=0,
\end{equation}
which deserves further study by deriving the corresponding
algebraic structures.
 Note that,
close to the interface at $y=y_0^+$, for $k \ge 2$ (the case $k=1$
is not oscillatory, as we have seen in Section \ref{S3}),
 a stabilization to a {\em periodic oscillatory component} $\phi(s)$, with the typical structure
 near the interface at $y=y_0^+$,
  \be
  \label{SS111}
   \tex{
   \tilde Y(y)=(y-y_0)^{2k} (\phi(s)+o(1)) \whereA s=-\ln(y-y_0) \to - \iy,
   }
   \ee
   of solutions of \ef{MM1} is most plausible. See
   \cite{Gl4, Gl6}, where such  oscillatory sign changing
   solutions such as \ef{SS111}, with a periodic component $\phi(s)$ satisfying a nonlinear
   higher-order ODE, have been found for the fourth- and sixth-order thin
   film equations.

 However, explicitly solving the ODE \ef{MM1}  for any $k \ge 2$
  in the positivity and negativity domains leads
to much more complicated algebraic systems, which do not admit
such a clear explicit resolving as for $k=1$. Moreover, even
shooting numerical methods lead to difficult and unclear
results, so we do not present such an analysis here.

\section{The second eigenfunction $\tilde Y_1(y)$ for $n=+\infty$:
algebraic approach}
 \label{S41}

Using the same method, the equations governing the behaviour as
$n\to\infty$, relating to the second  nonlinear eigenfunction
 $\tilde Y_1(y)$,
 can
easily be found. Namely and analogously, for $l=1$ in the ODE
\eqref{1stmomintcons}, with $k=1$, the limit equation at
$n=+\infty$ is given by
\begin{equation}
 \label{Y11}
\mbox{$\tilde{Y}^{\prime\prime}y - \tilde{Y}^\prime +
{\rm sign} \, \tilde Y \, y^2 = 0$}.
\end{equation}
This ODE \ef{Y11} can be written in a singular Sturm--Liouville
form
 \be
 \label{Y11N}
  \tex{
  \big( \frac {\tilde Y'}y \big)'=- {\rm sign} \, \tilde Y,
  }
  \ee
 where the weight $\rho(y)= \frac 1y \not \in L^1_{\rm loc}$.
 Therefore, $y=0$ is a singular inner point, where an extra
 condition must be posed. This is done by checking the functional
 weighted $L^2$-space corresponding to the operator in \ef{Y11N} and its available asymptotics
  as $y \to 0$:
  \be
  \label{SL1NN}
   \tex{
   \int_0 \rho(y) (\tilde Y'(y))^2\, {\mathrm d}y < \infty
   \LongA \tilde Y'(0)=0.
   }
   \ee
   Next,
   integrating \ef{Y11N} twice  yields, in the positivity and
  negativity domains, the following:
   \be
   \label{Y12}
   \tex{
 \tilde Y_\pm (y)= \mp \, \frac {y^3}3 + a y^2+b, \quad a,b, \in \re.}
  \ee
  Note that, unlike \ef{casesninfty} for $l=0$, here the linear
  term $\sim y$ is absent in the cubic polynomial.
  The analysis of the algebraic matching system corresponding to
  \ef{Y12} is similar in many places but a couple of ones, which
  we will concentrate upon now.

  \begin{theorem}
 \label{Th.l1}
 The problem $\ef{Y11N}$, $\ef{tr00}$ admits a unique nontrivial solution
 $\tilde Y_1(y)$ with transversal zeros at $\{y_i\}_{i \ge 1}$,
 where
  \be
  \label{i1}
   \tex{
  y_0=-1, \quad y_1= \frac 12, \quad y_2= \frac{5+3 \sqrt
  5}4=2.927050...\, ,
   }
   \ee
   and triples of zeros $\{y_{i-1},y_i,y_{i+1}\}_{i\ge 2}$
 satisfy a cubic homogeneous algebraic equation,
 \be
 \label{deg3}
 y_{i-1}y_{i+1}^2+ y_i y_{i+1}^2+ y_{i-1}^2 y_i
 +y_{i-1}^2y_{i+1}=y_{i-1}y_i y_{i-1}+ y_i^2 y_{i+1}+y_{i-1}
 y_i^2+ y_i^3
  \forA i \ge 2.
  \ee
 \end{theorem}

\noi{\bf Remark.} \ef{deg3} reduces to a 3D linear discrete
equation for
 $
 \a_{i,j,k}=y_iy_jy_k,
 $
 but, similar to \ef{tr15}, this does not essentially help to get any finite explicit
 expressions for zeros $\{y_i\}$.

 \ssk

 \noi{\em Proof.} The first step is elementary: on $(-1,y_1)$, in
 view of \ef{Y12} (cf. also \ef{SL1NN}),
  \be
  \label{p1}
  \tex{
  \tilde Y_+(y)=(y+1)^2(y-y_1)\big(- \frac 13\big) \whereA y_1=
  \frac 12.
  }
  \ee
  Next, on $(y_1,y_2)$ by \ef{Y12},
   \be
   \label{p2}
    \tex{
    \tilde Y_-(y)=(y-y_1)(y-y_2)\big(\frac 13\, y+c_1\big) \whereA
    c_1= \frac 13 \,\frac{y_1y_2}{y_1+y_2},
    }
    \ee
    so that the matching \ef{B2} at $y=y_1$ yields
     \be
     \label{p3}
      \tex{
     (y_1+1)^2\big(- \frac 13\big)=(y_1-y_2)\big(\frac 13\,
     y_1+c_1\big)
     \LongA 4 y_2^2-10 y_2-5=0,
 }
  \ee
  whence the desired value of the root $y_2$ in \ef{i1}.

Similarly, in the general case, on $(y_{i-1},y_i)$,
 \be
 \label{p4}
  \tex{
 \tilde Y_-(y)=(y-y_{i-1})(y-y_i)\big(\frac 13\, y+c_{i-1}\big)
 \whereA
  c_{i-1}= \frac 13 \, \frac{y_{i-1}y_i}{y_{i-1}+y_i},
  }
  \ee
  and
on $(y_{i},y_{i+1})$,
 \be
 \label{p5}
  \tex{
 \tilde Y_+(y)=(y-y_{i})(y-y_{i+1})\big(-\frac 13\, y+c_{i}\big)
 \whereA
  c_{i}=- \frac 13 \, \frac{y_{i}y_{i+1}}{y_{i}+y_{i+1}}.
  }
  \ee
 Then the standard matching via \ef{B2} of such $\tilde Y_\pm(y)$ yields
 the homogeneous cubic algebraic equation \ef{deg3} for the zero triple
 $\{y_{i-1},y_i,y_{i+1}\}$.
  $\qed$

 \ssk

 Figure \ref{FF3} shows the general structure of the second
eigenfunction $\tilde Y_1(y)$ on the large interval
$[y_0=-1,200]$. The envelope of the decaying oscillations are
governed by the algebraic curve
 \be
 \label{L0N}
  \tex{
  L_1(y) \approx \pm 0.89 \, y^{\frac 13} \asA y \to + \iy,
   }
   \ee
 which can be associated with the algebraic manipulations
 involved.


 \begin{figure}
\centering
\includegraphics[scale=0.85]{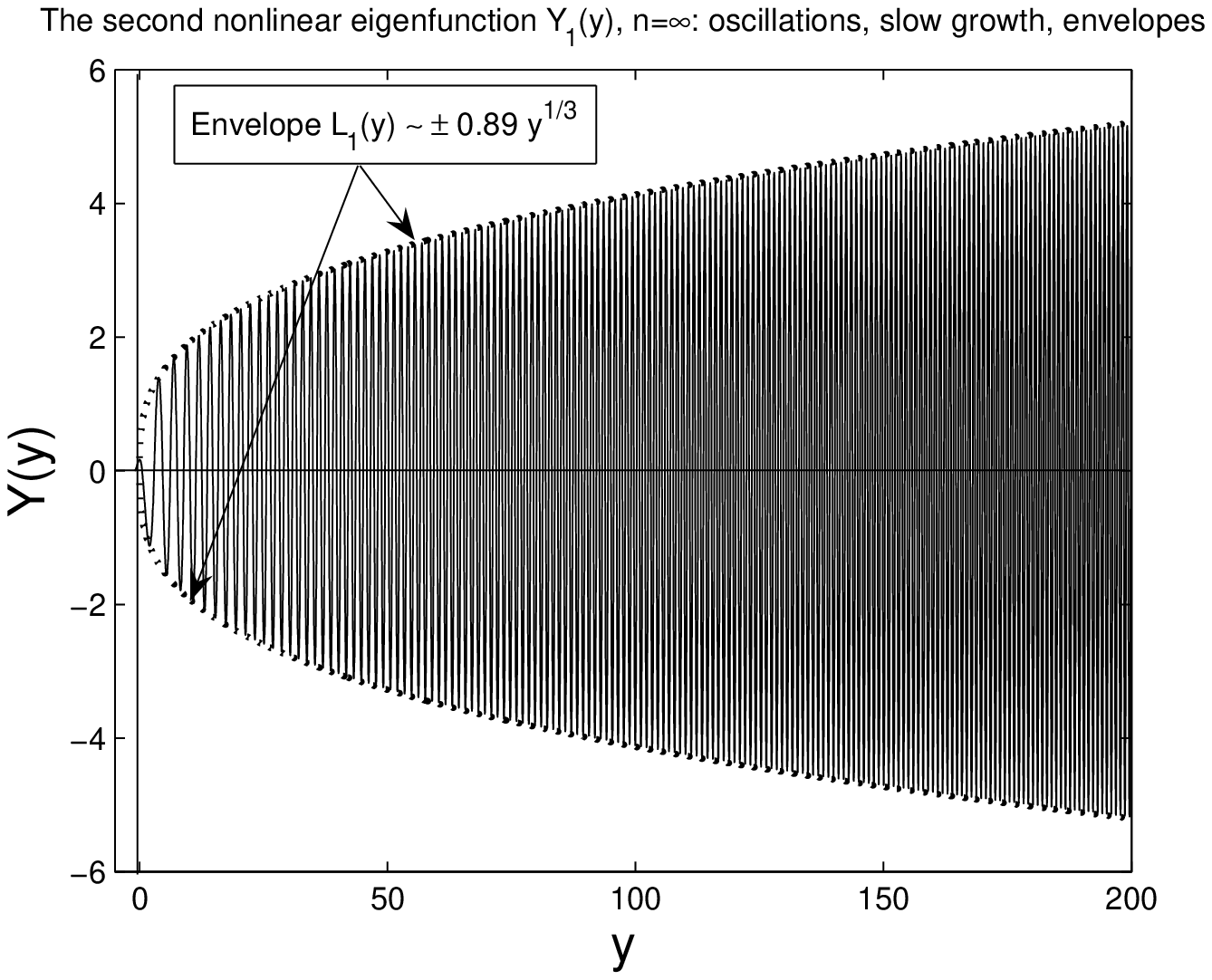}  
\vskip -.3cm
  \caption{The second nonlinear eigenfunction
 $\tilde Y_1(y)$ as a solution of the problem \ef{Y11N}, \ef{ninftyshtY};
 $k=1$, $l=1$.}
 \label{FF3}
\end{figure}

The next, Figure \ref{FF4}, shows first zeros of $\tilde Y_1(y)$ on
the interval $[-1,5]$.


 \begin{figure}
\centering
\includegraphics[scale=0.85]{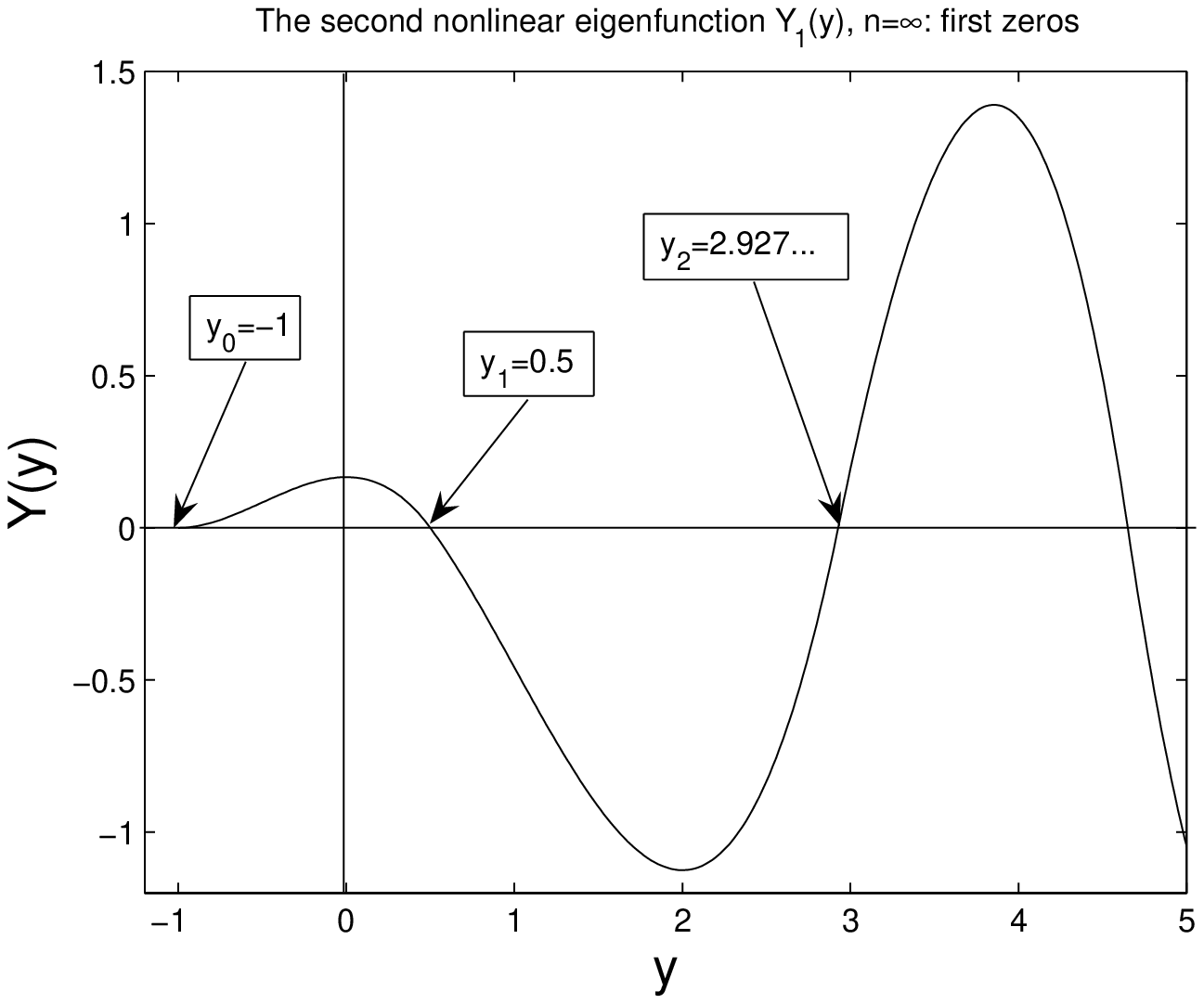}  
\vskip -.3cm
  \caption{The second nonlinear eigenfunction
 $\tilde Y_1(y)$ as a solution of the problem \ef{Y11N}, \ef{ninftyshtY};
 $k=1$, $l=1$.}
 \label{FF4}
\end{figure}

In general, for all values of $k \ge 1$, we have the limit
equation at $n=+\iy$ for $l=1$ in the form:
\begin{equation*}
\mbox{$(-1)^{k+1}D_y^{2k}\tilde{Y}y + (-1)^kD_y^{2k-1}\tilde{Y} +
 {\rm sign} \, \tilde Y
 = 0$}.
\end{equation*}
This can be integrated once, however, one cannot expect any
reasonably easy algebraic manipulations leading to some explicit
representation of asymptotic properties of $\tilde Y_1(y)$ for any
$k \ge 2$.

\section{The third eigenfunction $\tilde Y_2(y)$ for $n=+\infty$:
algebraic approach}
 \label{S42}

   Similarly, for the third eigenfunction, where $l=2$,
we have, in the lower-order case $k=1$, the limit ODE
\begin{equation}
 \label{Y22}
\mbox{$\tilde{Y}^{\prime\prime}y^2 - 2\tilde{Y}^\prime y +
2\tilde{Y} + {\rm sign} \, \tilde Y \, y^3 = 0$}.
\end{equation}
The corresponding Sturm--Liouville form
 \be
 \label{SL1}
 \tex{
 \big( \frac{\tilde Y'}{y^2}\big)'=- \frac 2{y^4}\, \tilde Y-  \frac 1y \,{\rm
 sign} \, \tilde Y,
 }
 \ee
 reveals even more singularity at $y=0$ as it used to be in
 \ef{Y11N} for $\tilde Y_1(y)$. The necessary extra condition at
 $y=0$ is presented below at \ef{p2N1}.

Writing \ef{Y22} down for $\tilde Y_{\pm}(y)$ in the form
 \be
 \label{Y221}
  \tex{
  Y''=\mp \, y+ 2\,\big(\frac {\tilde Y}y\big)'
  }
  \ee
  and integrating twice yields
   \be
   \label{Y222}
   \tex{
   Y_\pm(y)=\mp \, \frac {y^3}2+ay^2+by, \quad a,b \in \re.
   }
   \ee
Surprisingly, the matching at zeros $y=y_i$ by using the cubic
polynomials \ef{Y222}, having $y=0$ as a fixed zero always leads to
a simpler mathematics.

  \begin{theorem}
 \label{Th.3}
 The problem $\ef{Y22}$, $\ef{tr00}$ admits a unique nontrivial solution
 $\tilde Y_1(y)$ with transversal zeros at $\{y_i\}_{i \ge 1}$,
 where
  \be
  \label{i1N}
   \tex{
  y_0=-1, \quad y_1= 0, \quad y_2= 1+ \sqrt 2=2.4142...\,, \quad
  y_3=1+ 3 \sqrt 2=5.2426...\,,
   }
   \ee
   and further zeros are given by the second-order linear
discrete equation: for any $i \ge 3$,
 \be
 \label{i1N1}
  \tex{
  y_{i+1}-2 y_i+y_{i-1}=0 \LongA
  y_i=C_1+C_2 i \whereA C_1=1-3 \sqrt 2, \,\, C_2=2 \sqrt 2.
  }
  \ee
  In particular, the distribution of zeros is uniform:
   \be
   \label{i1NN}
   y_{i+1}-y_i=2 \sqrt 2=2.8284... \quad \mbox{for all} \,\,\, i
   \ge 2.
   \ee
 \end{theorem}

 \noi{\em Proof.} On $(-1,y_1)$, in
 view of \ef{Y222},
  \be
  \label{p1N}
  \tex{
  \tilde Y_+(y)=y(y+1)^2\big(- \frac 12\big) \LongA y_1=0.
  }
  \ee
  Next, on $(y_1=0,y_2)$ by \ef{Y222},
   \be
   \label{p2N}
    \tex{
    \tilde Y_-(y)=y(y-y_2)\big(\frac 12 \, y+c_1\big)
     \whereA
    c_1= \frac 1{2 y_2}.
    }
    \ee
 Unlike the previous cases, it is not possible to find $y_2$ and
 $c_1$ by using the standard matching conditions \ef{B2}. The
 point is that the differential  operator in \ef{SL1} is strongly singular
 at $y=0$,
 where the weight $\rho(y)= \frac 1{y^2} \not \in L^p(-1,1)$ for
 any $p \ge 1$.

 It then follows from \ef{Y222} due to the singular setting \ef{SL1} that
 the  two usual conditions at $y=0$
 such as the values of $\tilde Y(0)=0$ and of a given ``flux" $\tilde
 Y'(0)$ are not sufficient to determine a unique local solution for
 $y>0$ and $y<0$. To get a unique solution, the value of $\tilde
 Y''(0)$ should be prescribed.
 A further analysis shows that  a proper stronger continuity condition
 of matching at $y=0$ is necessary and this includes  the equality of the
 second-order derivatives:
  \be
  \label{p2N1}
   Y_+''(0^-)=Y_-''(0^+).
    \ee
 Overall, this yields the following quadratic equation for $y_2$:
  \be
  \label{t1}
   \tex{
   Y_+''(0)=-2=Y_-''(0)= 2 c_1-y_2, \,\,\, c_1= \frac 1{2 y_2}
    \LongA y_2^2-2 y_2-1=0,
   }
   \ee
 and this uniquely defines $y_2$ shown in \ef{i1N}.

 Next, we use \ef{p2N} with the obtained values of $y_2$ and $c_1$
 to match (now, in a standard way) with the solution representation on $(y_2,y_3)$,
  \be
  \label{yy11}
   \tex{
   \tilde Y_+(y)=y(y-y_2)(y-y_3)\big(- \frac 12\big),
   }
   \ee
   to get at $y=y_2$
    \be
    \label{yy12}
     \tex{
     y_2\big( \frac 12\, y_2+c_1\big)= y_2(y_2-y_3)\big(-\frac
     12\big) \LongA y_3=2y_2+2 c_1 =1+ 3 \sqrt 2.
     }
     \ee

 Finally, in the general case, on $(y_{i-1},y_i)$ for $i \ge 3$,
 \be
 \label{p4N}
  \tex{
 \tilde Y_+(y)=y(y-y_{i-1})(y-y_i)\big(-\frac 12\,\big),
  }
  \ee
  and
on $(y_{i1},y_{i+1})$,
 \be
 \label{p5N}
  \tex{
 \tilde Y_-(y)=y(y-y_{i})(y-y_{i+1})\frac 12.
  }
  \ee
 By
  the standard matching at $y=y_i$ via \ef{B2} of such $Y_\pm(y)$ yields
 \be
 \label{p6N}
  \tex{
  y_i(y_i-y_{i-1})\big(- \frac 12 \big)= y_i(y_i-y_{i+1})\frac 12,
  }
  \ee
  whence the linear difference relation \ef{i1N1} with some
  constants $C_1$ and $C_2$. These are uniquely obtained from
  the linear algebraic system,
   \be
   \label{ss1NN}
   \left\{
   \begin{matrix}
   y_2=1+\sqrt 2=C_1+2C_2,\,\,
     \ssk\\
   y_3=1+ 3 \sqrt 2= C_1+3 C_2
   \end{matrix}
   \right.
 \LongA
 \left\{
 \begin{matrix}
 C_1=1-3 \sqrt 2,\ssk\\
 C_2=2\sqrt 2, \quad\,\,\,
 \end{matrix}
   \right.
   \ee
  completing the proof.  $\qed$

 \ssk

 Figure \ref{FF5} shows the general structure of the third
eigenfunction $\tilde Y_2(y)$ on the large interval
$[y_0=-1,200]$. Directly connected with \ef{i1N1}, the envelope of
the decaying oscillations is linear
 \be
 \label{L0NN}
  \tex{
  L_2(y) \approx \pm  y \asA y \to + \iy.
   }
   \ee
First zeros  of $\tilde Y_2(y)$ on the interval $[-1,7]$ are shown
in
  Figure \ref{FF6}.


 \begin{figure}
\centering
\includegraphics[scale=0.85]{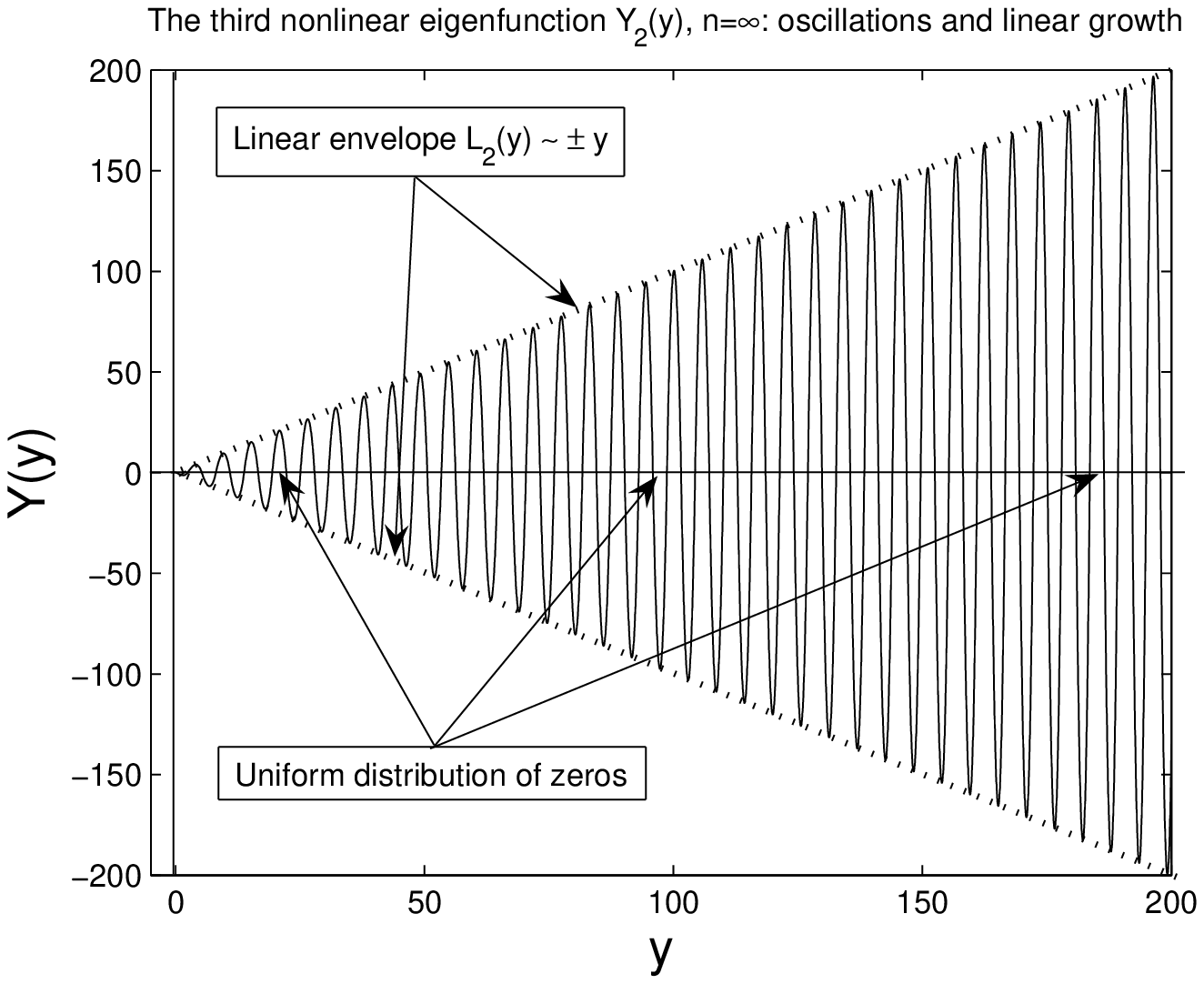}  
\vskip -.3cm
  \caption{The third nonlinear eigenfunction
 $\tilde Y_2(y)$ as a solution of the problem \ef{Y221}, \ef{ninftyshtY};
 $k=1$, $l=2$.}
 \label{FF5}
\end{figure}



 \begin{figure}
\centering
\includegraphics[scale=0.85]{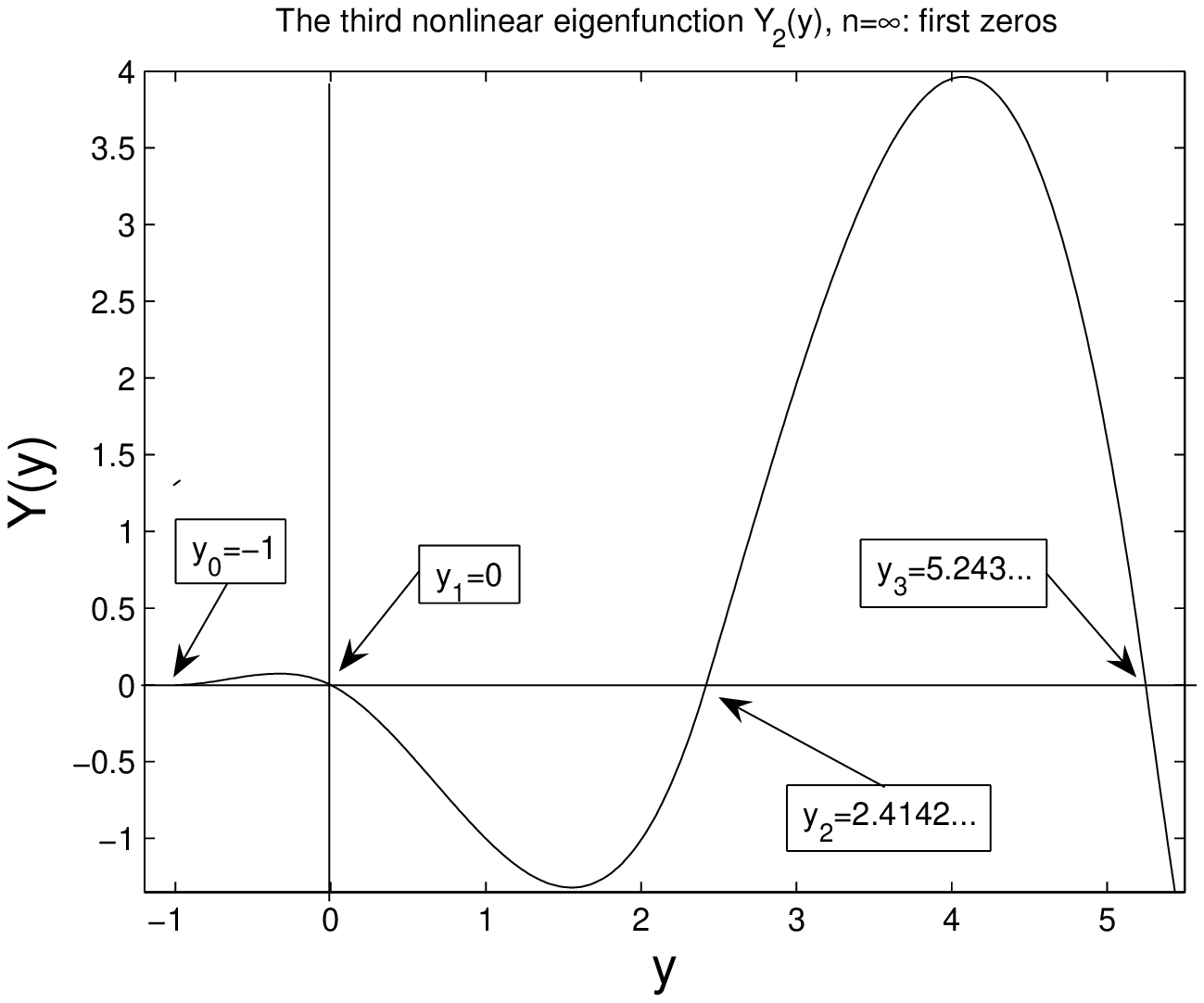}  
\vskip -.3cm
  \caption{The third nonlinear eigenfunction
 $\tilde Y_2(y)$ as a solution of the problem \ef{Y221}, \ef{ninftyshtY};
 $k=1$, $l=2$.}
 \label{FF6}
\end{figure}

\ssk

 For arbitrary  $k \ge 2$, the limit ODE for $\tilde Y_2(y)$ is
\begin{equation*}
\mbox{$(-1)^{k+1}D_y^{2k}\tilde{Y}y^2 +
2(-1)^kD_y^{2k-1}\tilde{Y}y + 2(-1)^{k+1}D_y^{2k-2}\tilde{Y} +
\frac{\tilde{Y}}{|\tilde{Y}|} \, y^3 = 0$},
\end{equation*}
which does not admit a clear geometric-algebraic method of
solution.


\ssk

Finally, we recall that exotic profiles such as $\tilde
Y_0(y)$ in Figure \ref{FF1}, $\tilde Y_1(y)$ in Figure \ref{FF3}
and $\tilde Y_2(y)$ in Figure \ref{FF5} are not just some
functions defining some self-similar solutions in the nonlinear
limit ($n= +\iy$) equation \ef{v12k1} for $k=1$,
 but these are the
most stable asymptotic pattern of such a  nonlinear PDE. Indeed,
$\tilde Y_0(y)$ is the most stable and is expected to attract,
as $t \to +\iy$, a.a. solutions, excluding only those with data
satisfying some extra orthogonality conditions. For instance,
those having zero mass (or also the zero moment for $Y_2(y)$ to
play a role), so that an extra time-scaling is necessary to get
convergence to the second nonlinear eigenfunction $\tilde Y_1(y)$,
etc. Of course, these stability questions are far beyond the scope
of this paper, are very difficult, and remain open for most
higher-order NDEs.

\section{Nonlinear dispersion equation with  absorption}
 \label{S5}

 \subsection{A full quasilinear NDE: homotopy deformation to linear PDEs}

The final natural progression, from the nonlinear model
\eqref{linnonl}, is to go to the odd-order NDE with absorption
\ef{nde33}.
This links up both the nonlinear model and the semilinear one
 \cite{RayGI},
\begin{equation}
\label{semilinear pde}
 u_t = (-1)^{k+1}D _x^{2k+1}u - |u|^{p-1}u
\quad \text{in} \quad \mathbb{R}\times\mathbb{R}_+ \quad (n=0).
\end{equation}
However, the NDE \ef{nde33} is indeed a more difficult quasilinear
equation than \ef{semilinear pde}, with not that well understood
phenomena of ``nonlinear" bifurcation and branching. Therefore, we
should
 pay here less effort towards establishing rigorously
 some of the key analytical aspects concerning  basic
similarity solutions of \ef{nde33}.

The NDE \eqref{nde33}, for any $n>0$, $p>n+1$, after similarity
scaling,
 \be
 \label{P11}
  \tex{
  u_{\rm gl}(x,t)=t^{-\frac 1{p-1}}f(y), \quad y =\frac x{t^\b} \whereA \b=
  \frac{p-(n+1)}{(p-1)(2k+1)},
   }
   \ee
 reduces to the
 ODE
\begin{equation}
\label{NDEabsODE}
    \mbox{$(-1)^{k+1}D^{2k+1}_y(|f|^nf) +
\frac{p-(n+1)}{(p-1)(2k+1)}\, f^\prime y + \frac{1}{p-1}\,f  -
|f|^{p-1}f=0$} \inB \re.
\end{equation}
As usual, a proper setting for \ef{NDEabsODE} assumes a finite
left-hand interface at some $y=y_0<0$ and an admissible
oscillatory behaviour as $y \to +\iy$, which was under scrutiny
above.

Recall that, using the reflections \ef{refl1}, simultaneously, we
construct blow-up solutions
 \be
 \label{P11bl}
  \tex{
  u_{\rm bl}(x,t)=(T-t)^{-\frac 1{p-1}}f(y), \quad y = - \frac x{(T-t)^\b} \whereA \b=
  \frac{p-(n+1)}{(p-1)(2k+1)}
   }
   \ee
   (since $\b>0$ for $p>n+1$, formally, this is a {\em single point
   blow-up})
of the {\em NDE with source},
 \be
 \label{Sbl}
 u_t=(-1)^{k+1}D_x^{2k+1}(|u|^n u)+|u|^{p-1} u,
  \ee
so that proper profiles $f(y)$ describe both global and blow-up
asymptotics of such NDEs.

\ssk

In the case $n=0$ in \ef{NDEabsODE}, we have a simpler semilinear
ODE corresponding to the  model \eqref{semilinear pde}. It is
important also to note that, in addition to the homotopy path $n
\to 0^+$ in \ef{nde33}, it is useful  to apply an extra limit $p
\to 1^+$. For the case $n=0$ and $p=1$ in \eqref{nde33}, we reduce
to the linear equation
\begin{equation*}
\mbox{$ u_t = (-1)^{k+1}D^{2k+1}_xu - u$}.
\end{equation*}
 Using the substitution
\begin{equation*}
u(x,t) = \mathrm{e}^{-t}v(x,t),
\end{equation*}
reduces it  to the standard linear dispersion  equation
\begin{equation}
 \label{v23}
v_t = (-1)^{k+1}D^{2k+1}_xv.
\end{equation}
 The extra scaling yields a semigroup (a group) with the
 infinitesimal generator $\BB$ in \ef{HermB}:
  \be
  \label{v30}
   \tex{
   v(x,t)=t^{-\frac 1{2k+1}}w(y,t), \quad y=x/t^{\frac 1{2k+1}},
   \quad \t= \ln t \LongA w_\t =\BB w.
    }
    \ee
 Therefore, linear Hermitian theory developed in \cite[\S~4]{RayGI} (see also
\cite[\S9]{2mSturm}) can be used, and this leads to an efficient
 approach based on comparison of linear eigenfunction structures
 of (\ref{v23}) and the nonlinear ones of (\ref{nde33}).
 Actually, this implies existence of a certain ``homotopy" of
 those PDEs as $n \to 0^+$ (already noted in Section \ref{S3.5}) and
 as $p \to 1^+$. Overall, this establishes a countable nature of
 the
 nonlinear eigenfunction family for (\ref{nde33}), which is
 difficult to prove rigorously.
 We refer to \cite{VSSSParaPDEs}, where such
 a homotopy is discussed for a related fourth-order PDE.

Thus, the continuous double limit, $n \to 0^+$ and $p \to 1^+$,
 by reducing to the LDE \ef{v23} with a {\em countable} set of
 rescaled eigenfunctions \ef{N4},
 in
general, justifies that the ODE \ef{NDEabsODE} admits a countable
set of so-called $p$-bifurcation branches of solutions, which
blow-up as $p \to (n+1)^+$. Further study of such nonlinear
phenomena is of importance. For $n=0$, this branching phenomenon
has been studied in \cite{RayGI}.

It is worth mentioning that, in general, such a {\em full}
branching approach assumes existence of a
 countable number of 2D bifurcation surfaces defined in the 2D parameter space of
 $(n,p)$, and the critical point $(0,1)$ is expected to be a complicated singular
 ``cusp" induced by such a collection of bifurcation surfaces.
 This is a very difficult bifurcation-branching phenomenon, which
 in 2D was not fully understood and remains an open problem.
 Therefore, below, we study another formal, but 1D,
``nonlinear bifurcation" phenomenon in order to explain existence
of a countable number of $p$-bifurcation branches in the problem
\ef{NDEabsODE}.

\subsection{Local ``nonlinear" bifurcations at critical exponents $p_l(n)$}
 \label{S5.2}

 We consider the VSS solutions \ef{P11} of \ef{nde33},
 and develop a formal {\em nonlinear} version of such a 1D
$p$-bifurcation (branching) analysis for any fixed $n>0$. As
usual, according to classic branching theory \cite{GeoMethNonAn,
VainbergTr}, a justification (if any) is performed for the
equivalent semilinear integral equation with compact operators in
suitable metrics. For simplicity, we present computations in  the
differential setting, which does not change anything essentially.
Note that, for nonlinear odd-order operators, some issues of
compactness can be rather tricky.

Let us compare the time-factor structure of the VSS \ef{P11} and
that for the pure NDE in \ef{linnonlss}. It follows that
 the critical bifurcation
exponents $\{p_l=p_l(n)\}$ are then determined from the equality
of the exponents:
 \be
 \label{pp1}
  \mbox{$
   -\frac 1{p_l-1}=- \a_l
 \quad \Longrightarrow \quad p_l(n)= 1 + \frac 1{
 \a_l(n)} \quad l=0,1,2, \, ... \, ,
  $}
  \ee
 where $\a_l(n)$ are the critical exponents as in \ef{alk}
 obtained explicitly and further fully nonlinear ones that cannot be
 determined dimensionally (via conservation laws).

In particular, for the semilinear case $n=0$, the eigenvalues
$\a_l(0)$ are given by \ef{N3}, and this leads to the critical
bifurcation exponents
 \be
 \label{KK1}
  \tex{
   p_l(0)= 1+ \frac{2k+1}{l+1}, \quad l=0,1,2,...\,,
   }
   \ee
 for the semilinear equation \ef{semilinear pde} studied in
 \cite[\S~7.4]{RayGI}.
 In the present nonlinear case with a fixed $n>0$, such a standard
 linearised approach is not suitable.

However,  first steps of this ``nonlinear" bifurcation theory
 are straightforward.
 We
next use an expansion
 relative to the small parameter $\e= p_l-p$, so that, as $\e \to 0$,
 $$
 \mbox{$
 \begin{matrix}
 \a= \frac 1{p-1}= \frac 1{p_l-1- \e}= \a_l +  \a_l^2 \e +O(\e^2) \, ,
 \ssk\ssk\\
 \b= \frac{p-(n+1)}{(p-1)(2k+1)} = \frac {1-\a_l n}{2k+1}- \frac {n
 \a_l^2}{2k+1}\, \e + O(\e^2), \ssk\ssk\ssk\\
 |f|^{p-1}f=|f|^{p_l-1-\e}f= |f|^{p_l-1}f (1-\e \ln|f| + O(\e^2)).
  \end{matrix}
   $}
    $$
    Note that, unlike
  the case \ef{fn0}, the last expansion has a clearer
 functional validity, since at $f=0$, there occurs standard
 issues of convergence, which
makes  sense suitable for
 passing to the limit in the integral operators.

Substituting these expansions
into (\ref{NDEabsODE}) and collecting $O(1)$ and $O(\e)$-terms
yield
  \be
 \label{eq4}
  \mbox{$
 {\bf A}(f,\a_l)- |f|^{p_l-1}f + \e {\mathcal L}f + \e |f|^{p_l-1}f \ln |f|+
 O(\e^2)=0, \,\,
  {\mathcal L}= - \frac {n \a_l^2}{2k+1}\,y D_y+ \a_l^2 I,
   $}
   \ee
  where $\AAA(f,\a_l)$ is the nonlinear operator in \ef{rednonlOde}.
 Recall that, at each nonlinear eigenvalue $\a=\a_l$, there exists the corresponding
 nonlinear eigenfunction $f_l$ such that \ef{BC2} holds.
 At least, we are going to use this conclusion, which was not completely
 proved.
 The fact is that the operator
 ${\bf A}(f,\a)$, with $\a=\a_l$ in (\ref{eq4}) of the
rescaled
 pure NDE, correctly describes the essence of a
 {\em ``nonlinear bifurcation phenomenon"} to be revealed.

To this end, we use the additional invariant scaling of the
operator ${\bf A}(f,\a)$ by introducing the new unknown function
$F(\cdot)$ as follows:
  \be
 \label{ep1}
 f(y)=b  F(y/b^{\frac n{2k+1}}),
  \ee
 where $b=b(\e)>0$ is also an unknown parameter  to be determined
 from a scalar branching equation
   and satisfying
  \be
  \label{bb1}
   b(\e) \to 0 \quad \mbox{as} \quad  \e \to 0.
   \ee
 Substituting (\ref{ep1}) into (\ref{eq4}) and, under natural (but
 not proved) regularity assumptions on such expansions,
 omitting
 all higher-order terms (including the one with the logarithmic multiplier $\ln
|b(\e)|$)
  yield
  \be
 \label{ep2}
 {\bf A}(F,\a_l) - b^{p_l-1} |F|^{p_l-1} F+
  \e {\mathcal L}
    F =0.
  \ee

 Finally, we perform linearization about the nonlinear
 eigenfunction $f_l(y)$ by setting
  $$
  F=f_l+Y.
   $$
    This yields the following linear
 non-homogeneous problem:
   \be
  \label{ep3}
  {\bf A}'(f_l,\a_l)Y =  b^{p_l- 1} |f_l|^{p_l-1} f_l -
   \e {\mathcal L} f_l.
  \ee
 Here, the derivative is given by
  \be
  \label{AA1}
   \mbox{$
  {\bf A}'(f_l,\a_l)Y= (n+1)(-1)^{k+1} D_y^{2k+1}(|f_l|^n Y) +
  \frac{1-\a_ln}{2k+1}\, Y' y + \a_l Y.
   $}
   \ee

   As usual in bifurcation theory,
 the rest of the analysis crucially  depends on assumed good  spectral properties of
 the linearised operator ${\bf A}'(f_l,\a_l)$, which are not easy at all, and,
 in fact, are much more
 complicated than that for the pair $\{\BB,\BB^*\}$ in \cite{RayGI}.
 Many aspects of such a theory remain quite obscure. However, we
 proceed to explain the key final results on possible bifurcations,
  and follow the same lines.
  A proper functional setting of this
 operator is more understandable in the present 1D case,
  where, using the behaviour of
 $f_l(y) \to 0$ as $ y \to y_0^+$ and $y \to +\iy$, it is, at least formally,
  possible to check whether the
 resolvent is compact in a suitable weighted $L^2$ space. In
 general, this is a difficult open problem, especially since we are not aware of precise asymptotic properties of all the
 eigenfunctions $f_l$.

We assume that such a proper functional setting is available
 for ${\bf A}'(f_l,\a_l)$. Therefore,
we deal with operators having solutions with ``minimal"
singularities at the boundary point of the support at $y=y_0<0$,
where the operator is degenerate and singular. The same is assumed
 at $y=+\iy$, where the necessary admitted bundle of solutions
should be identified to pose singular boundary conditions; see
Naimark's monographs on ordinary differential operators
\cite{NaimarkI, NaimarkII} as a guide.

  Namely (cf. \cite[\S~4]{RayGI}), we assume that the linear odd-order operator
   ${\bf A}'(f_l,\a_l)$ has a discrete spectrum, and a complete and closed set of
eigenfunctions
   denoted again by $\{\psi_\g\}$. We also assume that the kernel is finite-dimensional
   and we are able to
  determine the spectrum, eigenfunctions $\{\psi_\g^*\}$, and  the kernel of the
 ``adjoint" operator $({\bf A}'(f_l,\a_l))^*$ defined in  a natural
  way using the topology of the dual space $L^2$ (or, equivalently and possibly, of a space
   with an indefinite metric)
     and having
  the same point spectrum.
    The latter  is true for compact operators in
 suitable spaces in more standard setting,
  \cite[Ch.~4]{FuncAnaly}.
  We also require that
 the bi-orthonormal eigenfunction
  subset $\{\psi_\g\}$ of the operator ${\bf A}'(f_l,\a_l)$
  is complete and  closed in a suitable weighted $L^2$-space (for $n=0$,
   such results are available \cite{RayGI}). Note that, often,  this ``spectral collection" is
   too exhaustive in nonlinear operator theory; see Deimling
   \cite[p.~412]{NonFuncAnaly} for most general bifurcation results.


  Thus, by a typical Fredholm-like
   alternative, the unique solvability of \ef{ep3} requires the
   orthogonality of the inhomogeneous term therein to the ${\rm ker}\,
   \AAA'(f_l,\a_l)$. For simplicity,  let it be 1D with the eigenfunction $\phi_l$,
   so  the right-hand side in \ef{ep3} satisfies
  \be
  \label{ort1}
   b^{p_l- 1} |f_l|^{p_l-1} f_l -
   \e {\mathcal L} f_l \,\, \bot \,\,{\rm ker}\,
   \AAA'(f_l,\a_l)={\rm Span}\,\{\phi_l\}.
   \ee
 Then,
multiplying (\ref{ort1}) by $\phi_l^*$  in
$L^2$ (or within the equivalent  indefinite metric)
 yields  the
orthogonality condition
(Lyapunov-Schmidt's algebraic  branching equation
\cite[\S~27]{VainbergTr}):
  \be
 \label{ep41}
  \mbox{$
 b^{p_l-1}\langle |f_l|^{p_l-1} f_l, \phi_l^* \rangle =
 \e \langle {\mathcal L} f_l, \phi_l^* \rangle.
   $}
    \ee
    Similar to (\ref{M5}), one needs to check whether
    the constants are non-zero:
 \be
 \label{CC1}
 \mbox{$
\langle |f_l|^{p_l-1} f_l, \phi_l^* \rangle \not =0 \quad
\mbox{and} \quad \langle {\mathcal L} f_l, \phi_l^* \rangle \not =
0, $}
 \ee
which is not an easy problem and can lead to some restrictions
for such a behaviour, though is crucial for any hope to see a
bifurcation point.

Under the conditions (\ref{CC1}), the parameter $b(\e)$  in
(\ref{ep1}) for $p \approx p_l$ is given by
  \be
 \label{ep5}
 \mbox{$
 b(\e) \sim [\g_l(p_l-p)]^{\a_l(n)} \,\, \quad \big(\frac
 1{p_l-1}=\a_l\big), \quad \e=p_l-p, \quad \g_l= \frac{\langle {\mathcal L} f_l, \phi_l^*
 \rangle}
 {\langle |f_l|^{p_l-1} f_l, \phi_l^* \rangle}\, .
 $}
     \ee
The direction of developing in $p$ of each $p_l$-branch and
whether the bifurcation is sub- or supercritical depend on the
sign on the coefficient $\g_l$.
 This can be checked
numerically only, but, in general, we expect that $\g_l>0$, so
that  these nonlinear bifurcations are {\em subcritical} and  the
$p_l$-branches exist for $p < p_l$.

Overall, the above formal analysis detects a number of key
assumptions, which are necessary for such a nonlinear bifurcation
to occur at the critical exponents $p_l$ given by \ef{pp1}. Recall
again that,
    for $n=0$, a more rigorous justification of the corresponding linearised bifurcation
    analysis is done in \cite{RayGI},  where a
    countable number of $p$-branches was shown to be originated at bifurcation
    points (\ref{KK1}).

  \ssk

  Overall, we claim that, for any $n \ge 0$, the ODE \ef{NDEabsODE}, with proper
  setting as $y \to \pm \iy$,
   \be
   \label{Count1}
   \fbox{$
   \mbox{admits not more than {\em countable} number of
   $p$-branches of solutions},
    $}
    \ee
    which are originated at the critical exponents \ef{pp1}
    (no rigorous proof is available still).


\subsection{Numerical experiments for $k=1$}

We briefly attempt to find numerical solutions of the equation
\eqref{NDEabsODE}. As usual, we look at the lower-order case,
$k=1$. Once again, in order to remove the nonlinearity in the
highest (third)-order operator, the substitution
 $ 
Y = |f|^nf
 $ 
is used.  This yields the semilinear third-order  equation
\begin{equation}
 \label{Y33}
  f=|Y|^{-\frac n{n+1}}Y: \,\,\,
\mbox{$Y^{\prime\prime\prime} + \frac{p-(n+1)}{3(p-1)(n+1)}\,
|Y|^{-\frac{n}{n+1}}Y^\prime y  +
\frac{1}{p-1}\,|Y|^{-\frac{n}{n+1}}Y -
|Y|^{\frac{p-(n+1)}{n+1}}Y=0$}.
\end{equation}
Due to the complexity of the equation, which remains to be of the
third order, there is not much hope of solving this problem using
a shooting method, since, in fact, we do not know in detail the
``nonlinear bundle" as $ y \to +\iy$. As usual and as in the
semilinear case $n=0$, there is a difficulty in finding the
correct boundary points for $y>0$, such that the correct
oscillations are found.

 Hence, we return to the BVP
setting, trying to ``optimise" and ``minimise"  the oscillatory
bundle for $y \gg 1$. However, even using this approach, which was
rather  effectively implemented in the simpler semilinear case
$n=0$ in \cite[\S~6]{RayGI}, it is  difficult to produce reliable
numerics, due to the highly nonlinear/oscillatory nature of the
problem, even for smaller values of $p$ and $n$. In addition, we
must be careful when plotting, as we must avoid approaching any
nonlinear bifurcation points in $p$ given in \ef{pp1}, which
actually are not known explicitly.

\ssk

As a first example, in Figure \ref{FN1}, we present ``almost
converging" VSS profile for $p=5$ and small $n=0.1$. This example
is of a particular importance: its tail for $y \sim 35$ is clearly
not symmetric about $\{Y=0\}$ and is positive. According to
\cite[\S~6]{RayGI}, such a solution {\em is not} a VSS profile in
$\re$, and actually corresponds to some (obscure) BVP setting on a
bounded interval, which is of no interest here. We should avoid
such solutions in the future, even if these have been obtained
with a perfect convergence up to the tolerances $\sim 10^{-3}$.


 \begin{figure}
\centering
\includegraphics[scale=0.85]{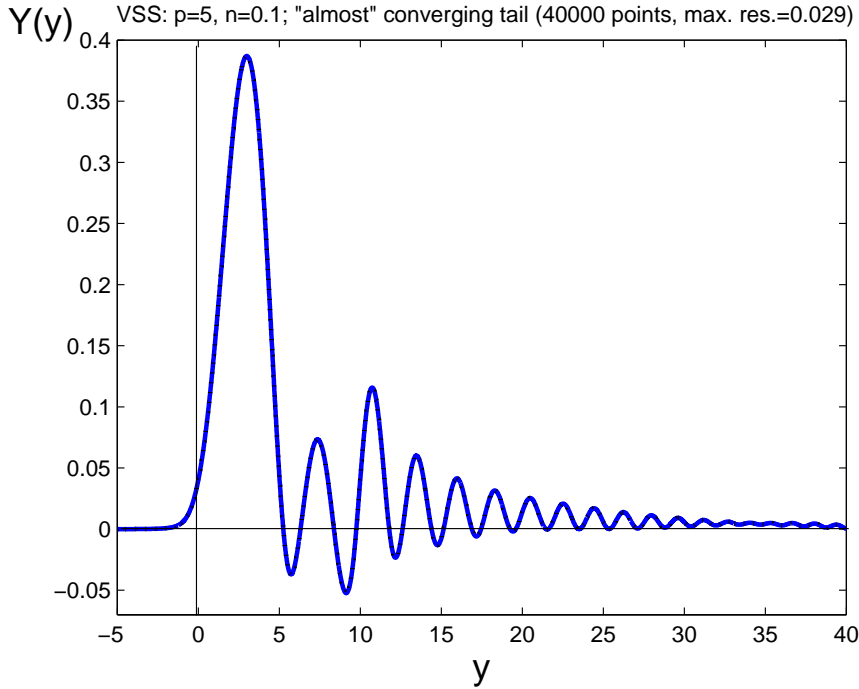}  
\vskip -.3cm
  \caption{A VSS similarity solution $Y(y)$ of the ODE \ef{Y33}
  for $p=5$ and $n=0.1$.}
 \label{FN1}
\end{figure}

The next Figure \ref{FN2} shows a better converging with a good
tail VSS profile with the same $p=5$ and a larger $n=0.6$. We see
that, for such larger $n$, the tail for $y \gg 1$ gets smaller and
keeps being symmetric.

For larger $n \ge 1$, the numerics gets less reliable, though we
have got a number of ``almost converging" results. For instance,
in Figure \ref{FN3}, this is done for $p=10$ and $n=1$. The tail
is now larger for $y \sim 12$, and remains rather symmetric.


 \begin{figure}
\centering
\includegraphics[scale=0.85]{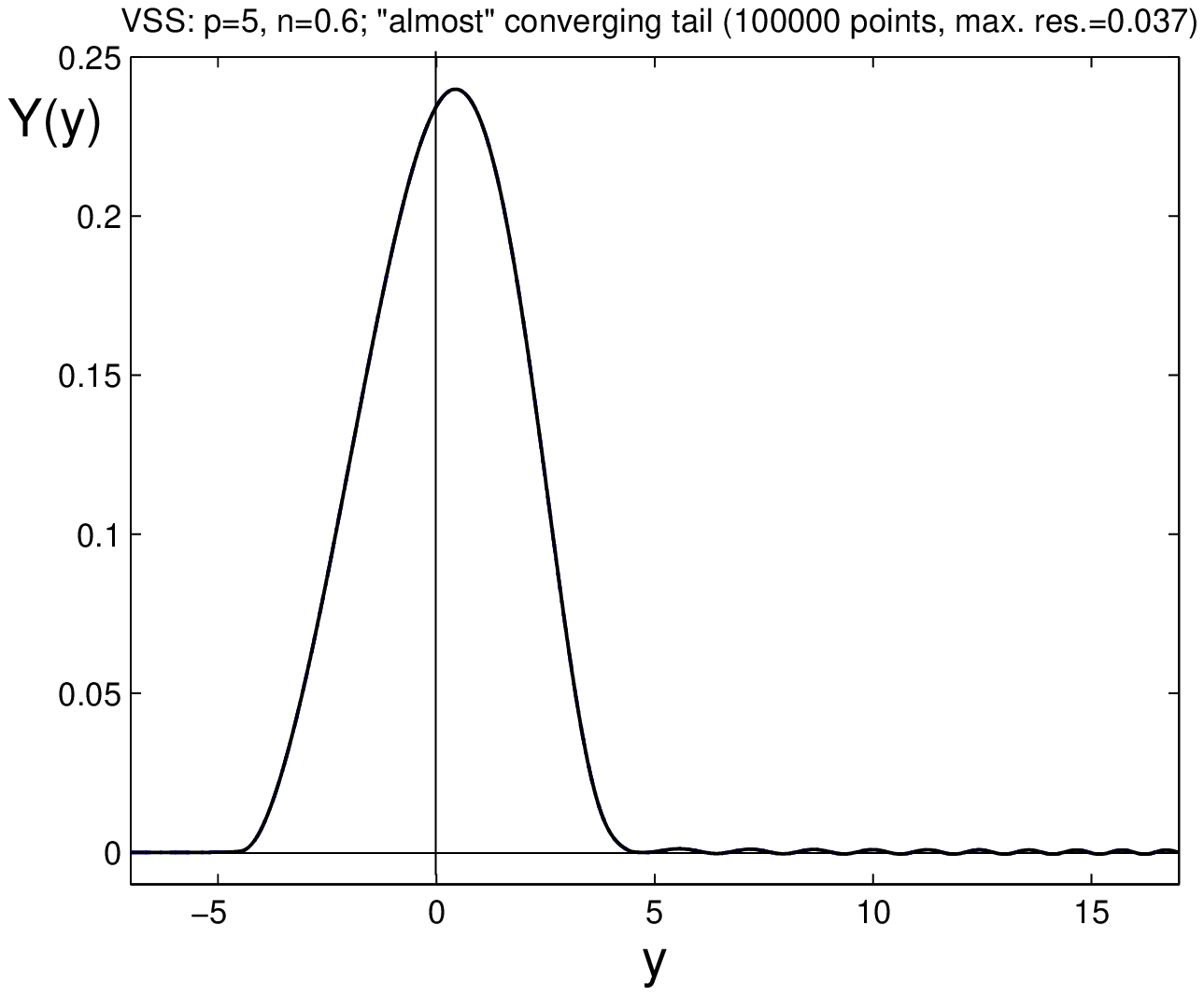}  
\vskip -.3cm
  \caption{A VSS similarity solution $Y(y)$ of the ODE \ef{Y33}
  for $p=5$ and $n=0.6$.}
 \label{FN2}
\end{figure}


 \begin{figure}
\centering
\includegraphics[scale=0.85]{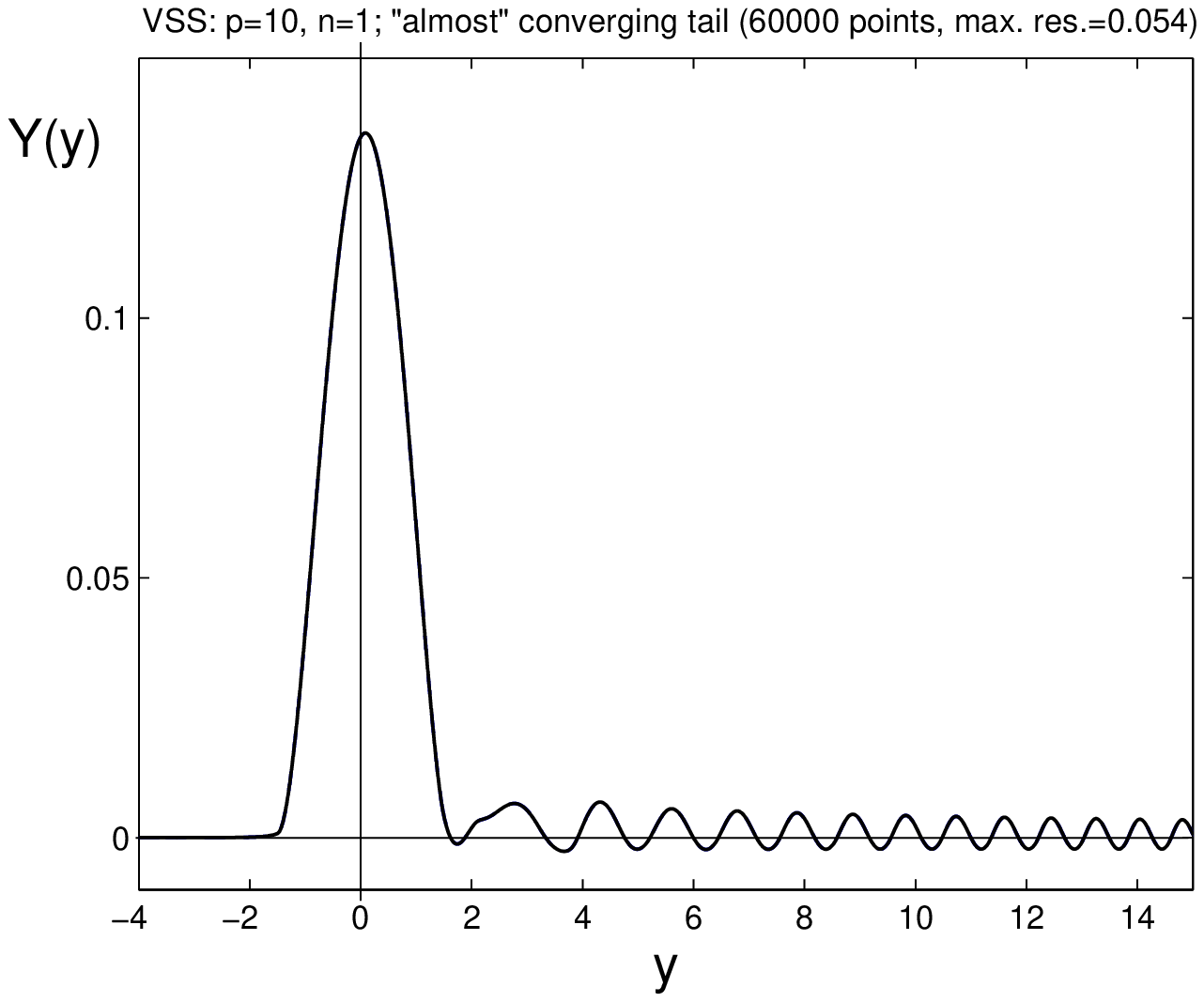}  
\vskip -.3cm
  \caption{A VSS similarity solution $Y(y)$ of the ODE \ef{Y33}
  for $p=10$ and $n=1$.}
 \label{FN3}
\end{figure}

Further increasing $n$ requires also increasing $p$. In the next
Figures \ref{FN4} and \ref{FN5}, we show VSS profiles for $p=12$,
$n=2$ and $p=16$, $n=3$, respectively. In the former one, Figure
\ref{FN4}, we present two different profiles (solid and dotted
lines), showing that non-converging tails do not affect the
convergence in the dominant positive part of the profiles.


 \begin{figure}
\centering
\includegraphics[scale=0.85]{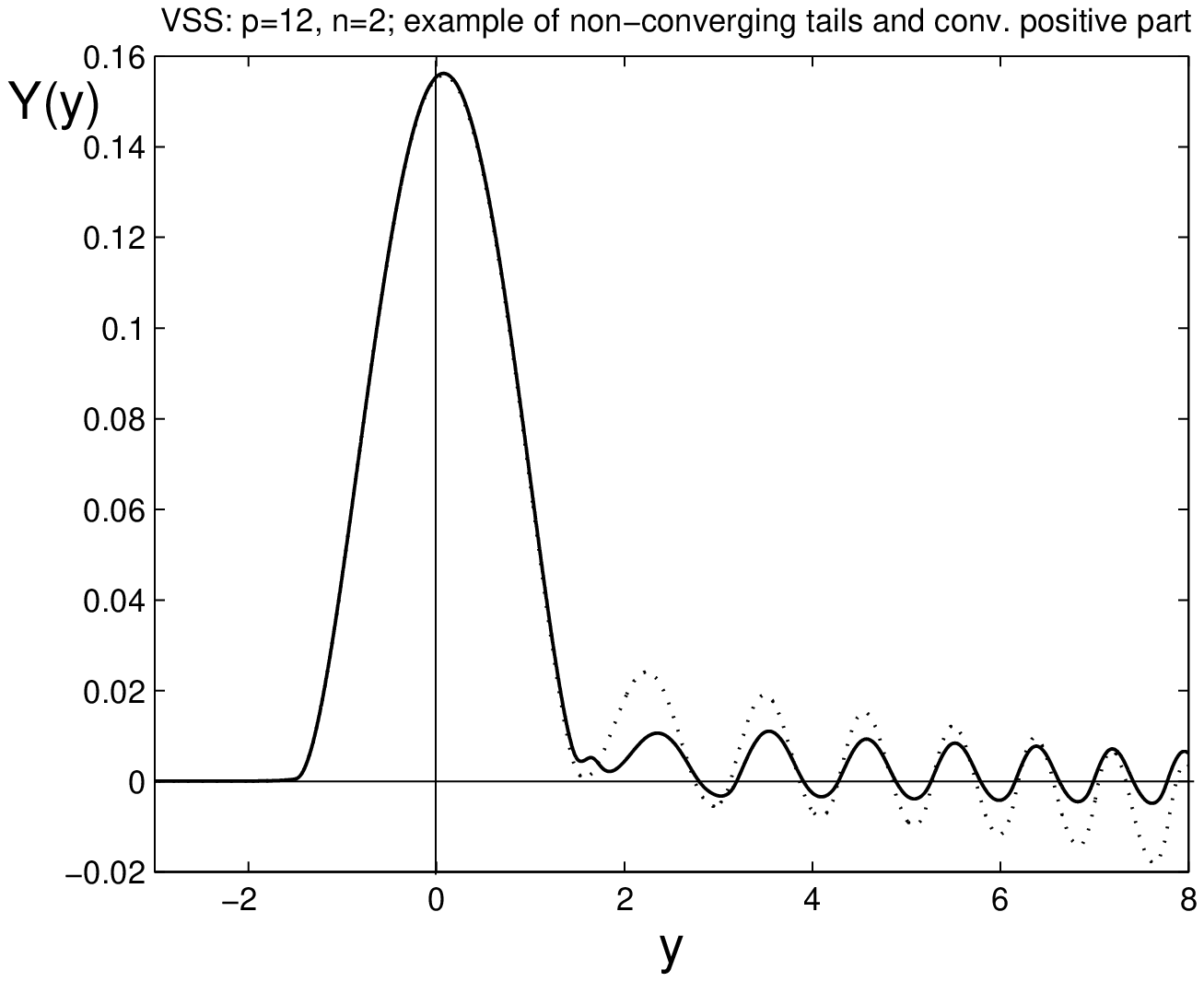}  
\vskip -.3cm
  \caption{A VSS similarity solution $Y(y)$ of the ODE \ef{Y33}
  for $p=12$ and $n=2$.}
 \label{FN4}
\end{figure}


 \begin{figure}
\centering
\includegraphics[scale=0.85]{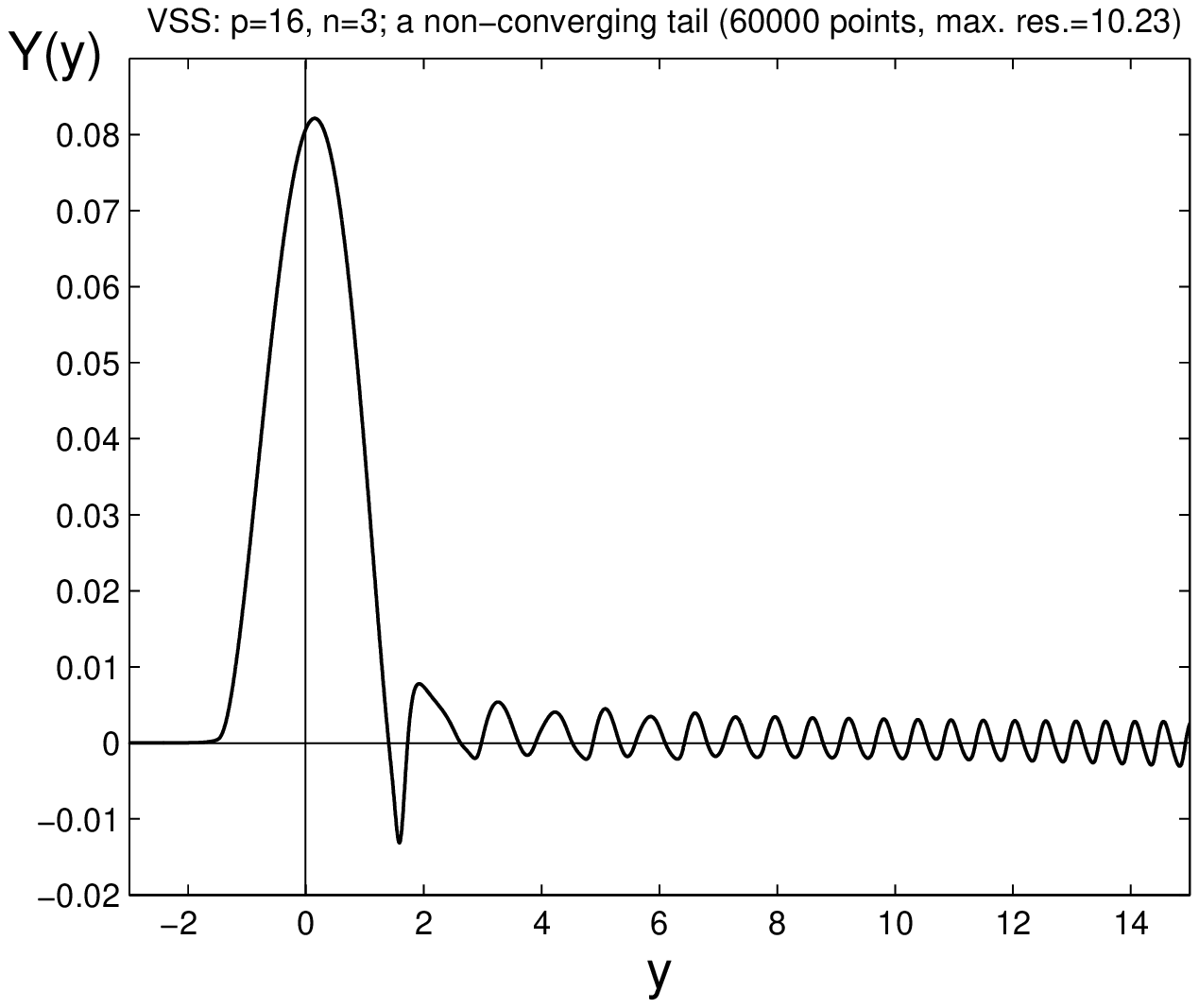}  
\vskip -.3cm
  \caption{A VSS similarity solution $Y(y)$ of the ODE \ef{Y33}
  for $p=16$ and $n=3$.}
 \label{FN5}
\end{figure}

Finally, we must admit that, since the correct ``oscillatory
bundle behaviour" as $y \to +\iy$ is still poorly understood,  we
cannot control it by choosing  proper ``minimal nonlinear
components" (a proper symmetry is not enough; see
\cite[\S~6.3]{RayGI}). Therefore, we must confess that, even
having good enough numerical convergence, we cannot guarantee that
the above numerical examples get into the countable family
indicated in \ef{Count1}. In other words, without a full use of
correct ``minimal oscillatory bundle" as $y \to +\iy$ (meaning
non-posing {\em any} condition at the singular end-point
$y=+\iy$), the family of VSS profiles becomes continuous. Then,
for any fixed value of $p>n+1$, ``solutions" $f(y)$ represent a
continuous  line, instead of a predicted at most countable
 subset of points, lying on the above
$p$-branches.

On the other hand, it is plausible that those Figures correctly
describe a general geometry of VSS profiles. In particular, we
have always observed a strong ``stability" of the first maximal
positive hump, which turned out to be rather independent of the
resulting tail for $y \gg 1$. In other words, we do note that,
whilst the oscillatory part is  difficult to obtain (and properly
justify mathematically), the non-oscillatory structure is rather
stable and does not change much, when changing the length in which
the problem is evaluated on. 
In particular,  the numerical profile in the last Figure \ref{FN5}
was obtained with the worst {\em maximal residual} $=10.23$
(actually meaning no convergence at all). However, we guarantee
that this non-convergence takes place in the tail {\em only},
while the dominant positive part remains stable and ``rather
convergent", so we do not hesitate to present such a Figure here.
However, in almost all other Figures, the convergence is much
better and is not worse than $2-5 \,\%$.  It may be also noted
that this confirms  such oscillatory (at least, symmetric) tails,
even when these are large in the amplitude,  are essentially {\em
zero} in some ``weak" sense (which we do not want to specify at
the moment).

\subsection{A nonlinear limit $n \to + \iy$: an example}

Finally, let us note that,  as in Section \ref{S4}, one way of
finding proper oscillatory patterns $f(y)$, would be to look at
the behaviour as $n\to +\infty$ and reduce the ODE to a simpler
one. We present an example of such a limit along the following
straight line on the $\{n,p\}$-plane:
  \be
  \label{pn21}
 p=n+4 \to + \iy \asA n \to + \iy.
  \ee
 Then the ODE \ef{Y33} reads
  \be
  \label{Y331}
   \tex{
  Y'''+ \frac 1{n+3}\, (|Y|^{-\frac n{n+1}}Y y)'-|Y|^{\frac 3{n+1}}Y=0,
  }
   \ee
   so that, on scaling, one gets
    \be
    \label{Z1}
    Y= (n+3)^{-\frac {n+1}n} \tilde Y
     \LongA
      \tilde Y'''+ (|\tilde Y|^{-\frac n{n+1}}\tilde Y y)'-
      (n+3)^{-\frac{3}n}\, |\tilde Y|^{\frac 3{n+1}}\tilde Y=0.
   \ee
Passing to the limit $n \to + \iy$, in the ODE in \ef{Z1}, in the
class of uniformly bounded solutions
 and using that $(n+3)^{-\frac{3}n} \to 1$,
 yield two terms as in
 the equation \ef{ninftyeqn} with an extra linear one:
 \be
 \label{Z11}
\tilde Y'''+ (|\tilde Y|^{-\frac n{n+1}}\tilde Y y)'-
      \tilde Y=0.
      \ee
However, since unlike \ef{ninftyeqn}, this is a third-order ODE,
 an algebraic treatment of the first nonlinear eigenfunction
$\tilde Y_0(y)$, as in Theorem \ref{Th.2}, becomes rather illusive.
Anyway, this shows a principal possibility to study the
``nonlinear" limits $n \to +\iy$.

On the other hand, scaling also the independent variable $y$,
 \be
 \label{Z12}
 Y= C \, \tilde Y, \quad y=a \, \tilde y, \andA C=(n+3)^{-\frac{n+1}n}
 a^{\frac{3(n+1)}n},
  \ee
  yields, instead of \ef{Z1},
   \be
   \label{Z13}
   \tilde Y'''+ (|\tilde Y|^{-\frac n{n+1}}\tilde Y \tilde y)'-
      (n+3)^{-\frac{3}n}\,a^{3+ \frac 9n} |\tilde Y|^{\frac 3{n+1}}\tilde Y=0.
   \ee
   Therefore, passing to the limit $n \to +\iy$ in \ef{Z13},
   after integrating once, we arrive precisely at the equation
\ef{ninftyeqn} provided that
 $$
 a=a(n) \to 0^+ \asA n \to +\iy.
 $$
Therefore, we obtain the same nonlinear eigenfunctions for
$l=0,1,2$ via the above algebraic-geometric approach, but,
according to \ef{Z12}, on expanding subsets in the independent
$\tilde y$-variable.

 Recall that the above
limit as $n \to +\iy$, with a possible study of branching as in
Section \ref{S4.4}, occurs along the straight line \ef{pn21} (or
in its ``small neighbourhood") in the 2D parametric
$\{n,p\}$-plane. Along other lines, the limits can be different
and lead to other patterns $\tilde Y_1(y)$, $\tilde Y_2(y)$, etc.

\ssk



Finally, overall, as we have seen, the ODE \ef{Y33} represents a
serious theoretical challenge with respect to both analytical
study ($n$-branching, $p$-bifurcation diagrams, and $p$-branches;
 the
latter are known for $n=0$ \cite[\S~7]{RayGI}, etc.), as well as
even a numerical one.


\addcontentsline{toc}{chapter}{Bibliography}

\bibliographystyle{amsplain}
\bibliography{biblioFG}

\providecommand{\bysame}{\leavevmode\hbox to3em{\hrulefill}\thinspace}
\providecommand{\MR}{\relax\ifhmode\unskip\space\fi MR }
\providecommand{\MRhref}[2]{%
  \href{http://www.ams.org/mathscinet-getitem?mr=#1}{#2}
}
\providecommand{\href}[2]{#2}
\begin{thebibliography}{10}

\bibitem{NonFuncAnaly}
K.~Deimling, \emph{\rm {N}onlinear {F}unctional {A}nalysis}, Springer-Verlag,
  Berlin/Heidelberg, 1985.

\bibitem{Gl4}
J.D. Evans, V.A. Galaktionov, and J.R. King, \emph{Source-type solutions of the
  fourth-order unstable thin film equation}, Euro J. Appl. Math. \textbf{18}
  (2007), 273--321.

\bibitem{Gl6}
\bysame, \emph{Unstable sixth-order thin film equation {II}. {G}lobal
  similarity patterns}, Nonlinearity \textbf{20} (2007), 1843--1881.

\bibitem{RayGI}
R.S. Fernandes and V.A. Galaktionov, \emph{Very singular similarity solutions
  and {H}ermitian spectral theory for semilinear odd-order {PDE}s}, Adv.
  Differ. Equat., submitted (arXiv:0910.4916).

\bibitem{Gal3NDENew}
V.A. Galaktionov, \emph{{O}n single point gradient blow-up and nonuniqueness
  for a third-order nonlinear dispersion equation}, Stud. Appl. Math., to
  appear {\rm (an earlier preprint in: ar{X}iv:0902.1635)}.

\bibitem{2mSturm}
\bysame, \emph{Sturmian nodal set analysis for higher-order parabolic equations
  and applications}, Adv. Differ. Equat. \textbf{12} (2007), 669--720.

\bibitem{GalPMEn}
\bysame, \emph{Countable branching of similarity solutions of higher-order
  porous medium type equations}, Adv. Differ. Equat. \textbf{13} (2008),
  641--680.

\bibitem{GPndeII}
\bysame, \emph{Nonlinear dispersion equations: smooth deformations, compactons,
  and extensions to higher orders}, Comput. Math. Math. Phys. \textbf{48}
  (2008), 1823--1856 (arXiv:0902.0275).

\bibitem{GalNDE5}
\bysame, \emph{Shock waves and compactons for fifth-order nonlinear dispersion
  equations}, Euro J. Appl. Math. \textbf{21} (2010), 1--50 (arXiv:0911.4446).

\bibitem{EvoCompNonlEigPME}
V.A. Galaktionov and P.J. Harwin, \emph{On evolution completeness of nonlinear
  eigenfunctions for the porous medium equation in the whole space}, Adv.
  Differ. Equat. \textbf{10} (2005), 635--674.

\bibitem{PetI}
\bysame, \emph{On centre subspace behaviour in thin film equations}, SIAM J.
  Appl. Math. \textbf{69} (2009), 1334--1358 (arXiv:0901.3995).

\bibitem{GMPSob}
V.A. Galaktionov, E.~Mitidieri, and S.I. Pohozaev, \emph{Variational approach
  to complicated similarity solutions of higher-order nonlinear evolution
  equations of parabolic, hyperbolic, and nonlinear dispersion types {\rm
  (ar}{\rm {x}}{\rm iv:0902.1425)}}, In: Sobolev Spaces in Mathematics. II,
  Appl. Anal. and Part. Differ. Equat., Series: Int. Math. Ser., p.~Vol. 9, V.
  Maz'ya Ed., Springer, New York, 2009.

\bibitem{GPnde}
V.A. Galaktionov and S.I. Pohozaev, \emph{Third-order nonlinear dispersive
  equations: shocks, rarefaction, and blow-up waves}, Comput. Math. Math. Phys.
  \textbf{48} (2008), 1784--1810 (arXiv:0902.0253).

\bibitem{GSVR}
V.A. Galaktionov and S.R. Svirshchevskii, \emph{\rm {E}xact {S}olutions and
  {I}nvariant {S}ubspaces of {N}onlinear {P}artial {D}ifferential {E}quations
  in {M}echanics and {P}hysics}, Chapman and Hall/CRC, Florida, 2007.

\bibitem{VSSSParaPDEs}
V.A. Galaktionov and J.F. Williams, \emph{On very singular similarity solutions
  of a higher-order semilinear parabolic equation}, Nonlinearity \textbf{17}
  (2004), 1075--1099.

\bibitem{Kaw85}
S.~Kawamoto, \emph{An exact transformation from the {H}arry {D}ym equation to
  the modified {K}d{V} equation}, J. Phys. Soc. Japan \textbf{54} (1985),
  2055--2056.

\bibitem{FuncAnaly}
A.N. Kolmogorov and S.V. Fomin, \emph{\rm {F}unctional {A}nalysis: {V}olume 1},
  4th ed., Graylock, Rochester, 1957.

\bibitem{GeoMethNonAn}
M.A. Krasnosel'ski\u{\i} and P.P. Zabre\u{\i}ko, \emph{\rm {G}eometrical
  {M}ethods of {N}onlinear {A}nalysis}, Springer-Verlag, Berlin/Tokyo, 1984.

\bibitem{NaimarkI}
M.A. Naimark, \emph{\rm {L}inear {D}ifferential {O}perators, {P}art {I}}, Ungar
  Publ. Comp., New York, 1967.

\bibitem{NaimarkII}
\bysame, \emph{\rm {L}inear {D}ifferential {O}perators, {P}art {II}}, Ungar
  Publ. Comp., New York, 1968.

\bibitem{Compactons}
P.~Rosenau and J.M. Hyman, \emph{Compactons: solitons with finite wavelength},
  Phys. Rev. Lett. \textbf{70} (1993), 564--567.

\bibitem{quasilin}
A.A. Samarskii, V.A. Galaktionov, S.P. Kurdyumov, and A.P. Mikhailov, \emph{\rm
  {B}low-up in {Q}uasilinear {P}arabolic {E}quations}, Walter de Gruyter,
  Berlin, 1995.

\bibitem{VainbergTr}
M.A. Vainberg and V.A. Trenogin, \emph{\rm {T}heory of {B}ranching of
  {S}olutions of {N}on-{L}inear {E}quations}, Noordhoff Int. Publ., Leiden,
  1974.

\bibitem{Zel56}
Ya.B. Zel'dovich, \emph{The motion of a gas under the action of a short term
  pressure shock}, Soviet Phys. Acoustics \textbf{2} (1956), 25--35.

\end{thebibliography}


\end{document}